\documentclass[journal]{IEEEtran}
\usepackage{indentfirst}
\usepackage{mathrsfs}
\usepackage{graphicx}
\usepackage{epstopdf}
\usepackage{amsmath, amssymb, amsthm, bm}
\usepackage{leftidx}
\usepackage{extarrows}
\usepackage{stfloats}
\usepackage{picinpar}
\usepackage{enumerate}
\usepackage{algorithm}
\usepackage{algorithmicx}
\usepackage{algpseudocode}
\usepackage{epsfig}
\usepackage{latexsym}
\usepackage{amsfonts}
\usepackage{enumerate}
\usepackage{graphics}
\usepackage{float}
\usepackage{pict2e}
\usepackage{tikz}
\usepackage{lineno}
\usepackage[pdftex,colorlinks,linkcolor=blue,anchorcolor=blue,citecolor=blue]{hyperref}
\usepackage{moreverb}
\usepackage{color}
\usepackage[nosort]{cite}
\usepackage{subfigure}
\usepackage{multirow}
\usetikzlibrary{arrows,automata}
\usepackage{colortbl}
\usepackage{wrapfig}

\definecolor{myGray}{rgb}{.95,.95,.95}

\newfont{\bb}{msbm10}

\newcommand{\tr}{^{\sf T}}
\newcommand{\M}[1]{{\bf{#1}}}

\newcommand{\m}[1]{{\mathrm{#1}}}

\allowdisplaybreaks[4]

\newtheorem{thm}{Theorem}
\newtheorem{ass}{Assumption}

\newtheorem{lem}{Lemma}

\newtheorem{rmk}{Remark}
\newtheorem{defn}{Definition}
\newtheorem{example}{Example}
\newtheorem{prop}{Proposition}

\ifCLASSINFOpdf

\else

\fi

\begin{document}
%
\title{DISA: A Dual Inexact Splitting Algorithm for Distributed Convex Composite Optimization}
%
%
%

\author{Luyao~Guo,
Xinli~Shi,~\IEEEmembership{Senior Member,~IEEE,}
Shaofu~Yang,~\IEEEmembership{Member,~IEEE,}
and Jinde~Cao,~\IEEEmembership{Fellow,~IEEE}
\thanks{This work was supported in part by the National Natural Science Foundation of China under Grant Nos. 61833005, 62003084 and 62176056, the Natural Science Foundation of Jiangsu Province of China under Grant No. BK20200355, and Young Elite Scientists Sponsorship Program by CAST, 2021QNRC001. (Corresponding author: Jinde Cao.)}

\thanks{Luyao Guo is with the School of Mathematics, Southeast University, Nanjing 210096, China
(e-mail: \href{mailto:ly_guo@seu.edu.cn}{{ly\_guo@seu.edu.cn; guo\_luyaoo@163.com}}).}
\thanks{Xinli Shi is with the School of Cyber Science \& Engineering, Southeast University, Nanjing 210096, China
(e-mail: \href{mailto:xinli_shi@seu.edu.cn}{{xinli\_shi@seu.edu.cn}}).}
\thanks{Shaofu Yang is with the School of Computer Science and Engineering, Southeast University, Nanjing 210096, China
(e-mail: \href{mailto:sfyang@seu.edu.cn}{{sfyang@seu.edu.cn}}).}
\thanks{Jinde Cao is with the School of Mathematics, Southeast University, Nanjing 210096, China, with Purple Mountain Laboratories, Nanjing 211111, China, and also with Yonsei Frontier Lab, Yonsei University, Seoul 03722, South Korea
(e-mail: \href{mailto:jdcao@seu.edu.cn}{{jdcao@seu.edu.cn}}).}
}

\maketitle

\begin{abstract}
In this paper, we propose a novel Dual Inexact Splitting Algorithm (DISA) for distributed convex composite optimization problems, where the local loss function consists of a smooth term and a possibly nonsmooth term composed with a linear mapping. DISA, for the first time, eliminates the dependence of the convergent step-size range on the Euclidean norm of the linear mapping, while inheriting the advantages of the classic Primal-Dual Proximal Splitting Algorithm (PD-PSA): simple structure and easy implementation. This indicates that DISA can be executed without prior knowledge of the norm, and tiny step-sizes can be avoided when the norm is large. Additionally, we prove sublinear and linear convergence rates of DISA under general convexity and metric subregularity, respectively. Moreover, we provide a variant of DISA with approximate proximal mapping and prove its global convergence and sublinear convergence rate. Numerical experiments corroborate our theoretical analyses and demonstrate a significant acceleration of DISA compared to existing PD-PSAs.
\end{abstract}


%
\IEEEpeerreviewmaketitle

\section{Introduction}
In this paper, we focus on the distributed convex Composite Optimization Problem (COP) over networks with $m$ agents:
\begin{align}\label{Problem}
J^*\triangleq\mathop{\min}\limits_{\m{x}\in \mathbb{R}^n}~\sum_{i=1}^m \big(f_i(\m{x})+g_i(U_i\m{x})\big),
\end{align}
where $U_i\in\mathbb{R}^{p\times n}$ is a matrix, and the local loss functions $f_i:\mathbb{R}^n\rightarrow(-\infty,+\infty)$, $g_i:\mathbb{R}^p\rightarrow(-\infty,+\infty]$ are convex and accessed only by agent $i$. Suppose that $f_i$ is $L_i$-smooth, $g_i$ is proper and closed (possibly nonsmooth), and each agent permits to local computation and communication with its immediate neighbors to obtain a consensus solution $\m{x}^*$ to problem \eqref{Problem}. Such composite structure covers a wide range of optimization scenarios \cite{Nedic20151}. Here, we give two examples.
\begin{example}
Support vector machine: $\min_{\m{x}} \frac{1}{m} \sum_{i=1}^{m}(\|\m{x}\|^2+ \max(0,1-\mathcal{B}_i\mathcal{A}_i\tr\m{x})$,
where $\{(\mathcal{A}_i,\mathcal{B}_i)\}^{m}_{i=1}$ are the samples in $n$-dimensional space, and $\mathcal{B}_i\in\{-1,1\}$ is the label of the $i$-th sample. Here, $f_i(\m{x})=\frac{1}{m}\|\m{x}\|^2$, $g_i(\cdot)=\max\{0,1-\cdot\}$, and $U_i=\mathcal{B}_i\mathcal{A}_i\tr$ is a feature matrix.
\end{example}
\begin{example}
Constrained Logistic regression problems:
$
\min _{\m{x}}\sum_{i=1}^{m}\big\{\frac{1}{m_i}\sum_{j=1}^{m_{i}} \ln \big(1+e^{-(\mathcal{U}_{i j}^{\top} \m{x})\mathcal{V}_{i j}}\big)+\frac{1}{2}\|\m{x}\|^2\big\}+\delta_{\mathcal{C}_i}(U_i\m{x}),
$
where $\mathcal{C}_i=\{\m{y}_i:\m{y}_i\leq w_i\}$, and any agent $i$ holds its own training date $\left(\mathcal{U}_{i j}, \mathcal{V}_{i j}\right) \in$ $\mathbb{R}^{n} \times\{-1,1\}, j=1, \cdots, m_{i}$, including sample vectors $\mathcal{U}_{i j}$ and corresponding classes $\mathcal{V}_{i j}$.
\end{example}
Additionally, several common regularizers in machine learning, such as the general LASSO regularizer \cite{Example3}, the fused LASSO regularizer \cite{Tibshirani2005}, the octagonal selection and clustering algorithm for regression (OSCAR) regularizer \cite{Howard2008}, and the group LASSO regularizer \cite{Jacob2009} can be formulated as $g_i(U_i\m{x})$. Moreover, several practical control problems provided in \cite{Aybat2016} and \cite{Aybat2022} also can be covered by problem \eqref{Problem}.

Problem \eqref{Problem} has been the subject of numerous distributed optimization algorithms, which have been proposed to address its special cases. For instance, when $g_i\equiv0$, several studies have investigated distributed gradient descent and its extensions, such as \cite{Nedic2009,Jako2014,Yuan2016}. These primal methods are unable to achieve the exact convergence under a fixed step-size, resulting in a convergence gap. To tackle this issue, many primal-dual algorithms have been proposed in recent literature \cite{DADMM,EXTRA,DIGing,Harnessing,PushPull,Xu2018,ExactDiffusion,Liang2019,Pan2022}. When $g_i\neq 0$ with $U_i=I$ and $g_i$ is proximable, i.e., the proximal mapping of $g_i$ has an analytical solution or can be computed efficiently, several primal-dual algorithms based on proximal gradient have been provided \cite{Sulaiman2021,Xu2021,Chang2015,Aybat2017,PGEXTR,NDIS,AMM2022}. However, for the general problem \eqref{Problem}, these primal-dual algorithms \cite{Sulaiman2021,Xu2021,Chang2015,Aybat2017,PGEXTR,NDIS,AMM2022} may not be suitable. The main challenge arises from the fact that $(g_i\circ U_i)$ may not be proximable, even though $g_i$ is proximable, which results in a high computation cost of the proximal operator of $(g_i\circ U_i)$. To address this challenge, several primal-dual algorithms have been developed, including the recent TPUS \cite{TPUS}, TriPD-Dist \cite{TriPD}, and PDFP-Dist \cite{PDFP-DIS}. Additionally, two primal-dual methods have been proposed in \cite{Aybat2016} for distributed COPs over private affine conic constraint sets. Furthermore, when the conic constraint sets are defined by nonlinear functions (covering COP \eqref{Problem}), DPDA and DPDA-TV have been developed in \cite{Aybat2022} for static and time-varying networks, respectively.

Note that for the COP of the form $\min_{\m{x}} \{f(\m{x})+g(U\m{x})\}$, where the convex function $f$ is smooth and $g$ is proximable, Primal-Dual Proximal Splitting Algorithms (PD-PSAs) \cite{Duality1,PDS2021} are generally effective. However, current research has primarily focused on centralized form \cite{TriPD,CV1,CV2,PDFP2O,PAPC,PDFP,AFBA,PD3O,DYS-T2}, and the choice of step-size is highly dependent on $\|U\tr U\|$. Specifically, when $\|U\tr U\|$ is large, the primal step-size $\tau$ and the dual step-size $\beta$ are forced to be small to guarantee the convergent step-size conditions. This requires more iterations and thus reduces the efficiency of these algorithms. In addition, since distributed algorithms TPUS, TriPD-Dist, and PDPF-Dist have been designed by PD3O \cite{PD3O}, TriPD \cite{TriPD}, and PDFP \cite{PDFP}, respectively, the choice of step-size also depends on $\|U_iU_i\tr\|$.

To overcome this dependency and solve problem \eqref{Problem} in a distributed manner, we propose a novel distributed PD-PSA called Dual Inexact Splitting Algorithm (DISA), with the convergent step-size condition
\begin{align}\label{Scond1}
0<\tau_i<{2}/{L_i},~0<\tau_i\beta<1,i=1,\cdots,m.
\end{align}
Compared to existing PD-PSAs \cite{TPUS,TriPD,CV1,CV2,PDFP2O,PAPC,PDFP,AFBA,PD3O}, the convergent step-size range of DISA is not restricted by any condition related to the Euclidean norm of the linear operator explicitly or implicitly. As a result, it can avoid tiny step-sizes of primal and dual updates even when $\|U_i U_i\tr\|$ is large thereby ensuring fast convergence, which is exactly the distinctive advantage of DISA. Additionally, there is no need to estimate $\|U_iU_i\tr\|$ when tuning the step-sizes. We consider DISA as a necessary and important supplement to PD-PSAs.

For the algorithm development, when $g_i$ is proximable, the inexact Forward-Backward Splitting method (FBS) \cite{FBS} and Fenchel-Moreau-Rockafellar (FMR) duality \cite[Chapter 19]{Duality1}, \cite{Duality2} have been effectively utilized. Firstly, problem \eqref{Problem} is formulated as a COP consisting of a smooth term and two nonsmooth terms. Then, since the proximal operator of the sum of two nonsmooth terms may not have a closed-form representation, we provide an efficient analytical estimation for the proximal operator via FMR duality and Moreau-Yosida regularization. This processing ensures that the convergent step-size range is independent of the network topology. Additionally, inspired by a key technique in the development of primal-dual algorithms, \emph{preconditioning}, or
equivalently changing the metric of the ambient Euclidean space \cite{PDS2021}, \cite[Sections 2 and 3]{BookYin2022}, we specifically design a non-diagonal preconditioner for the dual update of DISA to eliminate the dependence of the convergent step-sizes range on $\|U_i U_i\tr\|$. By this technique, the dual update for agent $i$ is equivalent to solving a positive definite system of linear equations. When $g_i$ is not proximable, to reduce the computational burden, we provide a variant of DISA (V-DISA), where it only needs to approximately solve a sequence of proximal mappings.

For the convergence analysis, with general convexity and condition \eqref{Scond1}, we prove the convergence of DISA, and establish an $O({1}/{k})$ non-ergodic convergence rate. In addition, under a stronger step-size condition that $0<\tau<{1}/{L_i}$ and $0<\tau_i\beta<1,i=1,\cdots,m$, we establish $O({1}/{k})$ ergodic convergence rates in the primal-dual gap, primal suboptimality, and consensus violation, respectively. With metric subregularity which is weaker than the strong convexity, we establish the linear convergence rate of DISA under the condition \eqref{Scond1}. Moreover, the global convergence and sublinear convergence rate of V-DISA are established under the same step-size condition as DISA and a summable absolute error criterion.

This paper is organized as follows. In Section \ref{Sec2}, we cast problem \eqref{Problem} into a constrained form, and based on the reformulation, DISA is provided with the help of inexact FBS and FMR duality, and we compare DISA with existing distributed optimization algorithms under various situations. Then, the convergence properties under general convexity and metric subregularity are investigated in Section \ref{SEC-Convergence Analysis}, respectively. Moreover, in Section \ref{SECIN}, a variant of DISA is provided, and the convergence analysis is given. Finally, two numerical simulations are implemented in Section \ref{Simulation} and conclusions are given in Section \ref{Conclusion}.

{Notations and preliminaries}: $\mathbb{R}^n$ denotes the $n$-dimensional vector space with inner-product $\langle \cdot,\cdot\rangle$. The $1$-norm, Euclidean norm and infinity norm are denoted as $\|\cdot\|_1$, $\|\cdot\|$ and $\|\cdot\|_{\infty}$, respectively. $0$ and $I$ denote the null matrix and the identity matrix of appropriate dimensions, respectively. Denote $1_n\in \mathbb{R}^n$ as a vector with each component being one. For a matrix $A\in \mathbb{R}^{p\times n}$, $\|A\|:=\max_{\|\m{x}\|=1} \|A\m{x}\|$. If $A$ is symmetric, $A\succ 0$ means that $A$ is positive definite. For a given matrix $\mathrm{H} \succ0$, non-empty closed convex subset $\mathcal{X}$ and vector $\m{x}$ in the same space, $\mathrm{H}$-norm is defined as $\|\m{x}\|^2_{\mathrm{H}}=\langle \m{x},\mathrm{H} \m{x}\rangle$, $\mathrm{dist}_{\mathrm{H}}(\m{x},\mathcal{X}):=\mathrm{inf}\{\|\m{x}-\m{x}'\|_{\mathrm{H}}:\m{x}'\in\mathcal{X}\}$ and $\mathcal{P}_{\mathcal{X}}^{\mathrm{H}}(\m{x}):=\arg\min_{\m{x}'\in \mathcal{X}}\|\m{x}-\m{x}'\|_{\mathrm{H}}$. Moreover, $\mathcal{B}_{r}(\bar{\m{x}}):=\{\m{x}:\|\m{x}-\bar{\m{x}}\|<r\}$ denotes the open Euclidean norm ball around $\bar{\m{x}}$ with radius $r>0$.

Let $\theta:\mathbb{R}^n\rightarrow (-\infty,+\infty]$. {Let $\m{H}\succ0$. The proximal operator of $ \theta$ relative to $\|\cdot\|_{\m{H}}$ is defined as $\mathrm{prox}^{\m{H}}_{ \theta}(\m{v})=\arg\min_{\m{x}}\{\theta(\m{x})+\frac{1}{2 \tau}\|\m{v}-\m{x}\|^{2}_{\m{H}}\}$, and the Moreau-Yosida regularization of $\theta$ is defined as $\widetilde{\theta}_{\m{H}}(\m{v})=\min_{\m{x}}\{\theta(\m{x})+\frac{1}{2}\|\m{v}-\m{x}\|^2_{\m{H}}\}$. If the proximal operator of $\theta$ has an analytical solution or can be computed efficiently, we say that $\theta$ is proximable. From \cite{Proximal}, the proximal operator of $\theta$ is unique, if $\theta$ is a proper closed convex function. Moreover, the Moreau-Yosida regularization of $\theta$ is a continuously differentiable convex function even $\theta$ is not, and the gradient is $\nabla \widetilde{\theta}_{\m{H}}(\m{x})=\m{H}(\m{x}-\mathrm{prox}^{\m{H}}_{ \theta}(\m{x}))$.}

\section{Dual Inexact Splitting Algorithm}\label{Sec2}
\subsection{Problem Statement}
Consider an undirected and connected network $\mathcal{G}(\mathcal{V},\mathcal{E})$, where $\mathcal{V}$ denotes the vertex set, and the edge set $\mathcal{E}\subseteq \mathcal{V}\times \mathcal{V}$ specifies the connectivity in the network, i.e., a communication link between agents $i$ and $j$ exists if and only if $(i,j)\in \mathcal{E}$. Let $\mathcal{N}_i=\{j\in\mathcal{N}:(i,j)\in\mathcal{E}\}\cup \{i\}$ denote the set of neighbours agent $i$ including itself. We introduce local copy $\mathrm{x}_{1,i}\in \mathbb{R}^n$ of $\m{x}$, which is the decision variable held by the agent $i$, and introduce the global mixing matrix $W$ associated with $\mathcal{G}$, where $W=[W_{ij}]\in \mathbb{R}^{m\times m}$. The mixing matrix $W$ is symmetric and doubly stochastic, and $\mathrm{null}(I-W)=\mathrm{Span}(1_m)$. Moreover, if $(i, j) \notin \mathcal{E}$ and $i \neq j$, $W_{i j}=W_{j i}=0$; otherwise, $W_{ij}>0$. By Perron-Frobenius theorem \cite{PFT}, we know that $W$ has a simple eigenvalue one and all the other eigenvalues lie in $(-1,1)$.
\begin{prop}
Let $\M{x}=[\m{x}_1\tr,\m{x}_2\tr]\tr\in\mathbb{R}^{m(n+p)}$, where $\m{x}_1=[\m{x}_{1,1}\tr,\cdots,\m{x}_{1,m}\tr]\tr $ and $\m{x}_2=[\m{x}_{2,1}\tr,\cdots,\m{x}_{2,m}\tr]\tr $, and $\M{U}=\mathrm{diag}\{U_1,\cdots,U_m\}$. Problem \eqref{Problem} is equivalent
to
\begin{equation}\label{Trans-P1}
\min_{\M{x}} \underbrace{ \overbrace{\sum_{i=1}^{m} f_i(\m{x}_{1,i})}^{:=F(\M{x})} + \overbrace{\sum_{i=1}^{m} g_i(\m{x}_{2,i})}^{:=G(\M{x})}}_{:=\Phi(\M{x})}+\delta_0(\underbrace{\left(
                                                                       \begin{array}{cc}
                                                                         \sqrt{\M{V}} & 0 \\
                                                                         \M{U} & -I \\
                                                                       \end{array}
                                                                     \right)}_{:=\M{B}}\M{x}),
\end{equation}
where $\M{V}=\frac{1}{2}(I-W)\otimes I_n$, and $\delta_0(\cdot)$ is an indicator function defined as $\delta_0(\M{B}\M{x})=0$ if $\M{B}\M{x}=0$; otherwise $\delta_0(\M{B}\M{x})=\infty$.
\end{prop}
\begin{IEEEproof}
Since $\mathrm{null}(I-W)=\mathrm{Span}(1_m)$, one has
$
\mathrm{x}_{1,1}=\cdots=\mathrm{x}_{1,m}\Leftrightarrow \M{V} \m{x}_1=0\Leftrightarrow \sqrt{\M{V}} \m{x}_1=0
$.
If $\M{x}^*=(\m{x}_1^*,\m{x}_2^*)$ solves problem \eqref{Trans-P1}, by the definition of $\delta_0(\M{B}\M{x})$, one has $\M{B}\M{x}^*=0$, i.e.,
$
\m{x}_1^*=\m{1}_m\otimes\m{x}^*$, $\mathrm{x}^*_{2,i}=U_i\m{x}^*, i\in\mathcal{V}
$, and
$\min_{\M{x}}( F(\M{x})+G(\M{x})+\delta_0(\M{Bx}))=\sum_{i=1}^{m}( f_i(\m{x}^*)+ g_i(U_i\m{x}^*))=J^*$. Therefore, $\m{x}^*$ solves problem \eqref{Problem}.
\end{IEEEproof}

Recall the Lagrangian of problem \eqref{Trans-P1}:
$\mathcal{L}(\mathbf{x},\M{y})=F(\M{x})+G(\mathbf{x})+\langle \M{y},{\M{B}}\M{x}\rangle,$
where $\M{y}\in \mathbb{R}^{m(n+p)}$ is the Lagrange multiplier. If the saddle-point problem exists a saddle-point, the set of saddle-points $\mathcal{M}^*$ can be characterized by $\mathcal{M}^*=\mathcal{X}^* \times \mathcal{Y}^*$ with
$\mathcal{X}^* =\arg\min_{\M{x}}\{F(\M{x})+G(\M{x})+\delta_0(\M{Bx})\}$ and
$\mathcal{Y}^*=\arg\min_{\M{y}}\{(F+G)^*(-\M{B}\tr\M{y})\}$,
where $\mathcal{X}^*$ is the optimal solution set to problem \eqref{Trans-P1}, $\mathcal{Y}^*$ is the optimal solution set to its dual problem, and $(F+G)^*$ is the Fenchel conjugate of $F+G$. Denote the KKT mapping as
$$
\mathcal{T}(\M{x},\M{y})=\left(\begin{array}{c}
                                             \partial G(\mathbf{x})+\nabla F(\mathbf{x})+{\M{B}}\tr\M{y} \\
                                            -{\M{B}}\mathbf{x}
                                          \end{array}
\right).
$$
It holds that $\left(\mathbf{x}^{*}, \M{y}^{*}\right) \in \mathcal{X}^* \times \mathcal{Y}^*$ if and only if $0 \in \mathcal{T}\left(\mathbf{\mathbf{x}^*,\M{y}^*}\right)$. Then, we give the following assumption.
\begin{ass}\label{ass1}
There exists a point $(\M{x}^*,\M{y}^*)\in \mathcal{M}^*$. The local loss functions $f_i:\mathbb{R}^n\rightarrow(-\infty,+\infty)$ and $g_i:\mathbb{R}^p\rightarrow(-\infty,+\infty]$ are convex. Moreover, $g_i$ is proper and closed (possibly nonsmooth), $f_i$ is differentiable and $\nabla f_i$ is $L_i$-Lipschitz continuous.
\end{ass}
\begin{lem}\label{NEWLEMMA1}
With Assumption \ref{ass1}, for any $\M{x},\M{y},\M{z}\in\mathbb{R}^{m(n+p)}$, we have the following two commonly used inequalities
\begin{align}
&\big\langle \M{x}-\M{y},\nabla F(\M{z})-\nabla F(\M{x})\big\rangle\leq \frac{1}{4} \|\M{y}-\M{z}\|_{\M{L}_F}^2,\label{SM1}\\
&\langle \M{x}-\M{y},\nabla F(\M{z})\rangle\leq F(\M{x})-F(\M{y})+\frac{1}{2}\|\M{y}-\M{z}\|_{\M{L}_F}^2,\label{SM2}
\end{align}
where $\M{L}_F=\mathrm{diag}\{L_1{I}_{n+p},\cdots,L_m{I}_{n+p}\}$.
\end{lem}
\begin{IEEEproof}
See Appendix \ref{APP45}.
\end{IEEEproof}

\subsection{Algorithm Development}
Considering problem \eqref{Trans-P1}, since $F(\M{x})$ is smooth and $G(\M{x})+\delta_{0}(\mathbf{B}\M{x})$ is nonsmooth, we can derive an algorithm based on FBS as follows
\begin{equation}\label{Al1}
\M{x}^{k+1}=\mathrm{prox}^{\Gamma^{-1}}_{(G+\delta_0\circ \M{B})}\big(\M{x}^k- \Gamma \nabla F(\M{x}^k)\big),
\end{equation}
where $\Gamma=\mathrm{diag}\{\Gamma_1,\Gamma_2\}$ is the step-size matrix with $\Gamma_1=\mathrm{diag}\{\tau_1I_n,\cdots,\tau_mI_n\}$ and $\Gamma_2=\mathrm{diag}\{\tau_1I_p,\cdots,\tau_mI_p\}$. By \cite[Corollary 28.9, pp. 522]{Duality1}, the algorithm is convergent when $0<\tau_i<{2}/{L_i}$. Unfortunately, the proximal operator of $G+\delta_0\circ \M{B}$ has no analytical solution, leading to a nontrivial derivation of $\mathrm{prox}^{\Gamma^{-1}}_{(G+\delta_0\circ \M{B})}(\cdot)$. To address this issue and inherit the favorable feature that the acceptable range of $\tau_i$ is network-independent, we will give an effective analytical estimation for $\mathrm{prox}^{\Gamma^{-1}}_{(G+\delta_0\circ \M{B})}(\cdot)$ based on FMR duality.

Recall the definition of $\mathrm{prox}^{\Gamma^{-1}}_{ (G+\delta_0\circ \M{B})}\big(\M{x}^k-\Gamma \nabla F(\M{x}^k)\big)$, which is equivalent to solving the following problem
\begin{equation}\label{Primal1}
\min_{\M{x}} \left\{ \begin{array}{c}
                       \frac{1}{2}\|\M{x}-(\M{x}^k-\Gamma\nabla F(\M{x}^k))\|_{\Gamma^{-1}}^2 \\
                       +G(\M{x})+\delta_0(\M{B}\M{x})
                     \end{array}
  \right\}.
\end{equation}
By \cite{Duality2}, the FMR dual of problem \eqref{Primal1} is
\begin{equation}\label{Dual1}
\max_{\M{y}'}  \left\{\begin{array}{l}-H(\M{y}')=
                       -\frac{1}{2}\|\M{B}\tr\M{y}'-\Gamma^{-1}\M{u}^k\|_{\Gamma}^2\\
                      \quad\quad\quad\quad\quad +\widetilde{G}_{\Gamma}(\M{u}^k-\Gamma\M{B}\tr\M{y}')
                    \end{array}
\right\},
\end{equation}
where $\M{y}'=\Gamma^{-1}\M{y}=\Gamma^{-1}[\m{y}_1\tr,\m{y}_2\tr]\tr$, $\M{u}^k=\M{x}^k-\Gamma\nabla F(\M{x}^k)$ and $\widetilde{G}_{\Gamma^{-1}}$ is the Moreau-Yosida regularization of $G$. Let $\M{y}_{opt}$ be solution to problem \eqref{Dual1}. By \cite[Proposition 3.4]{Duality2}, we have
\begin{align}\label{Al2}
\M{x}^{k+1}=\mathrm{prox}^{\Gamma^{-1}}_{G}(\M{u}^k-\Gamma \M{B}\tr\M{y}_{opt}).
\end{align}
Note that $\widetilde{G}_{\Gamma^{-1}}$, the Moreau-Yosida regularization of $G$, is continuously differentiable, and $\nabla \widetilde{G}_{\Gamma^{-1}}(\M{x})=\Gamma^{-1}(\M{x}-\mathrm{prox}^{\Gamma^{-1}}_{G}(\M{x}))$. One has $\nabla H(\M{y}')=-\Gamma^{-1}{\M{B}} \mathrm{prox}^{\Gamma^{-1}}_{ G}(\M{u}^k-\Gamma {\M{B}}\tr \M{y}')$.
On the other hand, by nonexpansivity of proximal operators,
$\|\nabla H(\M{y}_1)-\nabla H(\M{y}_2)\|\leq\|\M{BB}\tr\|  \|\M{y}_1-\M{y}_2\|, \M{y}_1,\M{y}_2\in \mathbb{R}^{m(n+p)}$.
Hence, $H(\M{y}')$ is $\|\mathbf{BB}\tr\|$-smooth. It implies that we can obtain the update of $\M{x}^{k+1}$ by solving the dual problem \eqref{Dual1} which is simpler than the primal problem \eqref{Primal1}. However, it is not practical to get an exact solution of \eqref{Dual1} as the coupling term $\frac{1}{2}\|\M{B}\tr\M{y}'-\Gamma^{-1}\M{u}^k\|^2$ requires infinite loops of consensus. To address this challenge, by the smoothness of $H(\M{y}')$, we give the approximation of $H(\M{y}')$ at $\M{y}^k$:
$
H(\M{y}')\approx H(\M{y}^{k})+\langle \nabla H(\M{y}^{k}),\M{y}'-\M{y}^{k}\rangle+\frac{1}{2}\|\M{y}'-\M{y}^k\|_{\Gamma^{-1}\M{Q}}^2
$,
where
$$
\M{Q}=
\left(
  \begin{array}{cc}
    \frac{1}{\beta}I & 0 \\
    0 & \M{S} \\
  \end{array}
\right),\M{S}=2\Gamma_1+\M{U}\Gamma_1\M{U}\tr+\frac{\beta}{1-\tau\beta}\M{U}\Gamma_1^2\M{U}\tr,
$$
$\beta>0$ is the dual step-size, and $\tau=\max_i\{\tau_i\}$. Based on the approximation, an inexact solution to problem \eqref{Dual1} is obtained by
$
\M{y}_{opt} \approx \mathop{\arg\min}_{\M{y}}\{\langle \nabla H(\M{y}^{k}),\M{y}-\M{y}^{k}\rangle+\frac{1}{2}\|\M{y}-\M{y}^k\|_{\Gamma^{-1}\M{Q}}^2\}
$.
Then, we have the updates of DISA
\begin{subequations}\label{DISAUP}
\begin{align}
\bar{\M{x}}^{k+1}&=\mathrm{prox}^{\Gamma^{-1}}_{ G}(\M{x}^k-\Gamma\nabla F(\M{x}^k)-\Gamma {\M{B}}\tr \M{y}^k),\label{barPrimalUP}\\
\M{y}^{k+1}&=\mathop{\arg\min}\limits_{\M{y}\in \mathbb{R}^{m(n+p)}}\big\{\frac{1}{2}\|\M{y}-\M{y}^k\|^2_{\M{Q}}-\langle\M{B}\tr\M{y},\bar{\M{x}}^{k+1}\rangle\big\}, \label{DualUP}\\
\M{x}^{k+1}&=\mathrm{prox}^{\Gamma^{-1}}_{ G}(\M{x}^k-\Gamma\nabla F(\M{x}^k)-\Gamma {\M{B}}\tr \M{y}^{k+1})\label{PrimalUP}.
\end{align}
\end{subequations}
Let $(\M{x}_{\m{fix}},\M{y}_{\m{fix}})$ be a fixed point of DISA \eqref{DISAUP}. It can be verified that $0\in \mathcal{T}(\M{x}_{\m{fix}},\M{y}_{\m{fix}})$, which implies that if the sequence $\{(\M{x}^k,\M{y}^k)\}$ generated by DISA \eqref{DISAUP} is convergent, it will converge to $\mathcal{M}^*$.
\begin{rmk}
In the development of DISA, instead of solving the dual problem \eqref{Dual1} directly, we minimize a quadratic approximation of $H(\M{y}')$ for the update of $\M{y}^{k+1}$.
This approach is a well-known technique in convex optimization, where successive quadratic approximations (or linearization) are used to approximate the objective function.
Various optimization methods, such as gradient descent, Linearized ALM (L-ALM), and Linearized ADMM (L-ADMM), can be interpreted from this perspective \cite{ADMM-D1}.
In the context of DISA, we intend to use the first-order gradient information of $H(\M{y})$ to determine the direction of the dual variable $\M{y}$ update.
\end{rmk}

The update \eqref{DISAUP} can be simplified by introducing $\tilde{\m{y}}^k_1=\sqrt{\M{V}}\m{y}^k_1$. Let $S_i=2\tau_i I+\frac{\tau_i(1-\tau\beta+\tau_i\beta)}{1-\tau\beta}U_iU_i\tr$, where $\tau=\max_{i}\{\tau_i\}$. With these notations and splitting the updates to the agents, we elaborate the implementation of DISA in Algorithm \ref{OurAlg2}. Note that \eqref{DISAUP} is a compact form of Algorithm \ref{OurAlg2} with $\tilde{\m{y}}^k_1=\sqrt{\M{V}}\m{y}^k_1$, and thus they are equivalent in the sense that they generate an identical sequence $\{\M{x}^k\}$. For communication costs, there is only one round of communication in each iteration, i.e., Step 4 requires neighboring variable $\bar{\m{x}}^{k+1}_{1,j}$.
\begin{algorithm}[!t]
  \caption{Dual Inexact Splitting Algorithm (DISA)} 
    \label{OurAlg2}
  \begin{algorithmic}[1]
    \Require
      Matrix $W$, step-sizes $0<\tau_i<{2}/{L_i}$ and $0<\tau_i\beta<1$.
    \State Initial $\m{x}_1^0\in \mathbb{R}^{mn}$, $\m{x}_2^0\in \mathbb{R}^{mp}$, $\tilde{\m{y}}_1^0=0$ and $\m{y}_2^0=0$.
    \For{$k=1,2,\ldots,K$}
      \State Agent $i$ computes $(\bar{\m{x}}^{k+1}_{1,i},\bar{\m{x}}^{k+1}_{2,i})$:
      \begin{align*}
      &\bar{\m{x}}^{k+1}_{1,i}=\m{x}^k_{1,i}-\tau_i \nabla f_i(\m{x}^k_{1,i})-\tau_i\tilde{\m{y}}^k_{1,i}-\tau_i U_i\tr \m{y}_{2,i}^k,\\
      &\bar{\m{x}}^{k+1}_{2,i}=\mathrm{prox}_{\tau_i g_i}(\m{x}^k_{2,i}+\tau_i \m{y}_{2,i}^k),
      \end{align*}
      and exchanges $\bar{\m{x}}^{k+1}_{1,i}$ with neighbors.
      \State Dual update $(\tilde{\m{y}}_{1,i}^{k+1},\m{y}_{2,i}^{k+1})$:
      \begin{equation*}
    \begin{aligned}
    &\tilde{\m{y}}^{k+1}_{1,i}=\tilde{\m{y}}^k_{1,i}+\frac{\beta}{2}\big(\bar{\m{x}}^{k+1}_{1,i}-\sum_{j\in \mathcal{N}_i}W_{ij}\bar{\m{x}}^{k+1}_{1,j}\big),\\
  &\m{y}^{k+1}_{2,i}=\m{y}^k_{2,i}+S_i^{-1}(U_i\bar{\m{x}}^{k+1}_{1,i}-\bar{\m{x}}^{k+1}_{2,i}).
    \end{aligned}
    \end{equation*}
      \State Primal update $(\m{x}_{1,i}^{k+1},\m{x}_{2,i}^{k+1})$:
      \begin{equation*}
      \begin{aligned}
    &\m{x}_{1,i}^{k+1}=\m{x}^k_{1,i}-\tau_i \nabla f_i(\m{x}^k_{1,i})-\tau_i\tilde{\m{y}}^{k+1}_{1,i}-\tau_i U_i\tr \m{y}_{2,i}^{k+1},\\
    &\m{x}_{2,i}^{k+1}=\mathrm{prox}_{\tau_i g_i}(\m{x}^k_{2,i}+\tau_i\m{y}_{2,i}^{k+1}).
      \end{aligned}
      \end{equation*}
    \EndFor
    \Ensure
      $\mathbf{x}^K$.
  \end{algorithmic}
\end{algorithm}
\begin{rmk}
When $g_i$ is proximable, the computational effort on DISA is generally dominated by the inverse of $S_i$. Consider the
update $\m{y}^{k+1}_{2,i}=\m{y}^k_{2,i}+S_i^{-1}(U_i\bar{\m{x}}^{k+1}_{1,i}-\bar{\m{x}}^{k+1}_{2,i})$. When $S_i$ is not large scale, since $S_i$ is positive definite, the inverse of it can be computed effectively, e.g., by the Cholesky decomposition. Moreover, in the whole iteration process, it only needs to be calculated once. However, when $S_i$ is large scale, the update can be reformulated as
$$
\m{y}_{2,i}^{k+1}=\mathop{\arg\min}\limits_{\m{y}_{2,i}\in\mathbb{R}^p}\{\frac{1}{2}\|\m{y}_{2,i}-\m{y}_{2,i}^k\|_{S_i}^2-\langle \m{y}_{2,i},U_i\bar{\m{x}}^{k+1}_{1,i}-\bar{\m{x}}^{k+1}_{2,i}\rangle\}.
$$
This is equivalent to a positive definite system of linear equations, which it can be solved effectively by the preconditioning conjugate gradient method or by Matlab function \texttt{quadprog}($\cdot$). Generally, we achieve the $\|\M{U}\M{U}\tr\|$-independent convergent step-size range at the cost of increasing the difficulty of dual updates.  Fortunately, since the dual update is equivalent to solving a strongly convex separable quadratic programming without constraints, the additional cost is acceptable. This is practically the distinctive merit of our proposed algorithm.
\end{rmk}
\begin{rmk}
The dual update \eqref{DualUP} can be rewritten as
$\M{y}^{k+1}=\arg\max_{\M{y}\in \mathbb{R}^{m(n+p)}}\{\mathcal{L}(\bar{\M{x}}^{k+1},\M{y})-\frac{1}{2}\|\M{y}-\M{y}^k\|_{\M{Q}}^2\}$.
Different from classic PD-PSAs, where the quadratic term is $\frac{1}{2\beta}\|\M{y}-\M{y}^k\|^2$, we specifically design a non-diagonal preconditioner $\M{Q}$ to change the metric of the ambient Euclidean space and enable the positive definiteness of the new metric matrix dependent of $\|\M{UU}\tr\|$, i.e.,
$$
\underbrace{\left(
  \begin{array}{cc}
    \tilde{\M{L}}_F & 0 \\
    0 & \frac{1}{\beta}I- \M{B}\Gamma\M{B}\tr \\
  \end{array}
\right)}_{\begin{array}{c}
            \frac{1}{2\beta}\|\M{y}-\M{y}^k\|^2 \\
            \M{y}^{k+1}=\M{y}^k+ \beta\M{B}\bar{\M{x}}^{k+1}
          \end{array}
}\longrightarrow\underbrace{\left(
  \begin{array}{cc}
    \tilde{\M{L}}_F & 0 \\
    0 & \M{Q}-\M{B}\Gamma\M{B}\tr \\
  \end{array}
\right):=\M{M}}_{\begin{array}{c}
                   \frac{1}{2}\|\M{y}-\M{y}^k\|_{\M{Q}}^2 \\
                   \M{y}^{k+1}=\M{y}^k+ \M{Q}^{-1}\M{B}\bar{\M{x}}^{k+1}
                 \end{array}
},
$$
where $\tilde{\M{L}}_F=\Gamma^{-1}-\frac{1}{2}\M{L}_F$. It can be easily verified that the metric matrix $\M{M}$ is positive definite, which is a crucial condition for ensuring the convergence of DISA, when the condition \eqref{Scond1} holds.
\end{rmk}
\begin{rmk}
Preconditioning is an important idea for the e development of primal-dual algorithms. In the case of Condat-Vu \cite{CV1,CV2}, the preconditioned FBS is employed to decouple updates of the primal and dual variables \cite[Section 3.3]{BookYin2022}. Additionally, PD3O \cite{PD3O} and PDDY \cite{DYS-T2}/AFBA \cite{AFBA} rely on the preconditioned Davis-Yin splitting \cite{DYS-T} to achieve a wider range of convergent step sizes than Condat-Vu \cite[Section 4]{DYS-T2}. On the other hand, preconditioning has been successfully used to accelerate the convergence of primal-dual hybrid gradient and ADMM \cite{Preconditioned1,Preconditioned2,Preconditioned3}. Moreover, preconditioning has been utilized to enhance certain properties of the equality constraint matrix in order to achieve acceleration. For example, in \cite{Preconditioned4} and \cite{Preconditioned5}, to overcome the sensitivity to ill conditioning, they reformulate the equality constraint $Ax=y$ as $EAx=Ey$, where $E$ is an invertible preconditioner and $EA$ is well-conditioned. In \cite{Preconditioned6}, to design the optimal edge weights for decentralized-ADMM \cite{DADMM}, it reformulates the equality constraint $Ax=By$ as $EAx=EBy$. In DISA, preconditioning is used to change the dual metric $\frac{1}{\beta}I-\M{B}\Gamma\M{B}\tr$ to $\M{Q}-\M{B}\Gamma\M{B}\tr$, whose positive definiteness is independent of $\|\M{UU}\tr\|$.
\end{rmk}
\begin{rmk}
For centralized COPs with one linear mapping, there are several algorithms that address the dependence of step-size on the linear mapping, such as the indefinite-proximal ALM \cite{InexactProximalALM1} and the indefinite-proximal ADMM \cite{InexactProximalALM1,InexactProximalALM2}. These algorithms implement convergence conditions that are independent of the linear mapping, but at the cost of increasing the difficulty of the primal update. As a result, the primal updates of these algorithms usually do not have an explicit iterative scheme and are no longer proximity-induced-featured. Thus, for their implementation, additional optimization algorithms need to be designed for the primal updates, depending on the specific characteristics of the problem. In our work, we address a more complex optimization problem (with two linear mapping $\M{V}$ and $\M{U}$, and needs to be solved in a distributed manner). To remove the dependence of $\M{V}$ (i.e., network-independent), we effectively use FBS, FMR duality, and linearization techniques. To remove the dependence of $\M{U}$, we design a non-diagonal preconditioner for the dual update. Our approach not only overcomes the dependence of $\M{V}$ and $\M{U}$, but also retains the proximity-induced feature of classic PD-PSAs avoiding additional inner iterations for subproblems. To the best of our knowledge, DISA is the first distributed algorithm to achieve this feature.
\end{rmk}

\subsection{Discussion}
\begin{table*}[!t]
\renewcommand\arraystretch{1.5}
\begin{center}
\caption{The Summarization of Iterations of Existing PD-PSAs and Algorithmic Frameworks}
\scalebox{1.1}{
\begin{tabular}{cc}
\hline
\textbf{Algorithm}&\textbf{Iteration Scheme}\\
\hline
\textbf{DPGA:}&
$
\M{x}^{k+1}=\mathrm{prox}^{\Gamma^{-1}}_{G}(\M{s}^k),
\M{s}^{k+1}=\M{s}^k+\M{W}(\M{x}^{k+1}-\M{x}^k)+(I-\M{W})\M{x}^{k+1}+\Gamma\nabla F(\M{x}^{k})-\Gamma\nabla F(\mathbf{x}^{k+1}).
$
\\
\textbf{PG-EXTRA:}&$
\mathbf{x}^{k+1}=\mathrm{prox}_{\tau G}(\mathbf{s}^{k}),
\mathbf{s}^{k+1}=\mathbf{s}^{k}-\M{x}^{k+1}+\tilde{\M{W}}(2\M{x}^{k+1}-\M{x}^{k})+\tau\nabla F(\M{x}^{k})-\tau\nabla F(\mathbf{x}^{k+1}).
$\\
\textbf{NIDS:}&
$
\mathbf{x}^{k+1}=\mathrm{prox}_{\tau G}(\mathbf{s}^{k}),
\mathbf{s}^{k+1}=\mathbf{s}^{k}-\M{x}^{k+1}+\tilde{\M{W}}\big(2\M{x}^{k+1}-\M{x}^{k}+\tau\nabla F(\M{x}^{k})-\tau\nabla F(\mathbf{x}^{k+1})\big).
$\\
\textbf{D-iPGM:}&
$
\tilde{\M{x}}^{k}=\mathrm{prox}_{\tau G}(\M{s}^{k}),
\M{x}^{k+1}=\tilde{\M{W}}\tilde{\M{x}}^{k},
\M{s}^{k+1}=\M{s}^{k}-\M{x}^{k}-\M{W}\tilde{\M{x}}^{k}+\tau\nabla F(\mathbf{x}^{k})-\tau\nabla F(\mathbf{x}^{k+1}).
$
\\
\hline
\hline
\textbf{Algorithmic Framework}&\textbf{Iteration Scheme}\\
\hline
$\mathrm{ABC}$-\textbf{Algorithm:}&$
\mathbf{x}^{k+1}=\mathrm{prox}_{\tau G}(\mathbf{s}^{k}),
\M{s}^{k+1}=(I-\mathrm{C})\M{s}^k+\mathrm{A}(\M{x}^{k+1}-\M{x}^k)+\tau \mathrm{B}(\nabla F(\M{x}^k)-\nabla F(\M{x}^{k+1})).
$\\
\textbf{PUDA}:& $
\mathbf{x}^{k+1}=\mathrm{prox}_{\tau G}(\mathrm{A}\mathbf{s}^{k}),
\mathbf{s}^{k+1}=(I-\mathrm{C}^2)\mathbf{s}^{k}+(I-\mathrm{C})(\mathbf{x}^{k+1}-\mathbf{x}^{k})
+\tau(\nabla F(\mathbf{x}^k)-\nabla F(\mathbf{x}^{k+1})).
$\\
\textbf{AMM:}&$
\M{x}^{k+1}={\arg\min}_{\M{x}\in \mathbb{R}^{nm}}\varpi^k(\M{x})+g(\M{x})+\frac{\rho}{2}\|\M{x}\|^2_{\mathrm{H}}+\langle \M{z}^k,\M{x}\rangle,
\M{z}^{k+1}=\M{z}^k+\rho \tilde{\mathrm{H}}\M{x}^{k+1}.
$\\
\hline
\end{tabular}}
\label{TableCOMPA}
\end{center}
\end{table*}
Consider another problem transformation of the problem \eqref{Problem}, i.e.,
$
\min_{\m{x}_1\in\mathbb{R}^{mn}}\{\sum_{i=1}^m  f_i(\m{x}_{1,i}) + \sum_{i=1}^m g_i(U_i\m{x}_{1,i}) + \delta_0(\sqrt{\M{V}}\m{x}_1)\}.
$
Here, the auxiliary variable $\m{x}_2$ has not been introduced. In fact, TPUS \cite{TPUS}, TriPD \cite{TriPD}, and PDFP-Dist \cite{PDFP-DIS} are provided based on this transformation. Let $\theta_1(\m{x}_1)=\sum_{i=1}^m  f_i(\m{x}_{1,i})$, $\theta_2(\tilde{\m{z}}_1,\tilde{\m{z}}_2)=\sum_{i=1}^m g_i(\tilde{\m{z}}_{1,i})+\delta_0(\tilde{\m{z}}_2)$, where $\tilde{\m{z}}_1\in \mathbb{R}^{mp},\tilde{\m{z}}_2\in\mathbb{R}^{mn}$, $\mathbf{C}=[\M{U}\tr,\sqrt{\M{V}}]\tr$, and $\M{C}\m{x}_1=(\M{U}\m{x}_1,\sqrt{\M{V}}\m{x}_1)\in \mathbb{R}^{(p+n)m}$. The problem has the form of $\min_{\m{x}_1}\{\theta_1(\m{x}_1)+\theta_2(\M{C}\m{x}_1)\}$, which can be well solved by Condat-Vu \cite{CV1,CV2} and PAPC-type algorithms (including PDFP \cite{PDFP}, AFBA \cite{AFBA}, PD3O \cite{PD3O}, and PDDY \cite{DYS-T2}). The parameter conditions for these algorithms solving this transformation are given in Table \ref{COMT}. Consider the convergent step-size range of L-ALM. Since $\|\M{B}\tr\M{B}\|\geq \max\{\|\M{U}\tr\M{U}+\M{V}\|,1\}$, the convergent step-size range of DISA is always larger than L-ALM. Then, consider the convergent step-size range of L-ADMM, Condat-Vu, TriPD, and PAPC-type. If $\|\M{U}\tr\M{U}\|\geq 1$, from $\|\M{U}\tr\M{U}+\M{V}\|\geq\|\M{U}\tr\M{U}\|\geq 1$, it deduces that the convergent step-size range of DISA is larger than these algorithms. Therefore, compared to these existing PD-PSAs, when $\|\M{U}\tr\M{U}\|\geq 1$, DISA extends the acceptable range of convergent step-size.
\begin{wrapfigure}{r}{4.3cm}
\centering
\includegraphics[width=0.22\textwidth]{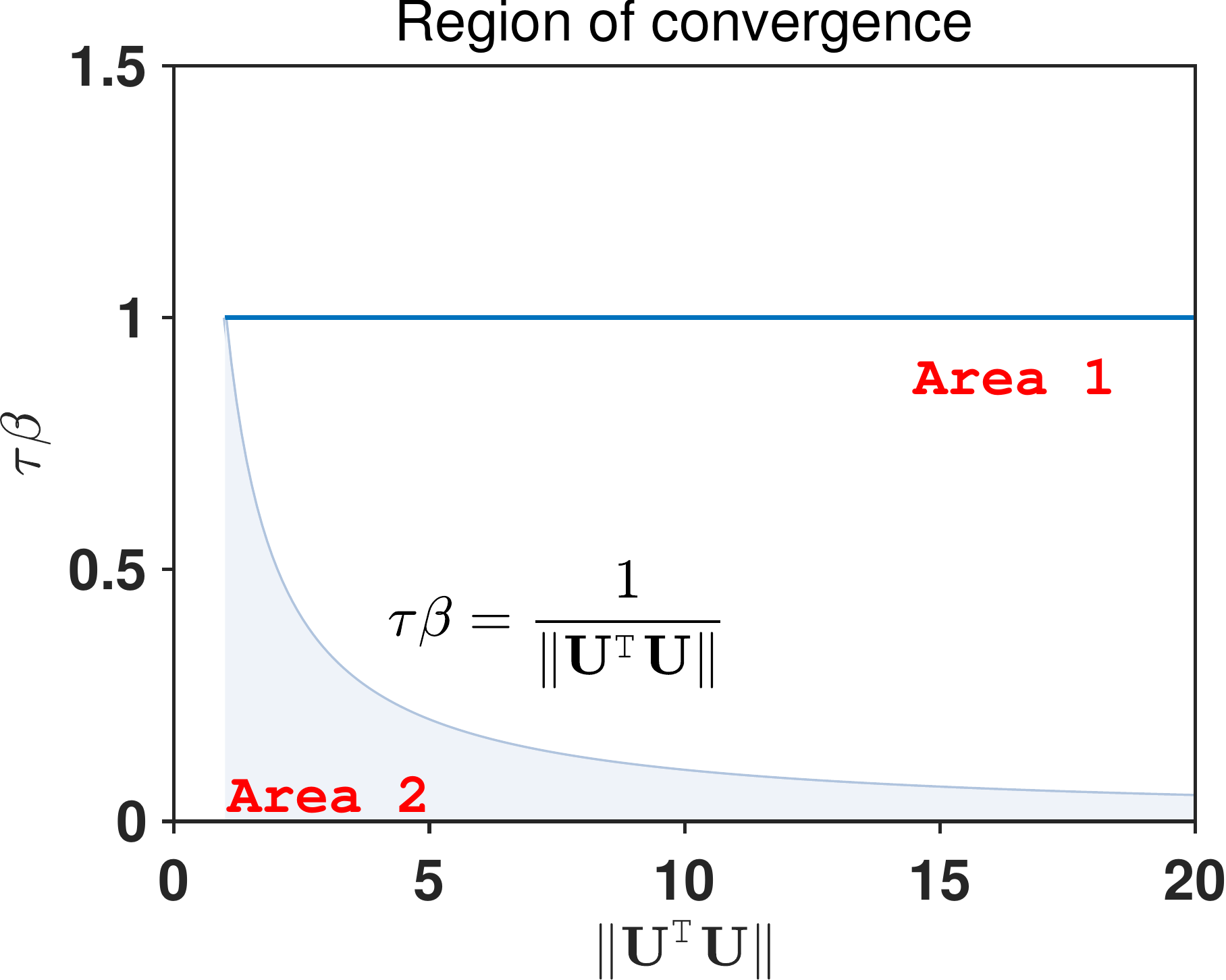}
\caption{Region of feasible step-sizes.}
\label{FigAR}
\end{wrapfigure}

Additionally, as presented in Table \ref{COMT} and Fig. \ref{FigAR}, it is apparent that the convergent step-sizes of the existing PD-PSAs are contingent on $\|\M{U}\M{U}\tr\|$. When $\|\M{U}\M{U}\tr\|$ is large, the efficacy of these PD-PSAs deteriorates. In particular, if $\|\M{U}\M{U}\tr\|$ is large, to satisfy the convergent step-size condition, $\tau\beta$ is forced to be small. Consequently, the step-sizes of the primal and dual updates decrease, necessitating more iterations. In contrast to these PD-PSAs that converge only in Area 2, DISA converges in both Area 1 and Area 2. Furthermore, the convergent step-size range is not limited by $\|\M{U}\M{U}\tr\|$, thus avoiding the use of small step-sizes even if $\|\M{U}\M{U}\tr\|$ is large. This feature of DISA allows it to be utilized in a wider range of applications, making it fundamentally distinct from these PD-PSAs.
\begin{table}[!t]
\renewcommand\arraystretch{1.5}
\begin{center}
\caption{The Comparison of Convergent Step-size Conditions}
\scalebox{0.9}{
\begin{tabular}{cc}
\hline
\textbf{~~~~Algorithm~~~~} &\textbf{Parameter Conditions}\\
\hline
\hline
L-ALM for problem \eqref{Trans-P1}&$0<\tau\beta\|\M{B}\tr\M{B}\|+\frac{\tau L}{2}<1$\\
L-ADMM, Condat-Vu, TriPD&$0<\tau\beta\|\M{U}\tr\M{U}+\M{V}\|+\frac{\tau L}{2}<1$\\
PAPC-type, TPUS, PDFP-Dist&$0<\tau<{2}/{L}\text{ and }0<\tau\beta\|\M{U}\tr\M{U}+\M{V}\|<1$\\
\hline
\cellcolor{myGray}\textbf{The proposed DISA}&\cellcolor{myGray} $0<\tau_i<{2}/{L_i}\text{ and }\tau_i\beta<1,i=1,\cdots,m$\\
\hline
\end{tabular}}
\label{COMT}
\end{center}
\end{table}
\subsection{Connections to Existing Algorithms}
This subsection compares DISA with several state-of-the-art primal-dual algorithms and algorithmic frameworks for distributed COP: $\m{x}^*\in\arg\min_{\m{x}}\{\sum_{i=1}^m (f_i(\m{x})+g_i(\m{x}))\}$. To this case, we do not introduce the auxiliary variable $\m{x}_2$. Thus, it holds that $\M{B}=\sqrt{\M{V}}$, $\M{x}=\m{x}_1$, $\M{y}=\m{y}_1$, and $\M{Q}=\frac{1}{\beta}I$. As in the previous derivation, introducing $\M{z}^k=\sqrt{\M{V}}\M{y}^k$, DISA updates in \eqref{DISAUP} take the following form:
\begin{align*}
\M{z}^{k+1}&=\M{z}^{k}+\beta\M{V}\mathrm{prox}^{\Gamma^{-1}}_{G}(\M{x}^k-\Gamma\nabla F(\M{x}^k)-\Gamma \M{z}^k), \\
\M{x}^{k+1}&=\mathrm{prox}^{\Gamma^{-1}}_{G}(\M{x}^k-\Gamma\nabla F(\M{x}^k)-\Gamma \M{z}^{k+1}).
\end{align*}
Eliminating the dual variable $\M{z}^k$, we have
\begin{align*}
\M{x}_c^{k}&=\Gamma\M{V}\mathrm{prox}_{G}^{\Gamma^{-1}}(\M{s}^k),\\
\M{x}^{k+1}&=\mathrm{prox}_{G}^{\Gamma^{-1}}(\M{s}^k-\beta\M{x}_c^{k}),\\
\M{s}^{k+1}&=\M{s}^{k}+\M{x}^{k+1}-\M{x}^k-\beta\M{x}^{k}_c+\Gamma\nabla F(\M{x}^{k})-\Gamma\nabla F(\M{x}^{k+1}).
\end{align*}
Recall several existing distributed primal-dual algorithms such as DPGA \cite{Aybat2017}, PG-EXTRA \cite{PGEXTR}, NIDS\cite{NDIS}, and D-iPGM \cite{Guo2022}. Similarly, by eliminating the dual variable, we outline their iterations in Table \ref{TableCOMPA}, where $\M{W}=W\otimes I_n$ and $\tilde{\M{W}}=\frac{1}{2}(I+\M{W})$. By comparison, we conclude that DISA is a novel algorithm distinct from these primal-dual algorithms. If $g_i\equiv0$, DISA updates in \eqref{DISAUP} take the following form:
\begin{align*}
\M{z}^{k+1}&=\M{z}^{k}+\beta\M{V}(\M{x}^k-\Gamma\nabla F(\M{x}^k)-\Gamma \M{z}^k), \\
\M{x}^{k+1}&=\M{x}^k-\Gamma\nabla F(\M{x}^k)-\Gamma \M{z}^{k+1}.
\end{align*}
Similarly, eliminating the dual variable $\M{z}^k$, we have
$$\M{x}^{k+1}=(I-\beta\Gamma\M{V})(2\M{x}^k-\M{x}^{k-1}+\Gamma\nabla F(\M{x}^{k-1})-\Gamma\nabla F(\M{x}^{k})).$$
When $\Gamma=\tau I$ and $\beta={1}/{\tau}$, it holds that $I-\beta\Gamma\M{V}=\frac{1}{2}(I+\M{W})=\tilde{\M{W}}$. Thus, DISA degenerates to NIDS ($g_i\equiv0$) \cite{NDIS}/D-iPGM ($g_i\equiv0$, $\varepsilon_k\equiv0$) \cite{Guo2022}/Exact Diffusion \cite{ExactDiffusion}.

Next, consider three unified frameworks of distributed primal-dual algorithms. When the nonsmooth function $g_i$ is common to all agents, we recall PUDA unified framework \cite{Sulaiman2021}, and $\mathrm{ABC}$-unified framework \cite{Xu2021}. Similarly, by eliminating the dual variable, the equivalent form of $\mathrm{ABC}$-unified framework and PUDA unified framework are presented in Table \ref{TableCOMPA}, where $\mathrm{A,B,C}$ are suitably chosen matrices. In contrast, DISA is independent of these two algorithmic frameworks. When the nonsmooth function $g_i$ may be distinct among these agents, a unifying approximate method of multipliers (AMM) has been proposed in \cite{AMM2022}, which is shown in Table \ref{TableCOMPA},
where $\varpi^k(\M{x})$ is a time-varying surrogate function, $\rho>0$ is a constant, $\mathrm{H}$ and $\tilde{\mathrm{H}}$ are suitably chosen matrices. For DISA, consider another equivalent form:
\begin{subequations}
\begin{align}
&\bar{\M{x}}^{k+1}=\mathop{\arg\min}\limits_{\M{x}\in \mathbb{R}^{nm}}\varpi^k(\M{x})+G(\M{x})+\langle \M{z}^k,\M{x}\rangle,\label{DISA-1}\\
&\M{z}^{k+1}=\M{z}^k+\beta\M{V}\bar{\M{x}}^{k+1}, \label{DISA-2}\\
&\M{x}^{k+1}=\mathop{\arg\min}\limits_{\M{x}\in \mathbb{R}^{nm}}\varpi^k(\M{x})+G(\M{x})+\langle \M{z}^{k+1},\M{x}\rangle, \label{DISA-3}
\end{align}
\end{subequations}
where $\varpi^k(\M{x})=\langle \nabla f(\M{x}^k),\M{x}-\M{x}^k\rangle+\frac{1}{2}\|\M{x}-\M{x}^k\|^2_{\Gamma^{-1}}$. Compared to AMM, DISA can be seen as a prediction-correction AMM, where the update \eqref{DISA-1} is the prediction step and the update \eqref{DISA-3} is the correction step. Thus, DISA cannot be recovered from AMM unified framework.

To summarize, compared to these existing state-of-the-art algorithms and algorithmic frameworks, DISA is a missing link in a group of distributed primal-dual algorithms.

\section{Convergence Analysis}\label{SEC-Convergence Analysis}
The convergence analysis will be conducted in the variational inequality context \cite{He2012}, and the forthcoming analysis of DISA is based on the following fundamental lemma.
\begin{lem}[\cite{He2012}]\label{LEMF}
Let $\theta_1(\m{x})$ and $\theta_2(\m{x})$ be proper closed convex functions. If $\theta_1$ is differentiable, and the solution set of the problem $\min\{\theta_1(\m{x})+\theta_2(\m{x}):\m{x}\in\mathbb{R}^n\}$ is nonempty, then it holds that $\m{x}^*\in \arg\min\{\theta_1(\m{x})+\theta_2(\m{x}):\m{x}\in \mathbb{R}^n\}$ if and only if
$
\theta_2(\m{x})-\theta_2(\m{x}^*)+\langle \m{x}-\m{x}^*,\nabla \theta_1(\m{x}^*)\rangle\geq 0, \forall \m{x}\in \mathbb{R}^n.
$
\end{lem}

A point $\mathbf{w}^*=[\M{x}^*,\M{y}^*]\in \mathcal{M}^* \subseteq \mathcal{M}:=\mathbb{R}^{m(n+p)}\times \mathbb{R}^{m(n+p)}$ is called a saddle point of the Lagrangian function $\mathcal{L}(\mathbf{x},\M{y})$ if
$$
\mathcal{L}\left(\mathbf{x}^{*}, \M{y}\right) \leq \mathcal{L}\left(\mathbf{x}^{*}, \M{y}^{*}\right) \leq \mathcal{L}\left(\mathbf{x}, \M{y}^{*}\right), \forall \M{w}:=[\mathbf{x}\tr,\M{y}\tr]\tr \in \mathcal{M},
$$
which can be alternatively rewritten as the following variational inequalities by Lemma \ref{LEMF}.
\begin{align*}
\text{\textbf{VI 1}: }&\Phi(\M{x})-\Phi(\M{x}^*)+\big\langle\mathbf{w}-\mathbf{w}^*,\mathcal{K}_1(\mathbf{w}^*)\big\rangle\geq 0, \forall \mathbf{w}\in \mathcal{M},\\
&\Longleftrightarrow\\
\text{\textbf{VI 2}: }&G(\M{x})-G(\M{x}^*)+\big\langle\mathbf{w}-\mathbf{w}^*,\mathcal{K}_2(\mathbf{w}^*)\big\rangle\geq 0, \forall \mathbf{w}\in \mathcal{M},
\end{align*}
where
$$
\mathcal{K}_1(\M{w})=\left(
                         \begin{array}{c}
                          {\M{B}}\tr\M{y} \\
                           -{\M{B}}\M{x} \\
                         \end{array}
                       \right),~\mathcal{K}_2(\M{w})=\left(
                         \begin{array}{c}
                           \nabla F(\M{x})+{\M{B}}\tr\M{y} \\
                           -{\M{B}}\M{x} \\
                         \end{array}
                       \right).
$$

\subsection{Global Convergence}
We first define some notations to simplify our analysis. More specifically, let $\M{w}^k=[(\mathbf{x}^k)\tr,(\M{y}^k)\tr]\tr$ and $\M{v}^k=[(\bar{\mathbf{x}}^k)\tr,(\M{y}^k)\tr]\tr$. Then, define two self-adjoint linear operators:
\begin{align*}
&\M{H}: \mathcal{M}\rightarrow \mathcal{M}:(\M{x},\M{y})\mapsto (\Gamma^{-1}\M{x},\M{Q}\M{y}),\\
&\M{M}: \mathcal{M}\rightarrow \mathcal{M}:(\M{x},\M{y})\mapsto ((\Gamma^{-1}-\frac{1}{2}\M{L}_F)\M{x},(\M{Q}- \M{B}\Gamma\M{B}\tr)\M{y}),
\end{align*}
where $\Gamma=\mathrm{diag}\{\Gamma_1,\Gamma_2\}$ is the step-size matrix with $\Gamma_1=\mathrm{diag}\{\tau_1I_n,\cdots,\tau_mI_n\}$ and $\Gamma_2=\mathrm{diag}\{\tau_1I_p,\cdots,\tau_mI_p\}$, and $\M{L}_F=\mathrm{diag}\{L_1{I}_{n+p},\cdots,L_m{I}_{n+p}\}$. In addition, letting $\tau=\max_{i\in \mathcal{V}}\{\tau_i\}$ and $L=\max_{i\in \mathcal{V}}\{L_i\}$. If the step-size condition \eqref{Scond1} holds, these two metric matrix are positive definite, and there exist two positive constants $c_1$ and $c_2$ such that for any $\M{w}_1,\M{w}_2\in \mathcal{M}$,
\begin{equation}\label{POSI}
\begin{aligned}
c_1\|\M{w}_1-\M{w}_2\|&\leq \|\M{w}_1-\M{w}_2\|_{\M{M}}\\
&\leq\|\M{w}_1-\M{w}_2\|_{\M{H}}\leq c_2\|\M{w}_1-\M{w}_2\|.
\end{aligned}
\end{equation}
To establish the global convergence and derive the convergence rate of DISA, we need the assertion in the following lemma.
\begin{lem}\label{F1}
Suppose that Assumption \ref{ass1} and the step-size condition \eqref{Scond1} hold. For $\forall (\M{x},\M{y})\in \mathcal{M}$, the sequence $\{\bar{\M{x}}^k,\M{x}^k,\M{y}^k\}$ generated by DISA \eqref{DISAUP} satisfies
\begin{align}\label{BASCI-IEQ1}
&G(\bar{\M{x}}^{k+1})-G(\M{x})+\langle \M{v}^{k+1}-\M{w},\mathcal{K}_2(\M{w})\rangle
+\frac{1}{2}\|\M{w}^k-\M{v}^{k+1}\|^2_{\M{M}}\nonumber\\
&\leq \frac{1}{2}\|\M{w}^k-\M{w}\|_{\M{H}}^2-\frac{1}{2}\|\M{w}^{k+1}-\M{w}\|_{\M{H}}^2, \forall (\M{x},\M{y})\in \mathcal{M}.
\end{align}
Furthermore, if $0<\tau_i<{1}/{L_i}$, it holds that
\begin{align}\label{BASCI-IEQ4}
&\mathcal{L}(\bar{\M{x}}^{k+1},\M{y})-\mathcal{L}(\M{x},\M{y}^{k+1})+\frac{1}{2}\|\M{w}^k-\M{v}^{k+1}\|^2_{\M{M}_1}\nonumber\\
&\leq \frac{1}{2}\|\M{w}^k-\M{w}\|_{\M{H}}^2-\frac{1}{2}\|\M{w}^{k+1}-\M{w}\|_{\M{H}}^2, \forall \M{w}\in \mathcal{M},
\end{align}
where $\M{M}_1: \mathcal{M}\rightarrow \mathcal{M}:(\M{x},\M{y})\mapsto ((\Gamma^{-1}-\M{L}_F)\M{x},(\M{Q}-\M{B}\Gamma\M{B}\tr)\M{y})$ is a self-adjoint linear operator.
\end{lem}
\begin{IEEEproof}
See Appendix \ref{APPA}.
\end{IEEEproof}
It follows from \eqref{BASCI-IEQ1} and \textbf{VI 2} that for any $\M{w}^*\in \mathcal{M}^*$
\begin{align}\label{BASCISS}
&\frac{1}{2}\|\M{w}^k-\M{w}^*\|_{\M{H}}^2-\frac{1}{2}\|\M{w}^{k+1}-\M{w}^*\|_{\M{H}}^2-\frac{1}{2}\|\M{w}^k-\M{v}^{k+1}\|^2_{\M{M}}\nonumber\\
&\geq G(\bar{\M{x}}^{k+1})-G(\M{x}^*)+\langle\mathbf{v}^{k+1}-\mathbf{w}^*,\mathcal{K}_2(\mathbf{w}^*)\rangle \geq 0,
\end{align}
which implies that $\|\M{w}^k-\M{w}^*\|_{\M{H}}^2\geq\|\M{w}^{k+1}-\M{w}^*\|_{\M{H}}^2,\forall \M{w}^*\in\mathcal{M}^*$, i.e., the sequence $\{\M{w}^k\}$ generated by DISA is a Fej\'{e}r monotone sequence with respect to $\mathcal{M}^*$ in $\mathbf{H}$-norm. Moreover, we have $\|\M{w}^k-\M{w}^*\|^2_{\M{H}}\leq\|\M{w}^0-\M{w}^*\|^2_{\M{H}},\forall k\geq0$. By the positive definiteness of $\M{H}$, we have the sequence $\{\M{w}^k\}$ generated by DISA is bounded.

Then, we give the following theorem to prove that the sequence $\{(\M{x}^k,\M{y}^k)\}$ generated by DISA \eqref{DISAUP} converges to a primal-dual solution of problem \eqref{Trans-P1}.
\begin{thm}\label{THE-1}
Suppose that Assumption \ref{ass1} and the step-size condition \eqref{Scond1} hold. Let $\{\M{w}^k\}$ be the sequence generated by DISA \eqref{DISAUP}. There exists $\M{w}^{\infty}\in \mathcal{M}^*$ such that
$\lim_{k\rightarrow\infty} \M{w}^k = \M{w}^{\infty}$. Since $\M{w}^k=((\m{x}_1^k,\m{x}_2^k),(\m{y}_1^k,\m{y}_2^k))$, we further have $\lim_{k\rightarrow\infty} \m{x}_1^{k}=1_m\otimes \m{x}^{\infty}$, where $\m{x}^{\infty}$ solves problem \eqref{Problem}.
\end{thm}
\begin{IEEEproof}
See Appendix \ref{APPB}.
\end{IEEEproof}

\begin{table*}[!t]
\renewcommand\arraystretch{1.5}
\begin{center}
\caption{Some Commonly Used Convex Loss Functions Making the Metric Subregularity of $\mathcal{T}$ Holds}
\scalebox{0.9}{
\begin{tabular}{ccccc}
\hline
\multirow{2}{*}{$h_i(\m{x}_{1,i})$}&\textbf{linear regression}&\textbf{logistic regression}&\textbf{likelihood estimation}&\textbf{poisson regression}\\
\cline{2-5}
~&$\frac{1}{2}\|\mathrm{x}_{1,i}-\mathrm{b}^1_i\|^2$&$\sum_{j=1}^{n}\log(1+e^{\mathrm{x}_{1,ij}})-(\mathrm{b}_i^2)^{\top}\mathrm{x}_{1,i}$&
$\sum_{j=1}^{n}\log(1+e^{(\mathrm{b}_i^3)^{\top}\mathrm{x}_{1,i}})$&$\sum_{j=1}^{n}e^{\mathrm{x}_{1,ij}}-(\mathrm{b}_i^4)^{\top}\mathrm{x}_{1,i}$\\
\hline
\multirow{2}{*}{$g_i(\m{x}_{2,i})$}&\textbf{$l_{\infty}$-norm}&\textbf{$l_{1}$-norm}&\textbf{elastic net}&\textbf{OSCAR}\\
\cline{2-5}
~&$\nu\|\mathrm{x}_{2,i}\|_{\infty}$&$\nu\|\mathrm{x}_{2,i}\|_{1}$&$\nu_1\|\mathrm{x}_{2,i}\|_1+\nu_2\|\mathrm{x}_{2,i}\|^2$&$\nu_{1}\|\mathrm{x}_{2,i}\|_{1}+\nu_{2} \sum_{k<l} \max \left\{|\mathrm{x}_{2,ik}|,|\mathrm{x}_{2,il}|\right\}$\\
\hline
\end{tabular}}
\label{Table2}
\end{center}
\end{table*}

\subsection{Sublinear Convergence Rate}
In this subsection, with general convexity, we establish $O({1}/{k})$ non-ergodic and ergodic convergence rates of DISA.

With the inequalities established in the previous subsections, one can obtain the following theorem.
\begin{thm}\label{THE-2}
Suppose that Assumption \ref{ass1} and the step-size condition \eqref{Scond1} hold. The sequence $\{\M{v}^k\}$ and $\{\M{w}^k\}$ generated by DISA \eqref{DISAUP} satisfies
\begin{enumerate}
  \item Running-average successive difference:
  \begin{align}\label{RA-RATE1}
  \frac{1}{K}\sum_{k=0}^{K-1}\|\M{w}^{k}-\M{v}^{k+1}\|^2_{\M{M}}=O\left({1}/{K}\right).
  \end{align}
  \item Running-best successive difference:
  \begin{align}\label{RB-RATE1}
  \min_{0\leq k\leq K-1} \{\|\M{w}^k-\M{v}^{k+1}\|^2_{\M{M}}\}= o\left({1}/{K}\right).
  \end{align}
  \item Running-average first-order optimality residual:
  \begin{align}\label{RA-RATE2}
  \frac{1}{K}\sum_{k=0}^{K-1}\mathrm{dist}^2(0,\mathcal{T}(\M{v}^{k+1}))=O\left({1}/{K}\right).
  \end{align}
  \item Running-best first-order optimality residuals:
  \begin{align}\label{RB-RATE2}
  \min_{0\leq k\leq K-1} \{\mathrm{dist}^2(0,\mathcal{T}(\M{v}^{k+1}))\} = o\left({1}/{K}\right).
  \end{align}
\end{enumerate}
\end{thm}
\begin{IEEEproof}
See Appendix \ref{APPC}.
\end{IEEEproof}
If $\M{w}^k=\M{w}^{k+1}$, i.e., $\M{w}^k$ is the fixed point of DISA, we have $0\in\mathcal{T}(\M{w}^k)$, which meets the first order optimality condition of problem \eqref{Trans-P1}. Thus $\M{x}^k$ is a minimizer of problem \eqref{Trans-P1}. On the other hand, it follows from $\|\M{x}^{k+1}-\bar{\M{x}}^{k+1}\|\leq\tau\|\sqrt{\mathbf{B}}(\M{y}^k-\M{y}^{k+1})\|$ that $\|\M{w}^k-\M{w}^{k+1}\|^2_{\M{M}}=0\Leftrightarrow \|\M{w}^k-\M{v}^{k+1}\|^2_{\M{M}}=0$. Hence, $\|\mathbf{w}^k-\mathbf{v}^{k+1}\|^2_{\M{M}}$ can be viewed as an error measurement after $k$ iterations of the DISA.

Then, we establish the convergence rates for primal-dual gap, primal suboptimality, and consensus violation.
\begin{thm}\label{THE-PDGAP}
Suppose that Assumption \ref{ass1} holds, and $0<\tau_i<1/L_i$, $0<\tau_i\beta<1, i\in\mathcal{V}$.
The sequence $\{\bar{\M{x}}^k,\M{x}^k,\M{y}^k\}$ generated by DISA \eqref{DISAUP} satisfies for $\forall \M{w}\in \mathcal{M}$
\begin{align}\label{PDGAP}
\mathcal{L}(\M{X}^{K},\M{y})-\mathcal{L}(\M{x},\M{Y}^{K})\leq \frac{1}{2K}\|\M{w}^0-\M{w}\|_{\M{H}}^2,
\end{align}
where $\M{X}^{K}=\frac{1}{K}\sum_{k=1}^{K}\bar{\M{x}}^k,\text{ and } \M{Y}^{K}=\frac{1}{K}\sum_{k=1}^{K}\M{y}^k$. Moreover, we have
\begin{align}
\Big|\sum_{i=1}^{m}(f_i(\m{X}_1^K)+g_i(\m{X}_2^K))-J^*\Big|&=O\left({1}/{K}\right),\\
\big\|((I-W)\otimes I_m)^{\frac{1}{2}}\m{X}_1^K\big\|&=O\left({1}/{K}\right) \label{consensusviolation},
\end{align}
where $\m{X}_1^{K}=\frac{1}{K}\sum_{k=1}^{K}\bar{\m{x}}_1^k$, $\m{X}_2^{K}=\frac{1}{K}\sum_{k=1}^{K}\bar{\m{x}}_2^k$, and $J^*$ is defined in \eqref{Problem}.
\end{thm}
\begin{IEEEproof}
See Appendix \ref{APP-PDGAP}.
\end{IEEEproof}

\subsection{Linear Convergence Rate}
In this subsection, the linear convergence of DISA is proven under metric subregularity \cite{va2004}, which is an important assumption on the establishment of linear convergence rate \cite{va2020,va2021}. With metric subregularity, \cite{Liang2019} and \cite{Pan2022} establish the linear convergence rate of primal-dual gradient dynamics and distributed ADMM, respectively. Similarly, we start to recall the basic concepts in variational analysis.
\begin{defn}
A set-valued mapping $\Psi:\mathbb{R}^{n}\rightrightarrows \mathbb{R}^{m}$ is metrically subregular at $(\bar{u},\bar{v})\in \mathrm{gph}(\Psi)$ if for some $\epsilon>0$ there exists $\kappa\geq0$ such that
$$
\mathrm{dist}(u,\Psi^{-1}(\bar{v}))\leq \kappa ~\mathrm{dist}(\bar{v},\Psi(u)), \forall u\in \mathcal{B}_{\epsilon}(\bar{u}),
$$
where $\mathrm{gph}(\Psi):=\{(u,v):v=\Psi(u)\}$, $\Psi^{-1}(v):=\{u\in\mathbb{R}^n:(u,v)\in\mathrm{gph}(\Psi)\}$.
\end{defn}

Then, we give the following theorem to present the local convergence rate of DISA.
\begin{thm}\label{THE-3}
Suppose that Assumption \ref{ass1} and the step-size condition \eqref{Scond1} hold, and the sequence $\{\M{w}^k\}$ generated by DISA \eqref{DISAUP} converges to $\M{w}^{\infty}$. If $\mathcal{T}$ is metrically subregular at $(\M{w}^{\infty},0)$ with modulus $\kappa_2$, there exists ${K}>0$ such that
\begin{align}\label{THE111}
\mathrm{dist}(\M{v}^{k+1},\mathcal{M}^*)\leq \kappa_1\kappa_2\|\M{v}^{k+1}-\M{w}^k\|_{\M{M}}, \forall k\geq {K},
\end{align}
where $\kappa^2_1=\frac{1}{c_1^2}\max\{3L^2+\frac{3}{\tau^2},3\|\M{B}\|^2+\|\M{Q}\|^2\}$. Moreover, it holds that
$$
\mathrm{dist}_{\M{H}}(\M{w}^{k+1},\mathcal{M}^*) \leq  \varrho ~ \mathrm{dist}_{\M{H}}(\M{w}^k,\mathcal{M}^*),\forall k\geq {K},
$$
where $\varrho=\sqrt{1-\frac{c_1^2}{c^2_2(\kappa_1\kappa_2c_2+1)^2}}<1$ and $c_1,c_2$ are constant given in \eqref{POSI}.
Furthermore, it holds that for $k>{K}$
$$
\|\M{w}^{k}-\M{w}^{\infty}\|_{\M{M}}\leq\frac{2\varrho^{k-{K}}}{1-\varrho} \mathrm{dist}_{\M{H}}(\M{w}^{{K}},\mathcal{M}^*).
$$
\end{thm}
\begin{IEEEproof}
See Appendix \ref{APPD}.
\end{IEEEproof}
\begin{rmk}
In \cite{Sulaiman2021}, it is reported that if each agent owns a different local nonsmooth term, a dimension dependent linear convergence may be attained in the worst case (the smooth term $f_i$ is strongly convex and the nonsmmoth term $g_i$ is proximable). However, to Theorem \ref{THE-3}, the dimension dependent global linear rate of DISA can be established, which thus does not contradict with the exiting results.
\end{rmk}
Next, referring to \cite{va2020} and \cite{va2021}, we give the sufficient condition such that $\mathcal{T}(\mathbf{w})$ is metrically subregular at $(\mathbf{w}^{\infty},0)$. A convex function $f(\m{x}_1)$ is said to satisfy the structured assumption if $f(\m{x}_1)=\sum_{i=i}^{m}f_i(\mathrm{x}_{1,i})=\sum_{i=1}^{m}(h_i(A_i\mathrm{x}_{1,i})+\langle q_i,\mathrm{x}_{1,i}\rangle)$,
where $A_i$ is a $m_i\times n$ matrix, $q_i$ is a vector in $\mathbb{R}^n$, and $h_i$ is smooth and essentially
locally strongly convex, i.e., for any compact and convex subset $\mathbb{K}$, $h_i$ is strongly convex on $\mathbb{K}$. Recall that a set-valued mapping is called a polyhedral multifunction if its graph is the union of finitely many convex polyhedra. Then, one has that a function $\phi:\mathbb{R}^n\rightarrow \mathbb{R}$ is called piecewise linear-quadratic if and only if $\partial \phi$ is a polyhedral multifunction. From \cite[Theorem 59]{va2020}, it holds that if $f(\m{x}_1)$ meets the structured assumption and $g(\m{x}_2)$ is convex piecewise linear-quadratic function, the mapping $\mathcal{T}$ is metrically subregular at $(\mathbf{w}^{\infty},0)$. Some commonly used loss functions and regularizers in machine learning satisfying the assumption that $\mathcal{T}$ is metrically subregular at $(\M{w}^{\infty},0)$ are summarized in Table \ref{Table2}, where $\mathrm{b}^1_i\in \mathbb{R}^n$,  $\mathrm{b}^2_i\in \{0,1\}^{n}$, $\mathrm{b}^3_i\in \{-1,1\}^{n}$, $\mathrm{b}_i^4\in\{0,1,2,\cdots\}$, and $\nu$, $\nu_1$ and $\nu_2$ are the given nonnegative parameters.

\section{DISA with Approximate Proximal Mapping}\label{SECIN}
In this section, we provide V-DISA that allows to use approximate proximal mappings.
\subsection{Algorithm Development}
Note that in the iteration of DISA \eqref{DISAUP}, two proximal operators of $g_i$ are required in each iteration. If the proximal operator of $g_i$ has no analytical solution, the additional cost of one proximal mapping of DISA can not be ignored. Let $\bm{\xi}^k=\M{x}^k-\Gamma\nabla F(\M{x}^k)-\Gamma {\M{B}}\tr \M{y}^k$. The update \eqref{PrimalUP} can be written as
\begin{align*}
\M{x}^{k+1}=\mathrm{prox}^{\Gamma^{-1}}_{G}(\bm{\xi}^k+\Gamma \M{B}\tr(\M{y}^k-\M{y}^{k+1})).
\end{align*}
To reduce the computational cost, we change it as follows.
\begin{align}\label{V-DISA-primalup}
\M{x}^{k+1}=\mathrm{prox}^{\Gamma^{-1}}_{ G}(\bm{\xi}^k)+\Gamma \M{B}\tr(\M{y}^k-\M{y}^{k+1}).
\end{align}
Because of this small modification, only one proximal mapping of $g_i$ is needed in each iteration.

If the proximal mapping of $g_i$ does not have an analytical solution, achieving high accuracy in solving subproblem \eqref{V-DISA-primalup} is inefficient. Instead, it is more practical to establish a suitable condition that, when met, terminates the subproblem procedure. To this end, we introduce an absolutely summable error criteria for subproblem \eqref{V-DISA-primalup}, and develop V-DISA with approximate proximal mapping (Algorithm \ref{OurAlg3}). This approach allows for the optimization of a problem in a distributed manner, where individual agents can work on their own independent subproblems with their own error criteria. The ability to work independently results in an improved overall efficiency in solving the problem. Additionally, to control the residual error $\m{d}_{i}^k$ of each subproblem at $k$-th iteration, a summable sequence $\{\varepsilon_k\}$ is introduced. This sequence can be selected in various ways, such as $\varepsilon_k=\varepsilon_0/(k+1)^r$ for any $\varepsilon_0>0$ and $r>1$, or $\varepsilon_k=r^{k}$ with $0<r<1$. The choice of $\varepsilon_k$ determines the rate at which the residual error decreases over time, which can be adjusted according to the specific requirements of the problem at hand.

Algorithm \ref{OurAlg3} is equivalent to the following iterations in the sense that
they generate an identical sequence $\{\M{x}^k\}$,
\begin{subequations}\label{Inexact}
\begin{align}
\tilde{\M{x}}^{k+1}&=\mathrm{prox}^{\Gamma^{-1}}_{ G}(\M{x}^k-\Gamma\nabla F(\M{x}^k)-\Gamma {\M{B}}\tr \M{y}^k+\Gamma \M{d}^k),\label{Inexact1}\\
\M{y}^{k+1}&=\mathop{\arg\min}\limits_{\M{y}\in\mathbb{R}^{m(n+p)}}\big\{\frac{1}{2}\|\M{y}-\M{y}^k\|^2_{\M{Q}}-\langle\M{B}\tr\M{y},\tilde{\M{x}}^{k+1}\rangle\big\}, \label{Inexact2}\\
\M{x}^{k+1}&=\tilde{\M{x}}^{k+1}+\Gamma \M{B}\tr(\M{y}^k-\M{y}^{k+1}), \label{Inexact3}
\end{align}
\end{subequations}
where $\M{d}^k=[0,(\m{d}^k)\tr]\tr$, $\m{d}^k=[(\m{d}_1^k)\tr,\cdots,(\m{d}_m^k)\tr]\tr$. Let $\M{d}^k\equiv0$, and $(\M{x}_{\mathrm{fix}},\M{y}_{\mathrm{fix}})$ be the fixed point of V-DISA. It holds that $0\in \mathcal{T}(\M{x}_{\m{fix}},\M{y}_{\m{fix}})$, which means that $(\M{x}_{\mathrm{fix}},\M{y}_{\mathrm{fix}})\in\mathcal{M}^*$.
\begin{algorithm}[!t]
  \caption{DISA with approximate proximal mapping} 
    \label{OurAlg3}
  \begin{algorithmic}[1]
    \Require
      Matrix $W$, step-sizes $0<\tau_i<{2}/{L_i}$ and $0<\tau_i\beta<1$, a summable sequence of nonnegative numbers $\{\varepsilon^k\}$.
    \State Initial $\m{x}_1^0\in \mathbb{R}^{mn}$, $\m{x}_2^0\in \mathbb{R}^{mp}$, $\tilde{\m{y}}_1^0=0$ and $\m{y}_2^0=0$.
    \For{$k=1,2,\ldots,K$}
      \State Agent $i$ computes
      $$
      \tilde{\m{x}}^{k+1}_{1,i}=\m{x}^k_{1,i}-\tau_i \nabla f_i(\m{x}^k_{1,i})-\tau_i\tilde{\m{y}}^k_{1,i}-\tau_i U_i\tr \m{y}_{2,i}^k,
      $$
      and computes $\tilde{\m{x}}^{k+1}_{2,i}\approx\mathrm{prox}_{\tau_i g_i}(\m{x}^k_{2,i}+\tau_i \m{y}_{i}^k)$, such that there exists a vector
      \begin{align*}
      &\m{d}_{i}^k\in \partial g_i(\tilde{\m{x}}^{k+1}_{2,i})-\m{y}_{2,i}^k+ \frac{1}{\tau_i}(\tilde{\m{x}}^{k+1}_{2,i}-\m{x}_{2,i}^k),\\
      &\text{satisfying } \|\m{d}_{2,i}^k\|\leq \varepsilon^k.
      \end{align*}
      Then exchanges $\tilde{\m{x}}^{k+1}_{1,i}$ with neighbors.
      \State Dual update $(\tilde{\m{y}}_{1,i}^{k+1},\m{y}_{2,i}^{k+1})$:
      \begin{equation*}
    \begin{aligned}
    &\tilde{\m{y}}^{k+1}_{1,i}=\tilde{\m{y}}^k_{1,i}+\frac{\beta}{2}\big(\tilde{\m{x}}^{k+1}_{1,i}-\sum_{j\in \mathcal{N}_i}W_{ij}\tilde{\m{x}}^{k+1}_{1,j}\big),\\
    &\m{y}^{k+1}_{2,i}=\m{y}^k_{2,i}+S_i^{-1}(U_i\tilde{\m{x}}^{k+1}_{1,i}-\tilde{\m{x}}^{k+1}_{2,i}).
    \end{aligned}
    \end{equation*}
      \State Primal update $(\m{x}_{1,i}^{k+1},\m{x}_{2,i}^{k+1})$:
      \begin{equation*}
      \begin{aligned}
    &\m{x}_{1,i}^{k+1}=\tilde{\m{x}}^{k+1}_{1,i}+\tau_i(\tilde{\m{y}}_{1,i}^k-\tilde{\m{y}}_{1,i}^{k+1} + U_i\tr(\m{y}_{2,i}^k-\m{y}_{2,i}^{k+1})),\\
    &\m{x}_{2,i}^{k+1}=\tilde{\m{x}}^{k+1}_{2,i}-\tau_i(\m{y}_{2,i}^k-\m{y}_{2,i}^{k+1}).
      \end{aligned}
      \end{equation*}
    \EndFor
    \Ensure
      $\mathbf{x}^K$.
  \end{algorithmic}
\end{algorithm}

\subsection{Convergence Analysis}
Let $\M{w}^k=[(\mathbf{x}^k)\tr,(\M{y}^k)\tr]\tr$ and $\M{v}^k=[(\tilde{\mathbf{x}}^k)\tr,(\M{y}^k)\tr]\tr$. To establish the global convergence and the convergence rate of V-DISA \eqref{Inexact}, we give the following lemma.
\begin{lem}\label{InLem}
Suppose that Assumption \ref{ass1} and the step-size condition \eqref{Scond1} hold. The sequence $\{\tilde{\M{x}}^k,\M{x}^k,\M{y}^k\}$ generated by \eqref{Inexact} satisfies that
\begin{align}\label{BASCI-IEQ66}
&2G(\tilde{\M{x}}^{k+1})-2G(\M{x})+2\langle \M{v}^{k+1}-\M{w},\mathcal{K}_2(\M{w})\rangle\nonumber\\
&\leq \|\M{w}^k-\M{w}\|_{\M{H}}^2-\|\M{w}^{k+1}-\M{w}\|_{\M{H}}^2-\|\M{w}^k-\M{v}^{k+1}\|^2_{\M{M}}\nonumber\\
&\quad+2\langle \tilde{\M{x}}^{k+1}-\M{x},\M{d}^k\rangle, \forall (\M{x},\M{y})\in \mathcal{M}.
\end{align}
Furthermore, if $0<\tau_i<{1}/{L_i}$, it holds that for $\forall \M{w}\in \mathcal{M}$
\begin{align}\label{INX66}
&2\mathcal{L}(\tilde{\M{x}}^{k+1},\M{y})-2\mathcal{L}(\M{x},\M{y}^{k+1})+\|\M{w}^k-\M{v}^{k+1}\|^2_{\M{M}_1}\nonumber\\
&\leq \|\M{w}^k-\M{w}\|_{\M{H}}^2-\|\M{w}^{k+1}-\M{w}\|_{\M{H}}^2+2\langle \tilde{\M{x}}^{k+1}-\M{x},\M{d}^k\rangle.
\end{align}
\end{lem}
\begin{IEEEproof}
See Appendix \ref{APPIn1}.
\end{IEEEproof}

Then, we prove the global convergence of V-DISA \eqref{Inexact} in the following theorem.
\begin{thm}\label{THE-4}
Suppose that Assumption \ref{ass1} and the step-size condition \eqref{Scond1} hold. Let $\{(\tilde{\M{x}}^k,\M{x}^k,\M{y}^k)\}$ be the sequence generated by V-DISA \eqref{Inexact}. If $\sum_{k=0}^{\infty}\varepsilon^k<\infty$, for any $\M{w}^*\in \mathcal{M}^*$, $\{\|\M{w}^k-\M{w}^*\|\}$ and $\{\|\M{v}^k-\M{w}^*\|\}$ are bounded. Moreover, there exists $\M{w}^{\infty}\in \mathcal{M}^*$ such that
$\lim_{k\rightarrow\infty} \M{w}^k = \M{w}^{\infty}$. Since $\M{w}^k=((\m{x}_1^k,\m{x}_2^k),(\m{y}_1^k,\m{y}_2^k))$, we have $\lim_{k\rightarrow\infty} \m{x}_1^{k}=1_m\otimes \m{x}^{\infty}$, where $\m{x}^{\infty}$ solves problem \eqref{Problem}.
\end{thm}
\begin{IEEEproof}
See Appendix \ref{APPG}.
\end{IEEEproof}

Different from the original DISA \eqref{DISAUP}, the sequence $\{\M{w}^k\}$ generated by \eqref{Inexact} satisfying
\begin{align*}
\|\M{w}^{k+1}-\M{w}^*\|_{\M{H}}\leq \|\M{w}^k-\M{w}^*\|_{\M{H}} +\psi \varepsilon^k, \psi>0,\M{w}^*\in \mathcal{M},
\end{align*}
which implies that the sequence $\{\|\M{w}^k-\M{w}^*\|_{\mathbf{H}}\}$ is quasi-Fej\'{e}r monotone. Note that, to V-DISA, the convergent step-size range is also independent of $\|\M{U}\M{U}\tr\|$.

Similar as Theorem \ref{THE-2}, with general convexity and a summable absolute error criterion we establish the sublinear convergence rate of the approximate version.
\begin{thm}\label{THE-5}
Suppose that Assumption \ref{ass1} and the step-size condition \eqref{Scond1} hold, and $\sum_{k=0}^{\infty}\varepsilon^k=\bar{\varepsilon}<\infty$. There exists a constant $0<\zeta<\infty$ such that sequence $\{\tilde{\M{x}}^k,\M{x}^k,\M{y}^k\}$ generated by DISA \eqref{Inexact} satisfies
\begin{align}\label{INRATA2}
\frac{1}{K}\sum_{k=0}^{K-1}\mathrm{dist}^2(\mathbf{0},\mathcal{T}({\M{v}}^{k+1}))\leq\frac{\psi_1\kappa^2_3+4\bar{\varepsilon}^2}{K},
\end{align}
where
\begin{align*}
&\kappa^2_3=\frac{1}{c_1^2}\max\{4L^2+\frac{4}{\tau^2},4\|\M{B}\|^2+\|\M{Q}\|^2\},\\
&\psi_1=\|\M{w}^{0}-\M{w}^*\|_{\M{H}}^2+2 \tau\sigma\sqrt{2m}\big(\|\M{w}^{0}-\M{w}^*\|_{\M{H}} +\bar{\varepsilon}\big) \bar{\varepsilon},\\
&\sigma=\max\{\sqrt{\|\M{H}^{-1}\|}\psi,\tau\sqrt{2m}(1+\|\M{Q}^{-1}\M{B}\|)\},\\ &\psi=\tau\sqrt{2m}\|\M{H}\|^{\frac{1}{2}}(\|\M{Q}^{-1}\M{B}\|(1+\tau\|\M{B}\|)+1).
\end{align*}
Furthermore, if $0<\tau_i<{1}/{L_i}$, we have
\begin{align}
\Big|\sum_{i=1}^{m}(f_i(\m{X}_1^K)+g_i(\m{X}_2^K))-J^*\Big|&\leq \frac{(1+\omega_1)\omega_2+\omega_1\bar{\varepsilon}}{2K},\label{TH6-28}\\
\big\|((I-W)\otimes I_m)^{\frac{1}{2}}\m{X}_1^K\big\|&\leq \frac{(1+\omega_1)\omega_2+\omega_1\bar{\varepsilon}}{2K}, \label{TH6-29}
\end{align}
where $\m{X}_1^{K}=\frac{1}{K}\sum_{k=1}^{K}\tilde{\m{x}}_1^k$, $\m{X}_2^{K}=\frac{1}{K}\sum_{k=1}^{K}\tilde{\m{x}}_2^k$, $\omega_1=2\tau\sigma\bar{\varepsilon}\sqrt{2m}$, $\omega_2=\|\M{x}^0-\M{x}^*\|_{\Gamma^{-1}}^2+\rho^2\|\M{Q}\|$, $J^*$ is defined in \eqref{Problem}, and $\rho=\max\{1+\|\M{y}^*\|,2\|\M{y}^*\|\}$.
\end{thm}
\begin{IEEEproof}
See Appendix \ref{APPGG}.
\end{IEEEproof}
\begin{rmk}
We assume the existence of the dual solution $\M{y}^*$, when assessing the primal suboptimality and consensus violation of V-DISA \eqref{Inexact}.
In terms of \cite[Section 4]{Hamedani2021}, if a slater point for \eqref{Trans-P1} is available, a bound of $\M{y}^*$ can be computed. In \cite{Hamedani2021}, an accelerated primal-dual (APD) algorithm has been proposed to solve bilinear saddle point problems and advanced technique, backtracking, has been provided. When APD is applied to convex optimization problems with nonlinear functional constraints, using this backtracking scheme, the optimal convergence rate can be achieved even when the dual domain is unbounded.
\end{rmk}

\begin{table*}[!t]
\renewcommand\arraystretch{1.3}
\begin{center}
\caption{Comparison of Distributed generalized Lasso problem solved by DISA, L-ALM, L-ADMM, Condat-Vu, TriPD-Dist, and TPUS}
\scalebox{0.86}{
\begin{tabular}{|c|c|cc|cc|cc|cc|cc|cc|}
\hline
\multirow{2}{*}{$n$}    & \multirow{2}{*}{$\|\M{UU}\tr\|$} & \multicolumn{2}{c|}{\textbf{DISA}} & \multicolumn{2}{c|}{\textbf{L-ALM}} & \multicolumn{2}{c|}{\textbf{L-ADMM}} & \multicolumn{2}{c|}{\textbf{Condat-Vu}} & \multicolumn{2}{c|}{\textbf{TriPD-Dist}} & \multicolumn{2}{c|}{\textbf{TPUS}} \\
\cline{3-14}
                      &                    & Iter. &Time(s)& Iter. &Time(s)& Iter. & Time(s)& Iter. & Time(s)&Iter. & Time(s)&Iter. & Time(s)\\
\hline
\hline
\multirow{5}{*}{200}  &3.4408&892&0.5942&975&0.9432&984&0.4456&973&1.0761&973&1.8929&968&3.7028\\
                      &331.9644&1576&0.9658&5346&5.6624&5113&2.3085&5153&5.8027&5175&9.4483 &5399&21.5687\\
                      &3.7126e+04&\cellcolor{myGray}\textbf{1315}&\cellcolor{myGray}\textbf{0.8465}&68728&71.6233&68876&29.9561&68931&85.7990&69654&122.3521&66817&251.6754\\
                      &3.3495e+06&\cellcolor{myGray}\textbf{1432}&\cellcolor{myGray}\textbf{0.9567}&719194&699.6877&711622&327.4182& 698351& 765.3428&737447&1.3558e+03 &720682 &2.6624e+03 \\
                      &3.4853e+08&\cellcolor{myGray}\textbf{1278}&\cellcolor{myGray}\textbf{0.8266}&$>$1e+06&$>$1e+03&$>$1e+06&$>$1e+03&$>$1e+06& $>$1e+03&$>$1e+06&$>$1e+03&$>$1e+06&$>$1e+03\\
\hline
\multirow{5}{*}{500}  &6.8988&584&1.3637&1555&9.4159&1557&7.9443&1532&12.9222&1560&22.7261&1561&40.5945\\
                      &466.0735&773&1.8274&4035&24.2059&4058&20.5233 & 3966& 33.3809&4211 & 62.5071&4024 &105.1999 \\
                      &2.5443e+04&\cellcolor{myGray}\textbf{770}&\cellcolor{myGray}\textbf{1.7541}& 91329& 540.0073& 99568& 491.3148&99584&786.3152& 99299&1.3500e+03& 99290&2.4566e+03 \\
                      &7.3258e+06&\cellcolor{myGray}\textbf{695}&\cellcolor{myGray}\textbf{1.7532}&$>$1e+05&$>$1e+03&$>$1e+05&$>$1e+03&$>$1e+05&$>$1e+03&$>$1e+05&$>$1e+03&$>$1e+05&$>$1e+03\\
                      &7.1088e+08&\cellcolor{myGray}\textbf{747}&\cellcolor{myGray}\textbf{1.8738}&$>$1e+05&$>$1e+03&$>$1e+05&$>$1e+03&$>$1e+05&$>$1e+03&$>$1e+05&$>$1e+03&$>$1e+05&$>$1e+03\\
\hline
\multirow{5}{*}{1000} &12.8915&572&4.2099& 2574&60.2779 &2570 &48.6809 &2579 & 86.6354&2575 &150.3896 &2569&250.5922\\
                      &322.2686&642&4.8546&3447&81.5100&3496&66.5924&3460&116.8618&3490&205.2825&3510 &331.3084 \\
                      &3.2946e+04&\cellcolor{myGray}\textbf{665}&\cellcolor{myGray}\textbf{4.9928}&152433&3.6276e+03&152348&2.8768e+03&152485&5.1073e+03&152405&8.8267e+03&152246&1.4327e+04\\
                      &3.2683e+06&\cellcolor{myGray}\textbf{645}&\cellcolor{myGray}\textbf{4.8568}&$>$1e+05&$>$3e+03&$>$1e+05&$>$3e+03&$>$1e+05&$>$5e+03&$>$1e+05&$>$8e+03&$>$1e+05&$>$1e+04\\
                      &3.1978e+08&\cellcolor{myGray}\textbf{651}&\cellcolor{myGray}\textbf{4.8740}&$>$1e+05&$>$3e+03&$>$1e+05&$>$3e+03&$>$1e+05&$>$5e+03&$>$1e+05&$>$8e+03&$>$1e+05&$>$1e+04\\
\hline
\end{tabular}}
\label{Table-com}
\end{center}
\end{table*}
\section{Numerical Simulation}\label{Simulation}
In this section, we conduct two numerical examples to validate the obtained theoretical results. All the algorithms are written in Matlab R2020b and implemented in a computer with 3.30 GHz AMD Ryzen 9 5900HS with Radeon Graphics and 16 GB memory.

\subsection{Distributed Generalized LASSO Problem}
\begin{figure}[!h]
\centering
\setlength{\abovecaptionskip}{-2pt}
\subfigure{
\includegraphics[width=0.31\linewidth]{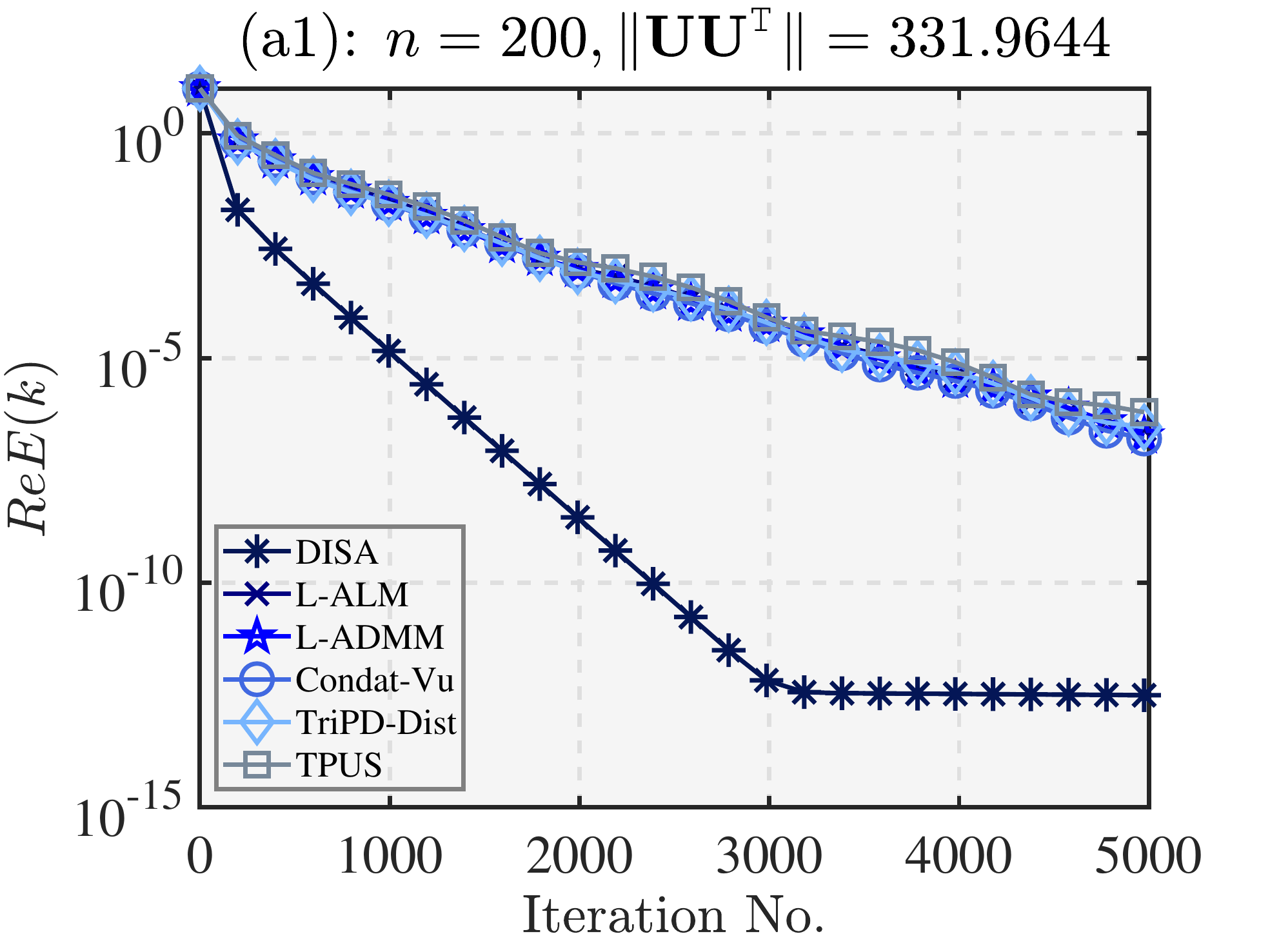}}
\subfigure{
\includegraphics[width=0.31\linewidth]{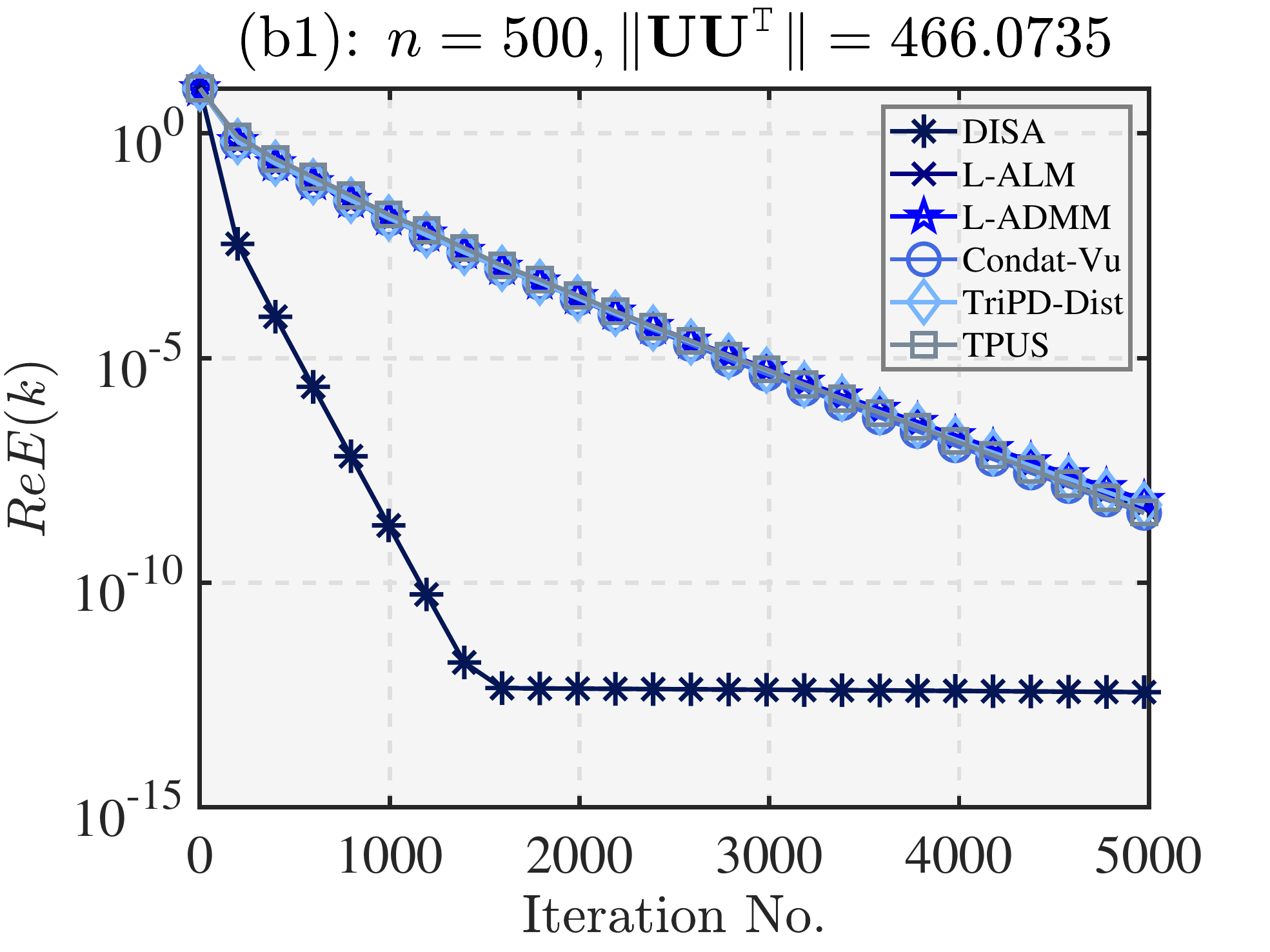}}
\subfigure{
\includegraphics[width=0.31\linewidth]{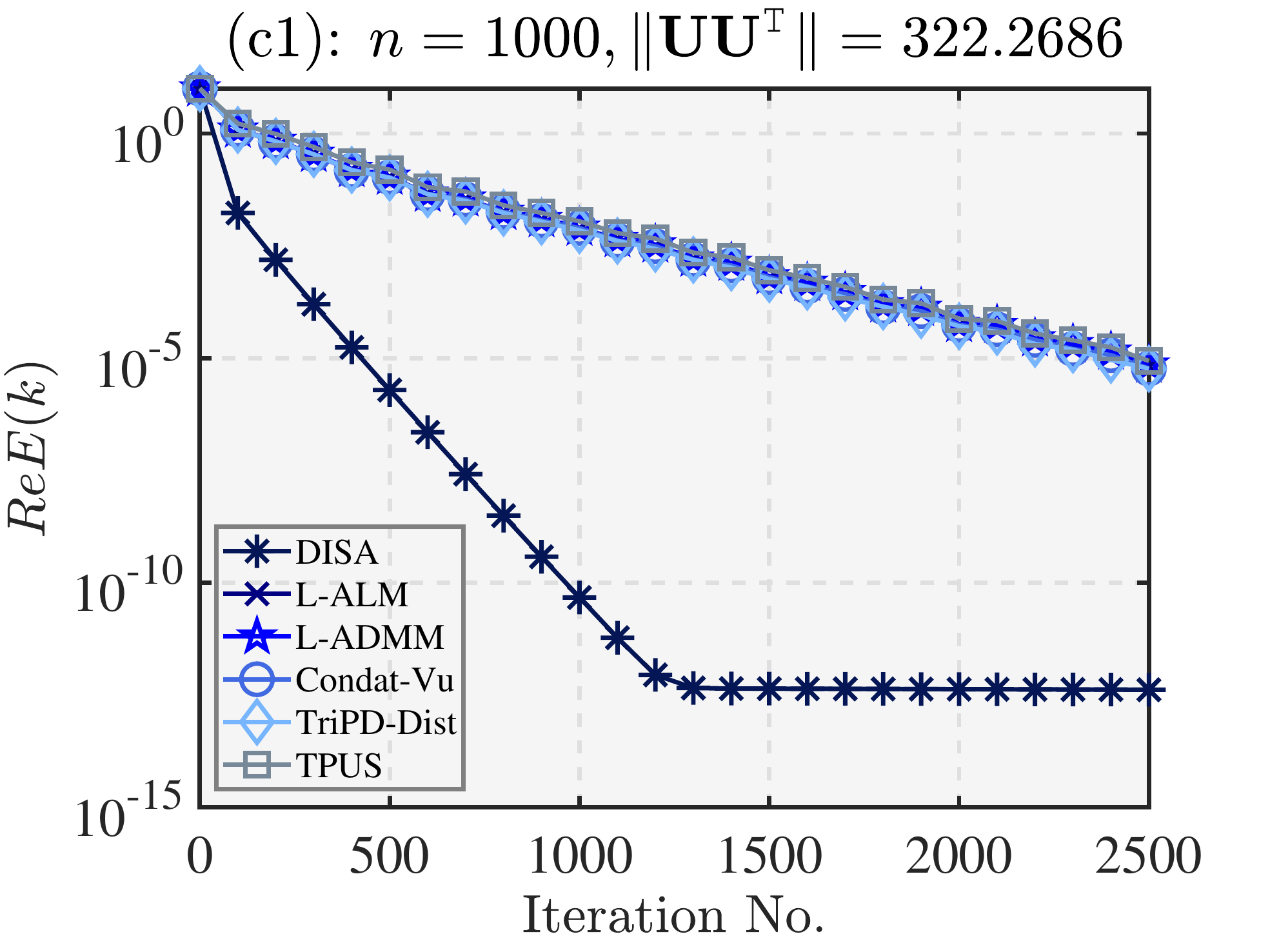}}
\subfigure{
\includegraphics[width=0.31\linewidth]{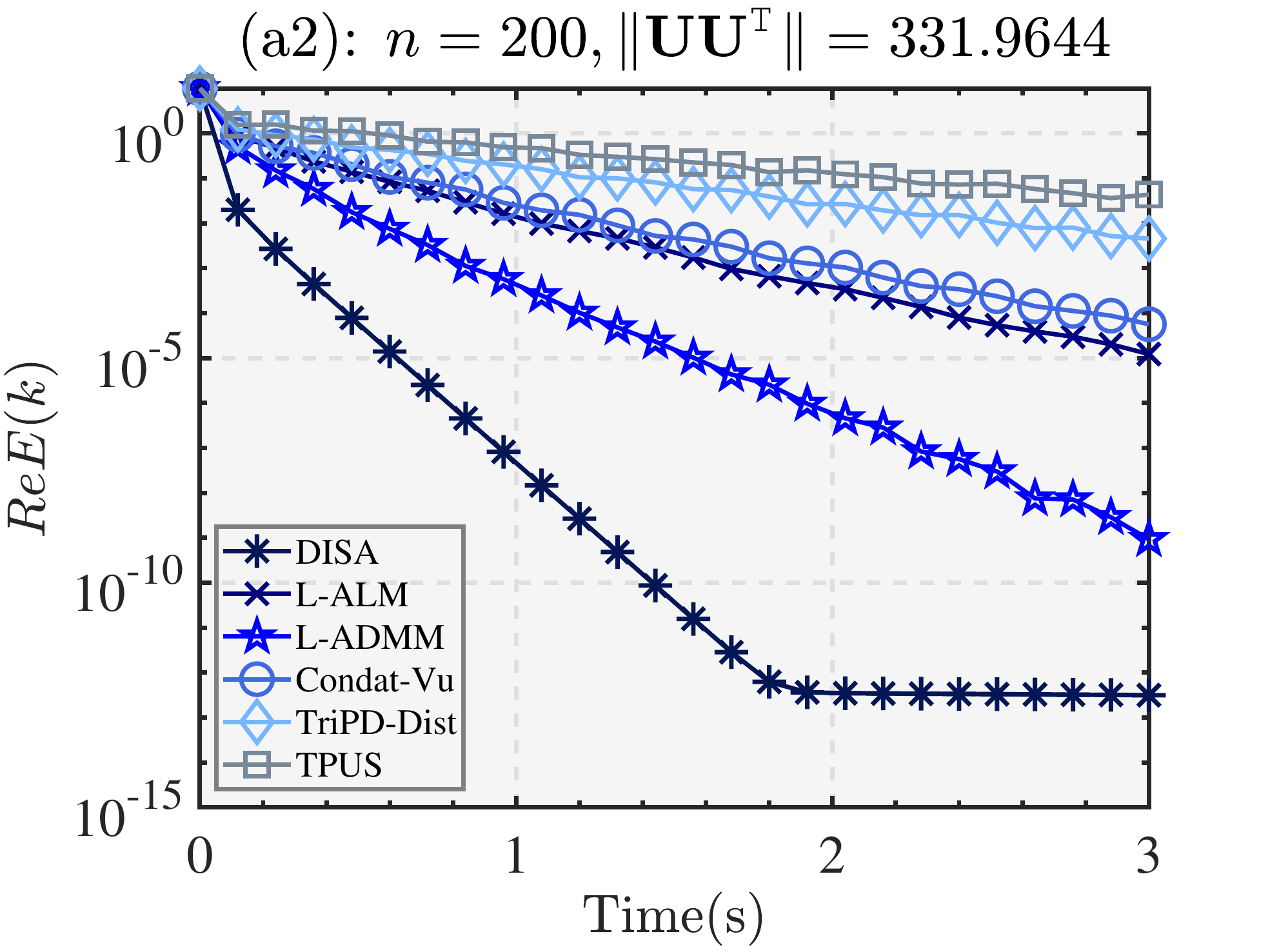}}
\subfigure{
\includegraphics[width=0.31\linewidth]{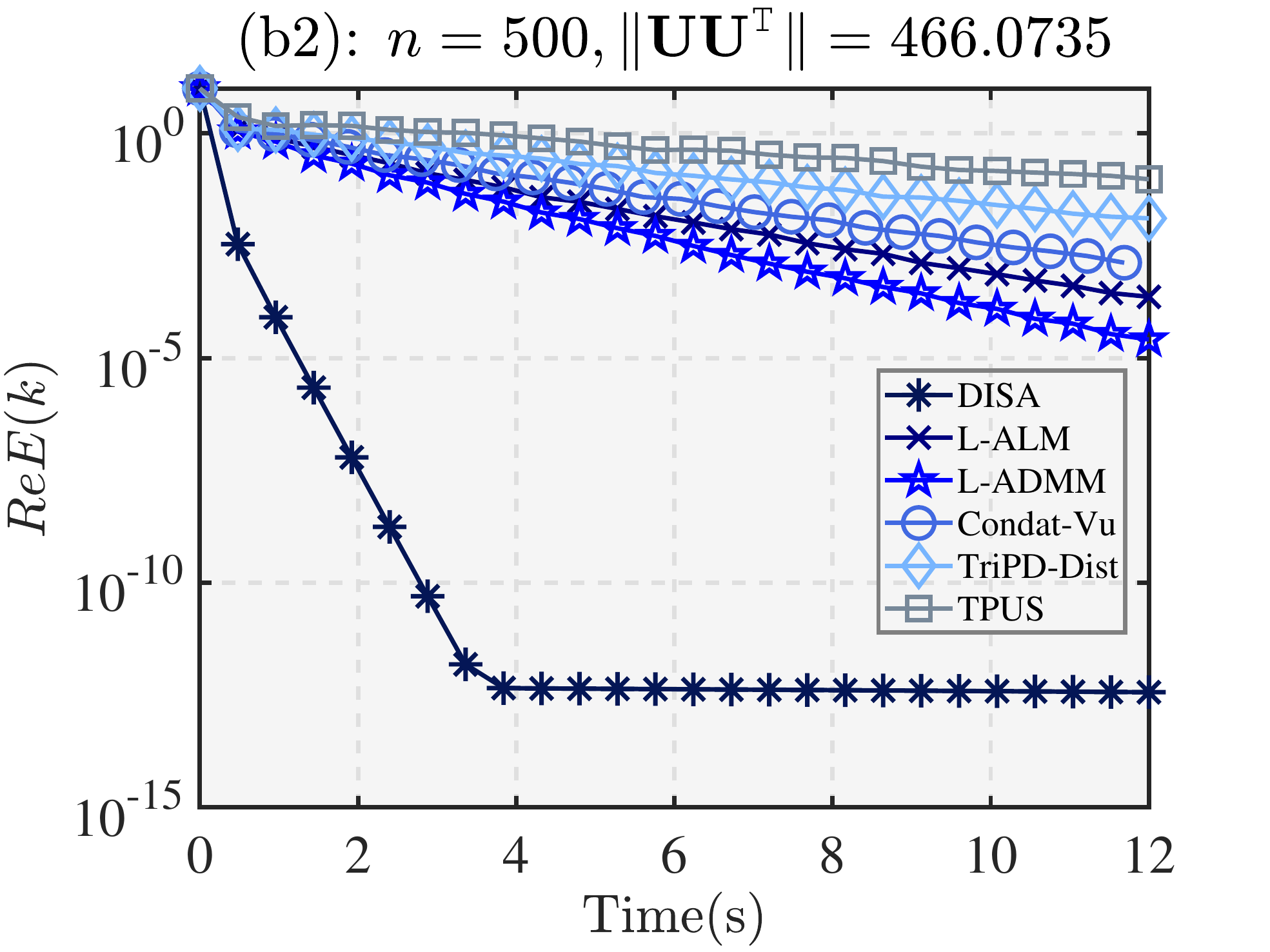}}
\subfigure{
\includegraphics[width=0.31\linewidth]{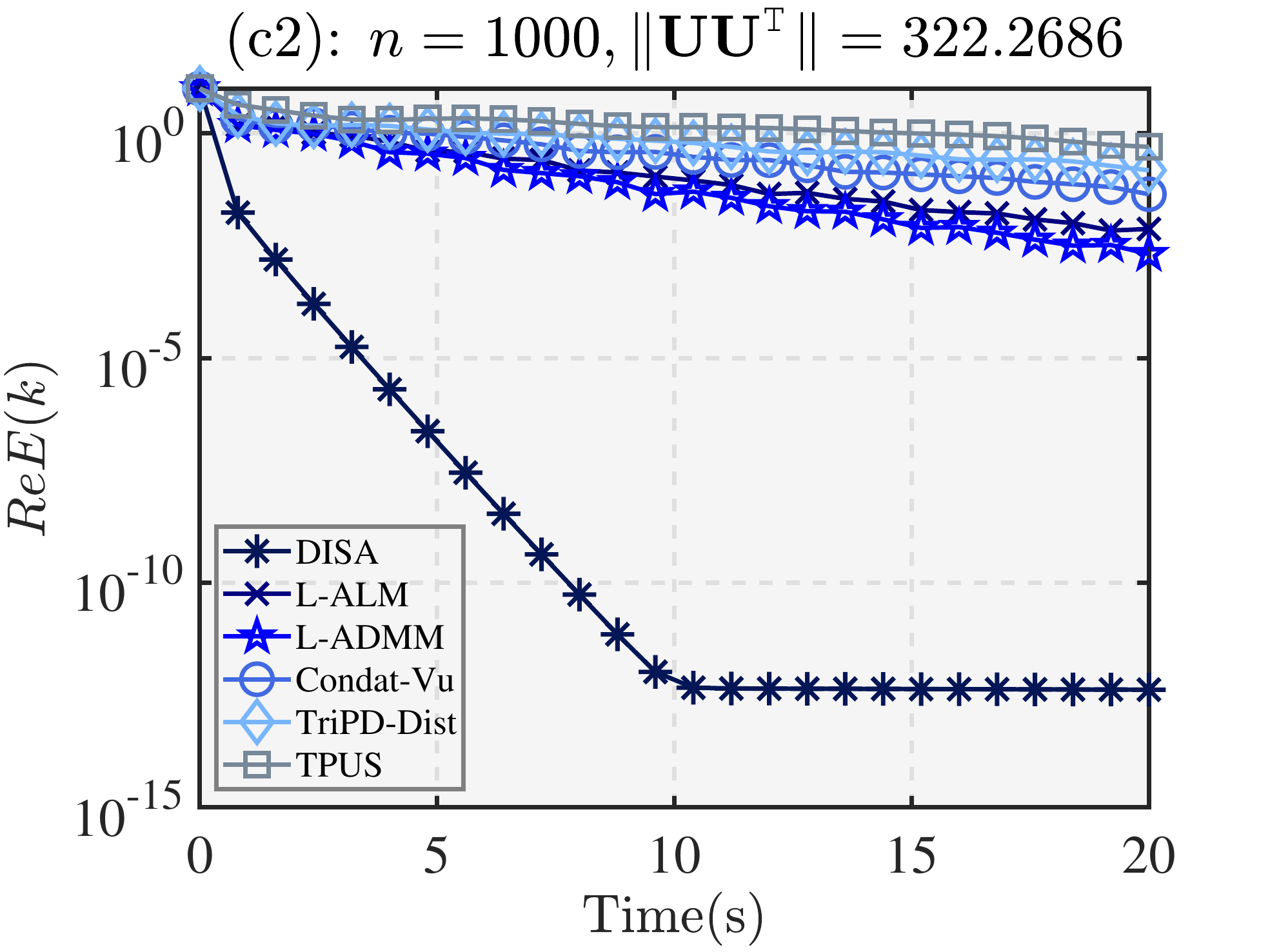}}
\caption{Convergence curves for distributed generalized Lasso solved by DISA, L-ALM, L-ADMM, Condat-Vu, TriPD-Dist, and TPUS.}
\label{case1}
\end{figure}

Consider the following decentralized generalized Lasso problem: $\sum_{i=1}^{m}\big\{\frac{1}{2}\|Q_i\m{x}-q_i\|^2+||U_i\m{x}||_1\big\},$
where each element in $Q_i\in\mathbb{R}^{2n\times n}$, $q_i\in\mathbb{R}^{2n}$, and $U_i\in\mathbb{R}^{20\times n}$ is drawn from the normal distribution. The communication topology is a line graph with 4 nodes. In this experiment, we solve the considered problem using DISA, L-ALM, L-ADMM, Condat-Vu, TriPD-Dist, and TPUS. Set the stopping criterion: $ReE(k):={\|\m{x}^k_1-\m{x}_1^*\|}/{\|\m{x}_1^*\|}<10^{-7}$.
To implement the aforementioned algorithms efficiently, we take the specific parameter settings as following. For DISA, we set $\tau_i={2}/{\|Q_i\tr Q_i\|}-0.0001$, and $\tau\beta=\frac{1}{2}$, where $\tau=\max_i\{\tau_i\}$. For L-ALM, L-ADMM, Condat-Vu, TriPD-Dist, and TPUS, we set $\tau = \min \{1/(\|Q_i\tr Q_i\|/2+\beta \|\M{U}\M{U}\tr\|)\}-0.0001$, where $\beta=0.5,0.01,0.001$ when $n=200,500,1000$, respectively. In Table \ref{Table-com}, for various values of $n$ and $\|\M{UU}\tr\|$, the required iteration number (Iter.) and the total computing time in seconds are presented. The associated convergence curves on some examples are plotted in Fig. \ref{case1}.

As shown in Table \ref{Table-com}, DISA has a significant acceleration compared with these existing PD-PSAs when $\|\M{UU}\tr\|$ is large. Moreover, $\|\M{UU}\tr\|$ does not affect the choice of step-sizes for DISA, which is consistent with our theoretical results. To further visualize the numerical results, in Fig. \ref{case1}, we plot the convergence curves versus both iteration numbers and run-times for the cases where $n=200$, $\|\M{U}\M{U}\tr\|=331.9644$, $n=500$, $\|\M{U}\M{U}\tr\|=466.0735$, and $n=1000$, $\|\M{U}\M{U}\tr\|=322.2686$, which can further demonstrate the numerical efficiency of DISA. In addition, we give Table \ref{TAB4} to present the required iteration number for V-DISA \eqref{Inexact} to solve the considered problem, when $n=500$. As shown in Table \ref{TAB4}, when $\varepsilon^k=\frac{1}{k}$,
the approximate iterative version of DISA is not convergent due to $\sum_{k=0}^{\infty}\varepsilon^k=\infty$. When $\varepsilon^k=\frac{1}{k^2},\frac{1}{k^3},e^{-k}$, the approximate iterative version are convergent. Moreover, it also has a significant acceleration and the choice of step-sizes is independent of $\|\M{U}\M{U}\tr\|$.
\begin{table}[!h]
\renewcommand\arraystretch{1.5}
\begin{center}
\caption{The required iteration number of V-DISA \eqref{Inexact}: $n=500$}
\scalebox{1}{
\begin{tabular}{ccccc}
\hline
$\|\M{UU}\tr\|$&$\varepsilon^k=\frac{1}{k}$&$\varepsilon^k=\frac{1}{k^2}$&$\varepsilon^k=\frac{1}{k^3}$&$\varepsilon^k=e^{-k}$\\
\hline
6.8988&N/A&664&583&577\\
466.0735&N/A&787&774&772\\
2.5443e+04&N/A&787&768&768\\
7.3258e+06&N/A&804&792&790\\
7.1088e+08&N/A&764&748&749\\
\hline
\end{tabular}}
\label{TAB4}
\end{center}
\end{table}

\begin{figure}[!h]
  \centering
  \setlength{\abovecaptionskip}{-2pt}
  \includegraphics[width=1\linewidth]{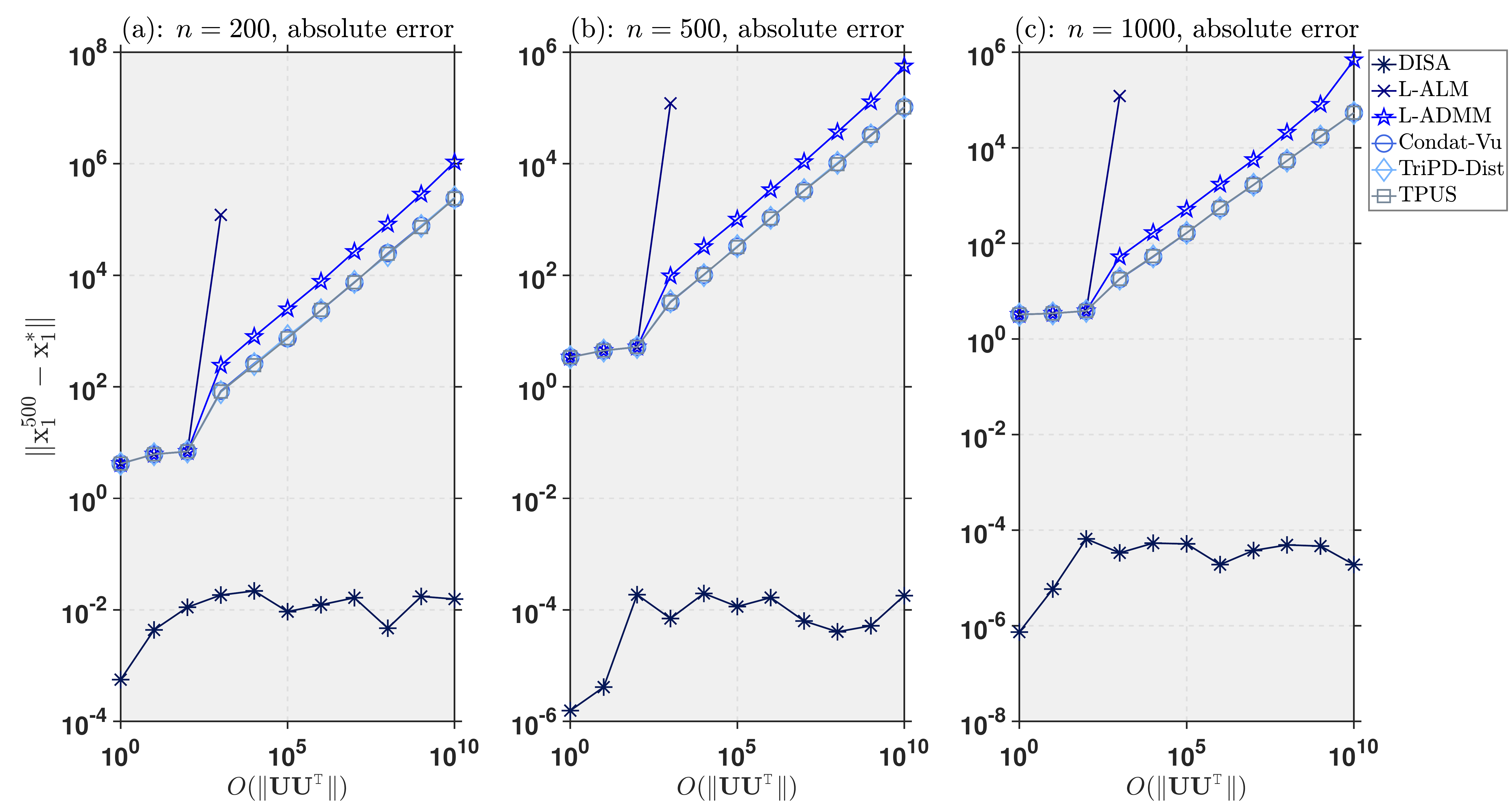}
  \caption{The absolute error $\|\m{x}_1^{500}-\m{x}_1^*\|$ against $O(\|\M{UU}\tr\|)$.}
  \label{CASE1D}
\end{figure}

Next, we investigate the impact of $\|\M{UU}\tr\|$ on the performance of various algorithms. To this end, we fix $\tau=\min_i{{1}/{\|Q_i\tr Q_i\|}}-0.0001$ and $\tau\beta=0.01$, for all algorithms. We then vary the values of $\|\M{UU}\tr\|$ from $O(10^0)$ to $O(10^{10})$, and evaluate the absolute error $\|\m{x}_1^k-\m{x}_1^*\|$ at iteration $k=500$ for each algorithm. Fig. \ref{CASE1D} shows that when $\|\M{UU}\tr\|>10^2$, the absolute error of L-ALM, L-ADMM, Condat-Vu, TriPD-Dist, and TPUS at iteration 500 increases with increasing $\|\M{UU}\tr\|$. This is because the condition $\tau\beta\|\M{U}\tr\M{U}\|>1$ required for their convergence (Table \ref{COMT}) is violated, leading to divergence. In contrast, DISA is not affected by the value of $\|\M{UU}\tr\|$.
\subsection{Distributed Logistic Regression on Real Datasets}
In this subsection, we consider the following distributed logistic regression problem
$$
\min _{\m{x}}\sum_{i=1}^{m}\big\{\frac{1}{m_i}\sum_{j=1}^{m_{i}} \ln \big(1+e^{-(\mathcal{U}_{i j}^{\top} \m{x})\mathcal{V}_{i j}}\big)+\frac{1}{2}\|\m{x}\|^2+\frac{1}{2}\|U_i\m{x}\|\big\},
$$
where $\m{x}\in \mathbb{R}^{n}$, $U_i\in \mathbb{R}^{20\times n}$ and each element in $U_i$ is drawn from the normal distribution $N(0,1)$. Any agent $i$ holds its own training date $\left(\mathcal{U}_{i j}, \mathcal{V}_{i j}\right) \in$ $\mathbb{R}^{n} \times\{-1,1\}, j=1, \cdots, m_{i}$, including sample vectors $\mathcal{U}_{i j}$ and corresponding classes $\mathcal{V}_{i j}$. It holds that the convex regularizer $\|\m{x}\|$ is nonsmooth, and $\mathrm{prox}_{\lambda\|\cdot\|}(\m{x})=\big(1-\frac{\lambda}{\max\{\|\m{x}\|,\lambda\}}\big)\m{x}$, where $\lambda>0$. Consider a circular
graph with 10 nodes, i.e., the agents form a cycle. We use four real datasets including a6a, a9a, covtype, and ijcnn1\footnote[1]{\href{https://www.csie.ntu.edu.tw/~cjlin/libsvmtools/datasets/}{https://www.csie.ntu.edu.tw/~cjlin/libsvmtools/datasets/}}, whose attributes are
$n=123$ and $\sum_{i=1}^{N}m_i=11220$, $n=123$ and $\sum_{i=1}^{N}m_i=32550$, $n=22$ and $\sum_{i=1}^{N}m_i=49990$, and $n=54$ and $\sum_{i=1}^{N}m_i=55500$, respectively. Moreover, the training samples are randomly and evenly distributed over all the $m$ agents.

\begin{figure}[!h]
\centering
\setlength{\abovecaptionskip}{-2pt}
\subfigure{
\includegraphics[width=0.455\linewidth]{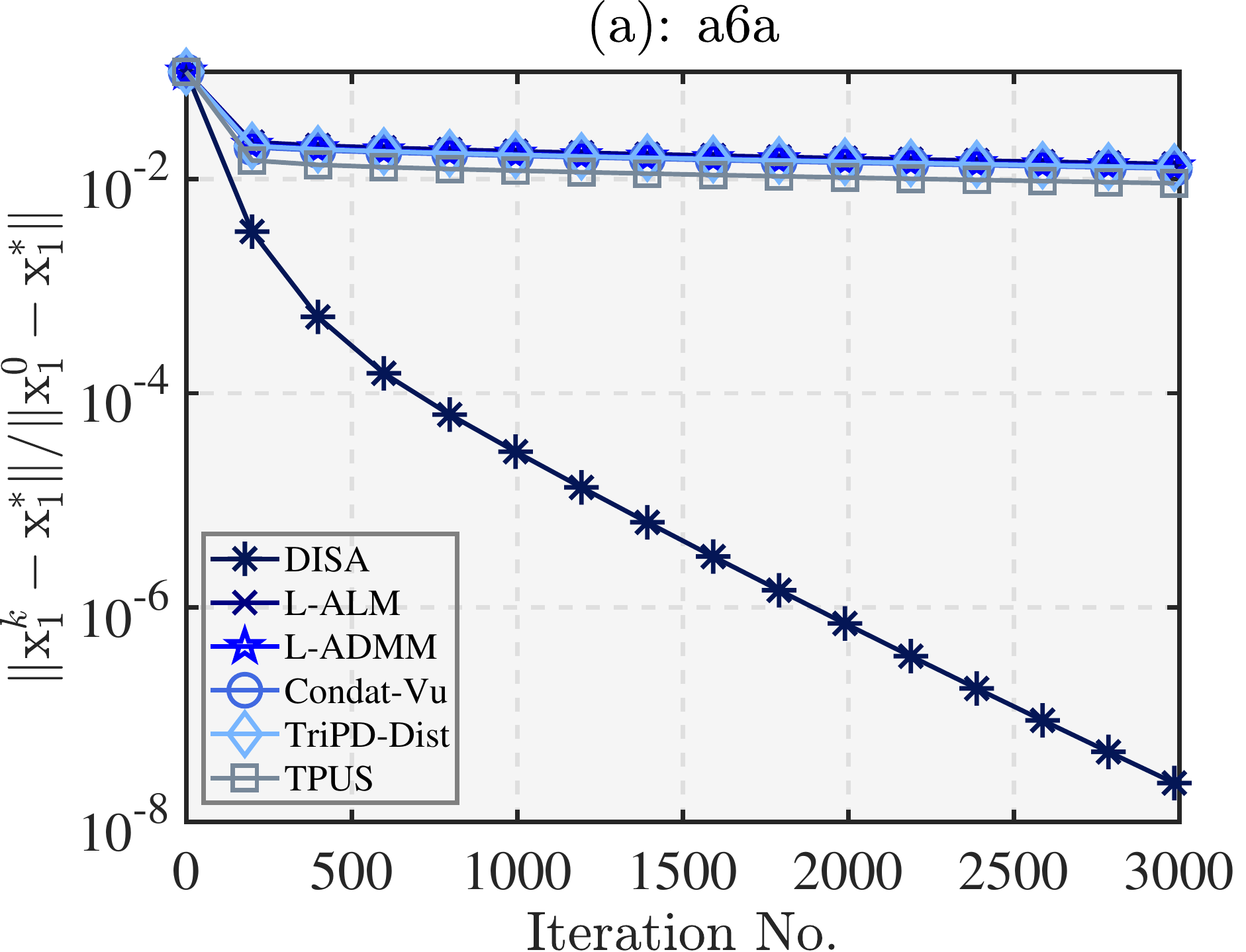}}
\subfigure{
\includegraphics[width=0.455\linewidth]{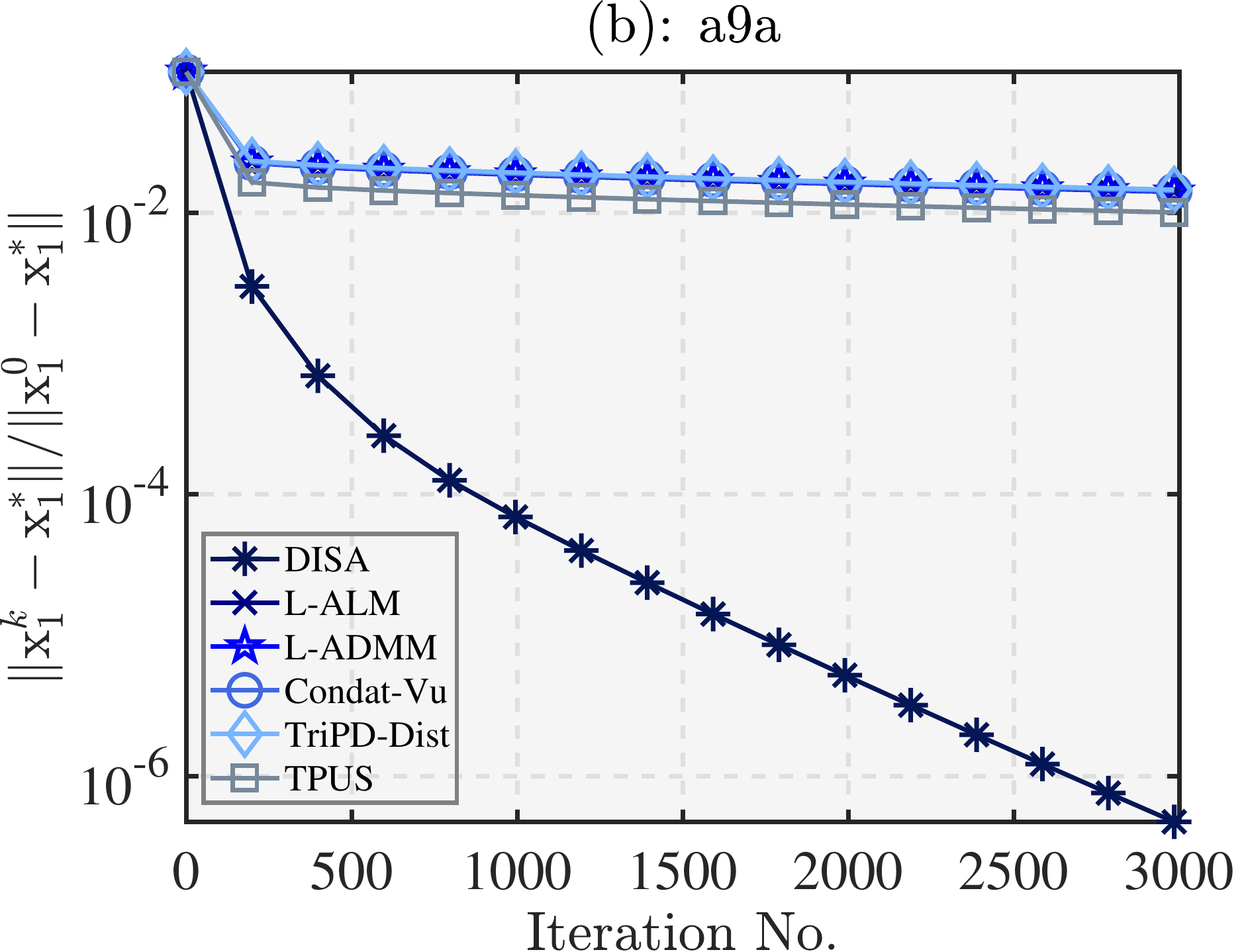}}
\subfigure{
\includegraphics[width=0.455\linewidth]{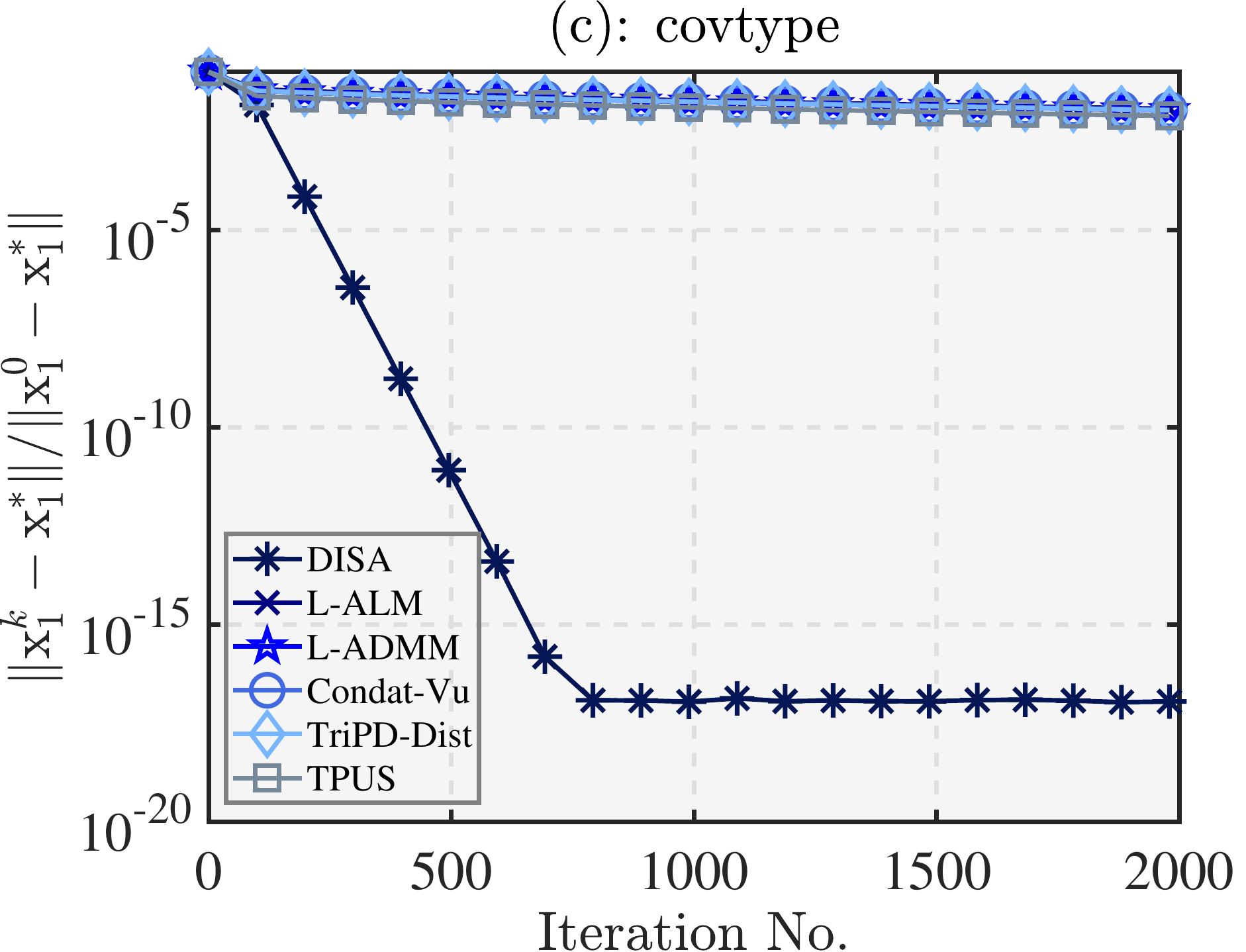}}
\subfigure{
\includegraphics[width=0.455\linewidth]{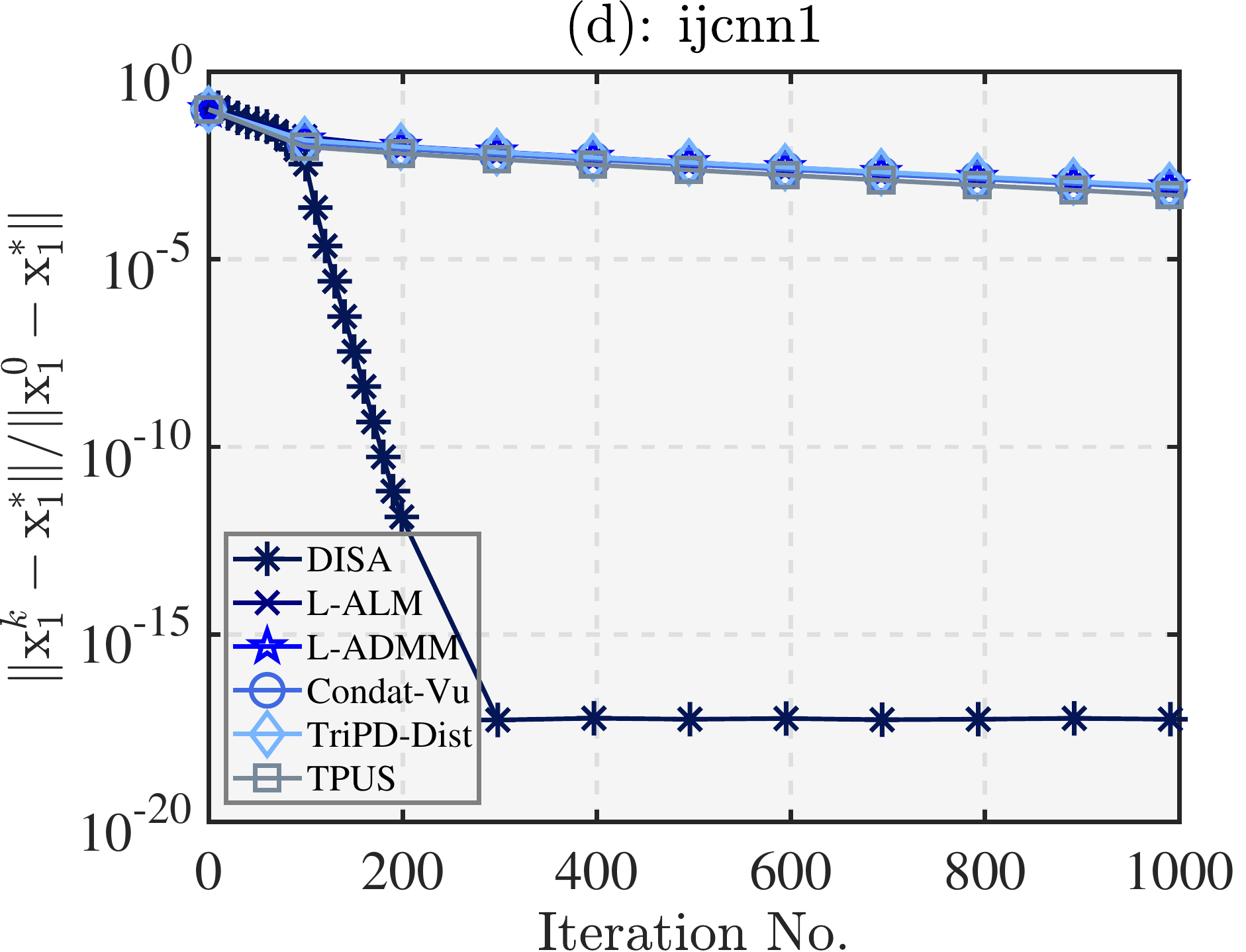}}
\caption{The results for distributed logistic regression over real datasets.}
\label{case2}
\end{figure}

In this experiment, we solve the distributed Logistic regression on real datasets still by DISA, L-ALM, L-ADMM, Condat-Vu, TriPD-Dist, and TPUS. To DISA, we set $\tau=0.25$ and $\tau\beta=0.5$. To L-ALM, L-ADMM, Condat-Vu, TriPD-Dist, and TPUS, we set $\tau=0.25$ and $\beta=0.01$. The performance is evaluated by the relative errors $\|\m{x}_1^k-\m{x}_1^*\|/\|\m{x}_1^0-\m{x}_1^*\|$. Given that the inverse of $S_i$ can be easily obtained (since $S_i\in\mathbb{R}^{20\times 20}$) and only needs to be computed once at the beginning of DISA's execution, the additional computational cost can be disregarded. Consequently, we report only the number of iterations. The results are illustrated in Fig. \ref{case2}. It is clear that DISA outperform L-ALM, L-ADMM, Condat-Vu, TriPD-Dist, and TPUS in all the four datasets.

\section{Conclusion}\label{Conclusion}
In this paper, we proposed DISA, an efficient algorithm for solving the COP \eqref{Problem} on a network. Our algorithm outperforms state-of-the-art PD-PSAs by having a convergent step-size range that is independent of the network topology and $\|U_i U_i\tr\|$, allowing it to converge quickly even for COPs with large $\|U_i U_i\tr\|$. We proved that DISA is sublinearly convergent under general convexity and linearly convergent under metric subregularity. Additionally, we proposed V-DISA with an approximate proximal mapping and established its convergence and sublinear convergence rate under a summable absolute error criterion and general convexity. Our numerical results demonstrated the advantages of DISA for distributed COP with a large Euclidean norm of the linear operator compared to existing algorithms.

We leave acceleration, asynchronous, and stochastic settings of DISA as future work. Moreover, the extension of uncoordinated dual step-sizes in DISA to circumvent a priori global coordination remains an open question.

\appendices

\section{Proof of Lemma \ref{NEWLEMMA1}}\label{APP45}
\begin{IEEEproof}
Since $f_i$ is convex and $L_i$-smooth, it holds that
\begin{align*}
&\langle \M{x}-\M{y},\nabla F(\M{z})-\nabla F(\M{x})\rangle\\
&=\langle \M{z}-\M{y}, \nabla F(\M{z})-\nabla F(\M{x})\rangle-\big\langle \M{z}-\M{x}, \nabla F(\M{z})-\nabla F(\M{x})\rangle \\
&\leq\langle \M{z}-\M{y}, \nabla F(\M{z})-\nabla F(\M{x})\rangle-\|\nabla F(\M{z})-\nabla F(\M{x})\|^2_{(\M{L}_F)^{-1}}\\
&=-\big\|\frac{(\M{L}_F)^{\frac{1}{2}}}{2}(\M{y}-\M{z})+(\M{L}_F)^{-\frac{1}{2}}(\nabla F(\M{z})-\nabla F(\M{x}))\big\|\\
&\quad+\frac{1}{4}\|\M{y}-\M{z}\|^2_{\M{L}_F}\leq\frac{1}{4}\|\M{y}-\M{z}\|^2_{\M{L}_F}.
\end{align*}
Hence, \eqref{SM1} holds. Additionally, \eqref{SM2} follows from
\begin{align*}
&\langle \M{x}-\M{y},\nabla F(\M{z})\rangle=\langle \M{x}-\M{z},\nabla F(\M{z})\rangle + \langle \M{z}-\M{y}, \nabla F(\M{z})\rangle\\
&\leq F(\M{x})-F(\M{z})+F(\M{z})-F(\M{y})+\frac{1}{2}\|\M{z}-\M{y}\|^2_{\M{L}_F}\\
&=F(\M{x})-F(\M{y})+\frac{1}{2}\|\M{z}-\M{y}\|^2_{\M{L}_F}.
\end{align*}
Therefore, the Lemma \ref{NEWLEMMA1} holds.
\end{IEEEproof}

\section{Proof of Lemma \ref{F1}}\label{APPA}
\begin{IEEEproof}
It follows from \eqref{barPrimalUP} and Lemma \ref{LEMF} that
\begin{align}\label{Proof-Lemma1-1}
&G(\M{x})-G(\bar{\M{x}}^{k+1})+\big\langle \M{x}-\bar{\M{x}}^{k+1},{\M{B}}\tr \M{y}^k+\nabla F(\M{x}^k) \nonumber \\
&\quad +\Gamma^{-1}(\bar{\M{x}}^{k+1}-\M{x}^k)\big\rangle\geq0, \forall \M{x}\in \mathbb{R}^{m(n+p)}.
\end{align}
Similar as \eqref{Proof-Lemma1-1}, from \eqref{PrimalUP}, one has
\begin{align}\label{Proof-Lemma1-2}
&G(\M{x})-G({\M{x}}^{k+1})+\big\langle \M{x}-{\M{x}}^{k+1},\M{B}\tr \M{y}^{k+1}+\nabla F(\M{x}^k) \nonumber \\
&\quad +\Gamma^{-1}(\M{x}^{k+1}-\M{x}^k)\big\rangle\geq0, \forall \M{x}\in \mathbb{R}^{m(n+p)}.
\end{align}
Setting $\M{x}=\M{x}^{k+1}$ in \eqref{Proof-Lemma1-1} and then adding to inequality \eqref{Proof-Lemma1-2}, one has that for $\forall \M{x}\in \mathbb{R}^{m(n+p)}$
\begin{align}\label{Proof-Lemma1-3}
&\big\langle \M{x}-\M{x}^{k+1},\M{B}\tr\M{y}^{k+1}\big\rangle+\big\langle \M{x}^{k+1}-\bar{\M{x}}^{k+1},\M{B}\tr \M{y}^k\big\rangle\nonumber\\
&-\big\langle \M{x}-\bar{\M{x}}^{k+1},\M{B}\tr \M{y}\big\rangle+\big\langle \M{x}-\bar{\M{x}}^{k+1} ,\nabla F(\M{x}^k)\big\rangle  \nonumber\\
&+\big\langle \M{x}^{k+1}-\bar{\M{x}}^{k+1},\Gamma^{-1}(\bar{\M{x}}^{k+1}-\M{x}^k)\big\rangle\nonumber\\
&+\big\langle \M{x}-\M{x}^{k+1},\Gamma^{-1}(\M{x}^{k+1}-\M{x}^k)\big\rangle  \nonumber\\
&\geq G(\bar{\M{x}}^{k+1})-G(\M{x})+\big\langle \bar{\M{x}}^{k+1}-\M{x}, \M{B}\tr\M{y}\big\rangle.
\end{align}
By \eqref{DualUP} and Lemma \ref{LEMF}, it holds that
\begin{align}\label{Proof-Lemma1-4}
&\big\langle \M{y}-\M{y}^{k+1},-\M{B}(\bar{\M{x}}^{k+1}-\M{x})+\M{Q}(\M{y}^{k+1}-\M{y}^k) \big\rangle \nonumber\\
&\quad \geq \big\langle \M{y}^{k+1}-\M{y},-{\M{B}}\M{x} \big\rangle,\forall \M{y}\in \mathbb{R}^{m(n+p)}.
\end{align}
On one hand, we have
\begin{align}\label{PROOFLEM233}
&\big\langle \M{x}-\M{x}^{k+1},{\M{B}}\tr \M{y}^{k+1}\big\rangle+\big\langle \M{x}^{k+1}-\bar{\M{x}}^{k+1},{\M{B}}\tr \M{y}^k\big\rangle\nonumber\\
&-\big\langle \M{x}-\bar{\M{x}}^{k+1},{\M{B}}\tr \M{y}\big\rangle+\big\langle \M{y}-\M{y}^{k+1},-\M{B}(\bar{\M{x}}^{k+1}-\M{x}) \big\rangle\nonumber\\
&=\big\langle \bar{\M{x}}^{k+1}-\M{x}^{k+1},\M{B}\tr (\M{y}^{k+1}-\M{y}^k).
\end{align}
On the other hand, applying the identity $\langle \m{p}_1-\m{p}_2, \mathrm{H}(\m{q}_1-\m{q}_2)\rangle=\frac{1}{2}\{\|\m{p}_1-\m{q}_2\|^2_{\mathrm{H}}-\|\m{p}_1-\m{q}_1\|^2_{\mathrm{H}}\}+\frac{1}{2}\{\|\m{p}_2-\m{q}_1\|^2_{\mathrm{H}}-\|\m{p}_2-\m{q}_2\|^2_{\mathrm{H}}\}$,
one gets that for $\forall (\M{x},\M{y})\in \mathrm{M}$
\begin{subequations}\label{ZYBDS}
\begin{align}
&\big\langle \M{y}-\M{y}^{k+1},\M{Q}(\M{y}^{k+1}-\M{y}^k) \big\rangle=\frac{1}{2}\|\M{y}^k-\M{y}\|_{\M{Q}}^2\nonumber\\
&-\frac{1}{2}\big(\|\M{y}^{k+1}-\M{y}\|_{\M{Q}}^2+\|\M{y}^k-\M{y}^{k+1}\|_{\M{Q}}^2\big),\\
&\big\langle \M{x}^{k+1}-\bar{\M{x}}^{k+1},\Gamma^{-1}(\bar{\M{x}}^{k+1}-\M{x}^k)\big\rangle = \frac{1}{2}\|\M{x}^k-\M{x}^{k+1}\|_{\Gamma^{-1}}^2 \nonumber\\ &-\frac{1}{2}\big(\|\M{x}^k-\bar{\M{x}}^{k+1}\|_{\Gamma^{-1}}^2+\|\M{x}^{k+1}-\bar{\M{x}}^{k+1}\|_{\Gamma^{-1}}^2\big),\\
&\quad\big\langle \M{x}-\M{x}^{k+1},\Gamma^{-1}(\M{x}^{k+1}-\M{x}^k)\big\rangle=\frac{1}{2}\|\M{x}-\M{x}^k\|_{\Gamma^{-1}}^2\nonumber\\
&-\frac{1}{2}\big(\|\M{x}-\M{x}^{k+1}\|_{\Gamma^{-1}}^2+\|\M{x}^k-\M{x}^{k+1}\|_{\Gamma^{-1}}^2\big).
\end{align}
\end{subequations}
Hence, summing \eqref{Proof-Lemma1-3} and \eqref{Proof-Lemma1-4}  and then substituting \eqref{PROOFLEM233} and \eqref{ZYBDS} into it, we have that for $\forall (\M{x},\M{y})\in \mathcal{M}$
\begin{align}\label{Proof-Lemma1-6}
&\big(\frac{1}{2}\|\M{x}^k-\M{x}\|_{\Gamma^{-1}}^2+\frac{1}{2}\|\M{y}^k-\M{y}\|_{\M{Q}}^2\big)\nonumber\\
&-\big(\frac{1}{2}\|\M{x}^{k+1}-\M{x}\|_{\Gamma^{-1}}^2+\frac{1}{2}\|\M{y}^{k+1}-\M{y}\|_{\M{Q}}^2\big)\nonumber\\
&\geq \frac{1}{2}\|\M{x}^k-\bar{\M{x}}^{k+1}\|_{\Gamma^{-1}}^2+\big\langle \bar{\M{x}}^{k+1}-\M{x}^{k+1},\M{B}\tr (\M{y}^{k}-\M{y}^{k+1})\big\rangle \nonumber\\
&\quad+\frac{1}{2}\|\M{x}^{k+1}-\bar{\M{x}}^{k+1}\|_{\Gamma^{-1}}^2+\frac{1}{2}\|\M{y}^k-\M{y}^{k+1}\|_{\M{Q}}^2 \nonumber\\
&\quad +G(\bar{\M{x}}^{k+1})-G(\M{x})+\big\langle \bar{\M{x}}^{k+1}-\M{x}, {\M{B}}\tr\M{y}\big\rangle\nonumber\\
&\quad + \big\langle \M{y}^{k+1}-\M{y},-{\M{B}}\M{x} \big\rangle + \big\langle \M{\bar{x}}^{k+1}-\M{x},\nabla F(\M{x}^k) \big\rangle.
\end{align}
By Cauchy-Schwarz inequality, it holds that
\begin{align*}
&\big\langle \bar{\M{x}}^{k+1}-\M{x}^{k+1},\M{B}\tr (\M{y}^{k}-\M{y}^{k+1})\big\rangle \\
&\geq -\frac{1}{2}\|\bar{\M{x}}^{k+1}-\M{x}^{k+1}\|_{\Gamma^{-1}}^2-\frac{1}{2}\|\M{y}^{k}-\M{y}^{k+1}\|^2_{\M{B}\Gamma\M{B}\tr}.
\end{align*}
Consequently, we have that for $\forall (\M{x},\M{y})\in \mathcal{M}$
\begin{align}\label{NEWPROOF1}
&\big(\frac{1}{2}\|\M{x}^k-\M{x}\|_{\Gamma^{-1}}^2+\frac{1}{2}\|\M{y}^k-\M{y}\|_{\M{Q}}^2\big)\nonumber\\
&-\big(\frac{1}{2}\|\M{x}^{k+1}-\M{x}\|_{\Gamma^{-1}}^2+\frac{1}{2}\|\M{y}^{k+1}-\M{y}\|_{\M{Q}}^2\big)\nonumber\\
&\geq \frac{1}{2}\|\M{x}^k-\bar{\M{x}}^{k+1}\|_{\Gamma^{-1}}^2+ \frac{1}{2}\|\M{y}^k-\M{y}^{k+1}\|^2_{\M{Q}- \M{B}\Gamma\M{B}\tr}\nonumber\\
&\quad +G(\bar{\M{x}}^{k+1})-G(\M{x})+\big\langle \bar{\M{x}}^{k+1}-\M{x}, {\M{B}}\tr\M{y}\big\rangle\nonumber\\
&\quad + \big\langle \M{y}^{k+1}-\M{y},-{\M{B}}\M{x} \big\rangle + \big\langle \M{\bar{x}}^{k+1}-\M{x},\nabla F(\M{x}^k) \big\rangle.
\end{align}
By \eqref{SM1}, one has that $\big\langle \M{x}-\M{\bar{x}}^{k+1},\nabla F(\M{x}^k)-\nabla F(\M{x})\big\rangle\leq \frac{1}{4} \|\bar{\M{x}}^{k+1}-\M{x}^k\|_{\M{L}_F}^2$. Therefore, it holds that
\begin{align}\label{NEWPROOF2}
&\big\langle \M{\bar{x}}^{k+1}-\M{x},\nabla F(\M{x}^k) \big\rangle\nonumber\\
&=\big\langle \M{\bar{x}}^{k+1}-\M{x},\nabla F(\M{x}^k)- \nabla F(\M{x})\big\rangle+\big\langle \M{\bar{x}}^{k+1}-\M{x},\nabla F(\M{x})\big\rangle\nonumber\\
&\geq-\frac{1}{4} \|\bar{\M{x}}^{k+1}-\M{x}^k\|_{\M{L}_F}^2+\big\langle \M{\bar{x}}^{k+1}-\M{x},\nabla F(\M{x})\big\rangle.
\end{align}
Combining it with \eqref{NEWPROOF1} and \eqref{NEWPROOF2}, the inequality \eqref{BASCI-IEQ1} holds.

From \eqref{SM2}, it holds that $\langle \M{x}-\bar{\M{x}}^{k+1},\nabla F(\M{x}^k)\rangle\leq F(\M{x})-F(\bar{\M{x}}^{k+1})+\frac{1}{2}\|\M{x}^k-\bar{\M{x}}^{k+1}\|_{\M{L}_F}^2$. Note that
$(F+G)(\bar{\M{x}}^{k+1})-(F+G)(\M{x})+\big\langle \bar{\M{x}}^{k+1}-\M{x}, {\M{B}}\tr\M{y}\big\rangle+\big\langle \M{y}^{k+1}-\M{y},-{\M{B}}\M{x} \big\rangle=\mathcal{L}(\bar{\M{x}}^{k+1},\M{y})-\mathcal{L}(\M{x},\M{y}^{k+1})$. By \eqref{NEWPROOF1}, \eqref{BASCI-IEQ4} is proven.
\end{IEEEproof}

\section{Proof of Theorem \ref{THE-1}}\label{APPB}
\begin{IEEEproof}
Let $\M{w}^*=(\M{x}^*,\M{y}^*)\in \mathcal{M}^*$ be an arbitrary saddle point of the Lagrangian of problem \eqref{Trans-P1}. Summing the inequality \eqref{BASCISS} over $k=0,1,\cdots,K-1$ yields that
\begin{align*}
&\sum_{k=0}^{K-1}\|\M{w}^k-\M{v}^{k+1}\|^2_{\M{M}}\leq \sum_{k=0}^{K-1}(\|\M{w}^k-\M{w}\|_{\M{H}}^2-\|\M{w}^{k+1}-\M{w}\|_{\M{H}}^2)\\
&=\|\M{w}^0-\M{w}^*\|_{\M{H}}^2-\|\M{w}^{K}-\M{w}^*\|_{\M{H}}^2,\quad \forall K\geq1.
\end{align*}
Thus, it holds that for any $K\geq1$,
\begin{align}\label{THPPPFF1}
\sum_{k=0}^{K-1}\|\M{w}^k-\M{v}^{k+1}\|^2_{\M{M}}\leq \|\M{w}^0-\M{w}^*\|_{\M{H}}^2,
\end{align}
which implies that $\sum_{k=0}^{\infty}\|\M{w}^k-\M{v}^{k+1}\|^2_{\M{M}}<\infty$. Since $\M{M}$ is positive definite, we have $\lim_{k\rightarrow \infty}\|\M{w}^k-\M{v}^{k+1}\|=0.$

Let $\mathbf{w}^{\infty}$ be an accumulation point of $\{\mathbf{w}^k\}$ and $\{\mathbf{w}^{k_j}\}$ be a subsequence converging
to $\mathbf{w}^{\infty}$. By the nonexpansivity of $\mathrm{prox}^{\Gamma^{-1}}_{ G}(\cdot)$, one has that
\begin{align}\label{Rta}
\|\M{x}^{k+1}-\bar{\M{x}}^{k+1}\|\leq\tau\|\mathbf{B}\tr(\M{y}^k-\M{y}^{k+1})\|,
\end{align}
where $\tau=\max_i\{\tau_i\}$.
Thus, it holds that
$
\lim_{k\rightarrow \infty}\|\M{v}^{k+1}-\M{w}^{k+1}\|=\lim_{k\rightarrow \infty} \big(\|\bar{\M{x}}^{k+1}-\M{x}^{k+1}\|+\|\M{y}^{k+1}-\M{y}^{k+1}\|\big)=0,
$
which implies that $\mathbf{w}^{\infty}$ is also an accumulation point of $\{\mathbf{v}^k\}$ and $\{\mathbf{v}^{k_j}\}$ is a subsequence converging
to $\mathbf{w}^{\infty}$. Hence, it deduces that
$\lim_{k\rightarrow \infty}\|\M{w}^{k}-\M{w}^{k+1}\|
\leq \lim_{k\rightarrow \infty}\|\M{w}^{k}-\M{v}^{k+1}\| + \lim_{k\rightarrow \infty}\|\M{v}^{k+1}-\M{w}^{k+1}\|=0.
$
From \eqref{Proof-Lemma1-1} and \eqref{Proof-Lemma1-4}, it holds that for $\forall \M{w}\in \mathcal{M}$
\begin{align*}
&G(\M{x})-G(\bar{\M{x}}^{k_j+1})+\big\langle \M{x}-\bar{\M{x}}^{k_j+1},\M{B}\tr \M{y}^{k_j} \\
&\quad +\nabla F(\M{x}^{k_j}) +\Gamma^{-1}(\bar{\M{x}}^{k_j+1}-\M{x}^{k_j})\big\rangle\geq0 ,\\
&\big\langle \M{y}-\M{y}^{k_j+1},-\M{B}\bar{\M{x}}^{k_j+1}+\M{Q}(\M{y}^{k_j+1}-\M{y}^{k_j})\big\rangle\geq0.
\end{align*}
Taking $k_j\rightarrow\infty$ in the above two inequalities, it holds that
$$
\Theta(\M{x},\M{y}^{\infty})-\Theta(\M{x}^{\infty},\M{y})\geq 0, \forall \M{w}\in \mathcal{M}.
$$
Compared to \textbf{VI 2}, we can deduce that $\M{w}^{\infty} \in \mathcal{M}^*$. Setting $\M{w}=\M{w}^{\infty}$ in \eqref{BASCI-IEQ1} yields that $\|\M{w}^k-\M{w}^{\infty}\|_{\mathbf{H}}^2-\|\M{w}^{k+1}-\M{w}^{\infty}\|_{\mathbf{H}}^2\geq 0$.
It implies that the sequence $\{\|\M{w}^k-\M{w}^{\infty}\|_{\mathbf{H}}\}$ converges to a unique limit point. Then, with $\M{w}^{\infty}$ being an accumulation point of $\{\M{w}^k\}$, we have that $\|\M{w}^k-\M{w}^{\infty}\|^2_{\M{H}}\rightarrow0$, i.e., $\lim_{k\rightarrow \infty}\M{w}^k=\M{w}^{\infty}$.
\end{IEEEproof}

\section{Proof of Theorem \ref{THE-2}}\label{APPC}
\begin{IEEEproof}
By \eqref{THPPPFF1} and \cite[Proposition 1]{PGEXTR}, the non-ergodic rates \eqref{RA-RATE1} and \eqref{RB-RATE1} hold, immediately. From the updates of $\bar{\M{x}}^{k+1}$ and $\M{y}^{k+1}$, we have
\begin{align*}
&0\in \partial G(\bar{\M{x}}^{k+1})+\nabla F(\M{x}^k)+{\M{B}}\tr\M{y}^k+\Gamma^{-1}(\bar{\M{x}}^{k+1}-\M{x}^k),\\
&0=-{\M{B}}\bar{\M{x}}^{k+1}+\M{Q}(\M{y}^{k+1}-\M{y}^{k}),
\end{align*}
which implies that
\begin{align*}
&\nabla F(\bar{\M{x}}^{k+1})-\nabla F(\M{x}^{k})+{\M{B}}\tr(\M{y}^{k+1}-\M{y}^{k})-\Gamma^{-1}(\bar{\M{x}}^{k+1}-\M{x}^k)\\
&\in\partial G(\bar{\M{x}}^{k+1})+\nabla F(\bar{\M{x}}^{k+1})+{\M{B}}\tr\M{y}^{k+1},\\
&-\M{Q}(\M{y}^{k+1}-\M{y}^{k})=-{\M{B}}\bar{\M{x}}^{k+1}.
\end{align*}
Therefore, it holds that
\begin{align*}
&\mathrm{dist}^2(0,\mathcal{T}(\M{v}^{k+1}))\leq \|\nabla F(\bar{\M{x}}^{k+1})-\nabla F(\M{x}^k)\\
&-\Gamma^{-1}(\bar{\M{x}}^{k+1}-\M{x}^k)+{\M{B}}\tr(\M{y}^{k+1}-\M{y}^k)\|^2+\|\M{Q}(\M{y}^{k+1}-\M{y}^k)\|^2\\
&\leq(3L^2+\frac{3}{\tau})\|\M{x}^k-\bar{\M{x}}^{k+1}\|^2+(3\|\M{B}\|^2+\|\M{Q}\|^2)\|\M{y}^{k+1}-\M{y}^k\|^2\\
&\leq\kappa^2_1\|\M{v}^{k+1}-\M{w}^k\|^2_{\M{M}},
\end{align*}
where $\kappa^2_1=\frac{1}{c_1^2}\max\{3L^2+\frac{3}{\tau^2},3\|\M{B}\|^2+\|\M{Q}\|^2\}$. Therefore, by \cite[Proposition 1]{PGEXTR}, the conclusions \eqref{RA-RATE1} and \eqref{RB-RATE1} hold
\end{IEEEproof}

\section{Proof of Theorem \ref{THE-PDGAP}}\label{APP-PDGAP}
\begin{IEEEproof}
Summing the inequality \eqref{BASCI-IEQ4} over $k=0,1,\cdots,K-1$, we obtain
$
2\sum_{k=0}^{K-1}(\mathcal{L}(\bar{\M{x}}^{k+1},\M{y})-\mathcal{L}(\M{x},\M{y}^{k+1}))\leq\|\M{w}^0-\M{w}\|_{\M{H}}^2-\|\M{w}^K-\M{w}\|_{\M{H}}^2.
$
By the convexity of $F$, $G$ and the definition of $(\M{X}^{K},\M{Y}^K)$, we have $$K(\mathcal{L}(\M{X}^{K},\M{y})-\mathcal{L}(\M{x},\M{Y}^{K}))\leq\sum_{k=0}^{K-1}(\mathcal{L}(\bar{\M{x}}^{k+1},\M{y})-\mathcal{L}(\M{x},\M{y}^{k+1})).$$
Therefore, the primal-dual gap \eqref{PDGAP} holds.

Note that the inequality \eqref{PDGAP} holds for any $\M{w}\in \mathcal{M}$, hence it is also holds for $\{\M{x}^*\}\times \mathcal{B}_{\rho}$, where $\mathcal{B}_{\rho}=\{\M{y}:\|\M{y}\|\leq\rho\}$ with $\rho$ being any given positive number. Since the mapping $\mathcal{K}_1$ is affine with a skew-symmetric matrix, we have
$$
\langle \M{w}_1-\M{w}_2,\mathcal{K}_1(\M{w}_1)- \mathcal{K}_1(\M{w}_2)\rangle\equiv0, \forall \M{w}_1,\M{w}_2\in\mathcal{M}.
$$
Letting $\M{w}=(\M{w}^*,\M{y})$, and $\M{\Lambda}^K=(\M{X}^K,\M{Y}^K)$ in \eqref{PDGAP}, it gives
\begin{align}\label{GAP2}
&\sup_{\M{y}\in \mathcal{B}_{\rho}}\{\mathcal{L}(\M{X}^{K},\M{y})-\mathcal{L}(\M{x}^*,\M{Y}^{K})\}\nonumber\\
&=\sup_{\M{y}\in \mathcal{B}_{\rho}}\{\Phi(\M{X}^{K})-\Phi(\mathbf{x}^*)+\langle \M{\Lambda}^{K}-\mathbf{w},\mathcal{K}_1(\mathbf{w})\rangle\}\nonumber\\
&=\sup_{\M{y}\in \mathcal{B}_{\rho}}\{\Phi(\M{X}^{K})-\Phi(\mathbf{x}^*)+\langle\M{\Lambda}^{K}-\mathbf{w},\mathcal{K}_1({\mathbf{W}}^{K})\rangle\}\nonumber\\
&=\sup_{\M{y}\in \mathcal{B}_{\rho}}\{\Phi(\M{X}^{K})-\Phi(\mathbf{x}^*)+\langle \M{X}^K-\M{x}^*,-\M{B}\M{Y}^K \rangle\nonumber\\
&\quad\quad\quad\quad+\langle \M{Y}^k-\M{y},\M{B}\M{X}^K \rangle\}\nonumber\\
&=\sup_{\M{y}\in \mathcal{B}_{\rho}}\{\Phi(\M{X}^{K})-\Phi(\mathbf{x}^*)-\langle\M{y},\M{B}\M{X}^K\rangle\}\nonumber\\
&=\Phi(\M{X}^{K})-\Phi(\mathbf{x}^*)+\rho\|\M{BX}^K\|.
\end{align}
Note that $\M{y}^0=0$. Together with \eqref{PDGAP}, one has
\begin{align*}
\Phi(\M{X}^K)-\Phi(\M{x}^*)+\rho\|\M{B}\M{X}^K \|
\leq \frac{\|\M{x}^0-\M{x}^*\|_{\Gamma^{-1}}^2+\rho^2\|\M{Q}\|}{2K}.
\end{align*}
Letting $\rho=\max\{1+\|\M{y}^*\|,2\|\M{y}^*\|\}$, and applying \cite[Lemma 2.3]{ADMM-D1}, it holds that
\begin{align*}
\|\M{B}\M{X}^K\|&\leq \frac{\|\M{x}^0-\M{x}^*\|_{\Gamma^{-1}}^2+\rho^2\|\M{Q}\|}{2K}=O(1/K),\\
|\Phi(\M{X}^K)-\Phi(\M{x}^*)|&\leq \frac{\|\M{x}^0-\M{x}^*\|_{\Gamma^{-1}}^2+\rho^2\|\M{Q}\|}{2K}=O(1/K).
\end{align*}
Since $\|\M{BX}^K\|=\|((I-W)\otimes I_m)^{\frac{1}{2}}\m{X}_1^K\|+\|\M{U}\m{X}_1^K-\m{X}^K_2\|$, the consensus violation \eqref{consensusviolation} holds.
\end{IEEEproof}

\section{Proof of Theorem \ref{THE-3}}\label{APPD}
\begin{IEEEproof}
Since $\mathrm{dist}^2(0,\mathcal{T}(\M{v}^{k+1}))\leq\kappa^2_1\|\M{v}^{k+1}-\M{w}^k\|^2_{\M{M}}$ and $\mathcal{T}$ is metrically subregular at $(\M{w}^{\infty},0)$ and $\mathcal{T}^{-1}(0)=\mathcal{M}^*$, there exists $\kappa_2>0$ and $\epsilon>0$ such that when $\M{v}^{k+1}\in\mathcal{B}_{\epsilon}(\M{w}^{\infty})$
\begin{align*}
&\mathrm{dist}(\M{v}^{k+1},\mathcal{M}^*)=\mathrm{dist}(\M{v}^{k+1},\mathcal{T}^{-1}(0))\\
&\leq \kappa_2 \mathrm{dist}(0,\mathcal{T}(\M{v}^{k+1}))\leq \kappa_1\kappa_2\|\M{v}^{k+1}-\M{w}^k\|_{\M{M}}.
\end{align*}
By \eqref{Rta}, it holds that $\lim_{k\rightarrow\infty}\M{v}^k=\lim_{k\rightarrow\infty}\M{w}^k=\M{w}^{\infty}$. Hence, there exists ${K}\geq0$ such that $\M{v}^{k+1}\in\mathcal{B}_{\epsilon}(\M{w}^{\infty})$. It implies that
\begin{align}\label{TH3-Q1}
\mathrm{dist}(\M{v}^{k+1},\mathcal{M}^*)\leq \kappa_1\kappa_2\|\M{v}^{k+1}-\M{w}^k\|_{\M{M}}, \forall k\geq {K}.
\end{align}
Next, from \eqref{POSI}, we have
\begin{align}\label{TH3-Q2}
&\mathrm{dist}(\M{v}^{k+1},\mathcal{M}^*)\geq \frac{1}{c_2} \mathrm{dist}_{\M{M}}(\M{v}^{k+1},\mathcal{M}^*)\nonumber\\
&\quad\geq\frac{1}{c_2}(\mathrm{dist}_{\M{M}}(\M{w}^{k},\mathcal{M}^*)-\|\M{v}^{k+1}-\M{w}^{k}\|_{\M{M}}).
\end{align}
Combining \eqref{TH3-Q1} and \eqref{TH3-Q2}, it holds that
\begin{align}\label{TH3-Q3}
&\mathrm{dist}_{\M{H}}(\M{w}^k,\mathcal{M}^*)\leq\frac{c_2}{c_1} \mathrm{dist}_{\M{M}}(\M{w}^k,\mathcal{M}^*)\nonumber\\
&\quad \leq\frac{c_2(\kappa_1\kappa_2c_2+1)}{c_1}\|\M{w}^k-\M{v}^{k+1}\|_{\M{M}}, \forall k\geq {K}.
\end{align}
Note that for any $\M{w}^*\in\mathcal{M}^*$, from \textbf{VI 2}, we have $\Theta(\M{x},\M{y}^*)-\Theta(\M{x}^*,\M{y})>0,\forall \M{w}\in \mathcal{M}$. Thus, by \eqref{BASCI-IEQ1}, it holds that
$
\|\M{w}^k-\M{w}^{\infty}\|_{\M{H}}^2-\|\M{w}^{k+1}-\M{w}^{\infty}\|_{\M{H}}^2\geq \|\M{w}^k-\M{v}^{k+1}\|^2_{\M{M}}.
$
Together with \eqref{TH3-Q3}, we have
\begin{align*}
&\mathrm{dist}_{\M{H}}^2(\M{w}^k,\mathcal{M}^*)-\mathrm{dist}_{\M{H}}^2(\M{w}^{k+1},\mathcal{M}^*)\geq \|\M{w}^k-\M{v}^{k+1}\|^2_{\M{M}}\\
&\quad \geq \frac{c_1^2}{c^2_2(\kappa_1\kappa_2c_2+1)^2}  \mathrm{dist}^2_{\M{H}}(\M{w}^k,\mathcal{M}^*), \forall k\geq {K},
\end{align*}
which implies that
$$
\mathrm{dist}_{\M{H}}(\M{w}^{k+1},\mathcal{M}^*) \leq  \varrho \mathrm{dist}_{\M{H}}(\M{w}^k,\mathcal{M}^*),\forall k\geq {K},
$$
where $\varrho=\sqrt{1-\frac{c_1^2}{c^2_2(\kappa_1\kappa_2c_2+1)^2}}<1$.

Finally, we prove that the sequence $\{\M{w}^k\}$ converges to $\M{w}^{\infty}$ R-linearly. Since the sequence $\{\M{w}^k\}$ is a Fej\'{e}r monotone sequence with respect to $\mathcal{M}^*$ in $\mathbf{H}$-norm, i.e.,
$\|\M{w}^k-\M{w}^*\|_{\M{H}}^2\geq\|\M{w}^{k+1}-\M{w}^*\|_{\M{H}}^2,\forall \M{w}^*\in\mathcal{M}^*,$
we have that, when $k>{K}$
\begin{align*}
&\|\M{w}^{k}-\M{w}^{k+1}\|_{\M{H}}\\
&\leq \|\M{w}^{k}-\mathcal{P}_{\mathcal{M}^*}^{\M{H}}(\M{w}^k)\|_{\M{H}}+\|\M{w}^{k+1}-\mathcal{P}_{\mathcal{M}^*}^{\M{H}}(\M{w}^k)\|_{\M{H}}\\
&\leq 2\|\M{w}^{k}-\mathcal{P}_{\mathcal{M}^*}^{\M{H}}(\M{w}^k)\|_{\M{H}}=2 \mathrm{dist}_{\M{H}}(\M{w}^k,\mathcal{M}^*)\\
&\leq 2 \varrho^{k-{K}} \mathrm{dist}_{\M{H}}(\M{w}^{{K}},\mathcal{M}^*).
\end{align*}
Therefore, it holds that, when $k>{K}$
\begin{align*}
&\|\M{w}^{k}-\M{w}^{\infty}\|_{\M{H}}=\|\sum_{s=k}^{\infty}(\M{w}^s-\M{w}^{s+1})\|_{\M{H}}\\
&\leq\sum_{s=k}^{\infty}\|\M{w}^{s}-\M{w}^{s+1}\|_{\M{H}}\leq \sum_{s=k}^{\infty}2 \varrho^{s-{K}} \mathrm{dist}_{\M{H}}(\M{w}^{{K}},\mathcal{M}^*)\\
&=\frac{2\varrho^{k-{K}}}{1-\varrho} \mathrm{dist}_{\M{H}}(\M{w}^{{K}},\mathcal{M}^*).
\end{align*}
This completes the proof.
\end{IEEEproof}

\section{Proof of Lemma \ref{InLem}}\label{APPIn1}
\begin{IEEEproof}
From \eqref{Inexact1}, it holds that
\begin{align}\label{Proof-Lemma3-1}
&G(\M{x})-G(\tilde{\M{x}}^{k+1})+\big\langle \M{x}-\tilde{\M{x}}^{k+1},{\M{B}}\tr \M{y}^k+\nabla F(\M{x}^k) \nonumber \\
&\quad +\Gamma^{-1}(\tilde{\M{x}}^{k+1}-\M{x}^k)-\M{d}^k\big\rangle\geq0, \forall \M{x}\in \mathbb{R}^{m(n+p)}.
\end{align}
Rearranging \eqref{Proof-Lemma3-1}, it holds that
\begin{align}\label{Proof-Lemma3-1-1}
&\big\langle \M{x}-\tilde{\M{x}}^{k+1},\M{B}\tr(\M{y}^k-\M{y})-\M{d}^k+\Gamma^{-1}(\tilde{\M{x}}^{k+1}-\M{x}^k)\big\rangle\nonumber\\
&\geq G(\tilde{\M{x}}^{k+1})-G(\M{x})+\big\langle \tilde{\M{x}}^{k+1}-\M{x},\M{B}\tr\M{y}+\nabla F(\M{x}^k) \big\rangle.
\end{align}
By \eqref{Inexact2} and similar as \eqref{Proof-Lemma1-4}, we have
\begin{align}\label{Proof-Lemma3-4}
&\big\langle \M{y}-\M{y}^{k+1},-\M{B}(\tilde{\M{x}}^{k+1}-\M{x})+\M{Q}(\M{y}^{k+1}-\M{y}^k) \big\rangle \nonumber\\
&\quad \geq \big\langle \M{y}^{k+1}-\M{y},-{\M{B}}\M{x} \big\rangle,\forall \M{y}\in \mathbb{R}^{m(n+p)}.
\end{align}
Summing \eqref{Proof-Lemma3-1-1} and \eqref{Proof-Lemma3-4}, we have
\begin{align}
&\big\langle \M{x}-\tilde{\M{x}}^{k+1},\M{B}\tr(\M{y}^k-\M{y})+\Gamma^{-1}(\tilde{\M{x}}^{k+1}-\M{x}^k)\big\rangle\nonumber\\
&\quad+\big\langle \M{y}-\M{y}^{k+1},-\M{B}(\tilde{\M{x}}^{k+1}-\M{x})+\M{Q}(\M{y}^{k+1}-\M{y}^k) \big\rangle\nonumber\\
&\geq G(\tilde{\M{x}}^{k+1})-G(\M{x})+\big\langle \M{v}^{k+1}-\M{w}, \mathcal{K}_1(\M{w})\big\rangle\nonumber\\
&\quad+ \big\langle \tilde{\M{x}}^{k+1}-\M{x},\nabla F(\M{x}^k) -\M{d}^k\big\rangle.
\end{align}
Since $
\langle \M{x}-\tilde{\M{x}}^{k+1},\M{B}\tr(\M{y}^k-\M{y})\rangle+\big\langle \M{y}-\M{y}^{k+1},-\M{B}(\tilde{\M{x}}^{k+1}-\M{x})\rangle
=\langle \M{x}-\tilde{\M{x}}^{k+1},\M{B}\tr(\M{y}^k-\M{y}^{k+1}) \rangle$, it holds that
\begin{align}\label{PPLEMMA3-1-1}
&\big\langle \M{x}-\tilde{\M{x}}^{k+1},\Gamma^{-1}(\tilde{\M{x}}^{k+1}-\M{x}^k)\big\rangle+\big\langle \M{y}-\M{y}^{k+1},\M{Q}(\M{y}^{k+1}-\M{y}^k) \big\rangle\nonumber\\
&\quad+\langle \M{x}-\tilde{\M{x}}^{k+1},\M{B}\tr(\M{y}^k-\M{y}^{k+1}) \rangle\nonumber\\
&\geq G(\tilde{\M{x}}^{k+1})-G(\M{x})+\big\langle \M{v}^{k+1}-\M{w}, \mathcal{K}_1(\M{w})\big\rangle\nonumber\\
&\quad+ \big\langle \tilde{\M{x}}^{k+1}-\M{x},\nabla F(\M{x}^k) -\M{d}^k\big\rangle.
\end{align}
Note that
\begin{align*}
&\|\tilde{\M{x}}^{k+1}-\M{x}\|_{\Gamma^{-1}}^2-2\langle \M{x}-\tilde{\M{x}}^{k+1},\M{B}\tr(\M{y}^k-\M{y}^{k+1}) \rangle\\
&=\|\tilde{\M{x}}^{k+1}-\M{x}-\Gamma\mathbf{B}\tr(\M{y}^{k+1}-\M{y}^k)\|_{\Gamma^{-1}}^2-\|\M{y}^k-\M{y}^{k+1}\|_{\M{B}\Gamma\M{B}\tr}^2\\
&=\|\M{x}^{k+1}-\M{x}\|_{\Gamma^{-1}}^2-\|\M{y}^k-\M{y}^{k+1}\|_{\M{B}\Gamma\M{B}\tr}^2,
\end{align*}
where the last equality follows from \eqref{Inexact3}. Combining the above equality, \eqref{ZYBDS}, and \eqref{PPLEMMA3-1-1}, we have
\begin{align}\label{InLemeq}
&\frac{1}{2}\|\M{w}^k-\M{w}\|_{\mathbf{H}}^2-\frac{1}{2}\|\M{w}^{k+1}-\M{w}\|_{\mathbf{H}}^2 \geq G(\tilde{\M{x}}^{k+1})-G(\M{x}) \nonumber\\
&+ \frac{1}{2}\|\M{x}^k-\tilde{\M{x}}^{k+1}\|_{\Gamma^{-1}}^2+ \frac{1}{2}\|\M{y}^k-\M{y}^{k+1}\|^2_{\M{Q}- \M{B}\Gamma\M{B}\tr}\nonumber\\
& +\big\langle \M{v}^{k+1}-\M{w}, \mathcal{K}_1(\M{w})\big\rangle+ \big\langle \tilde{\M{x}}^{k+1}-\M{x},\nabla F(\M{x}^k) -\M{d}^k\big\rangle.
\end{align}
Similar as Lemma \ref{F1}, by \eqref{SM1} and \eqref{SM2}, the inequalities \eqref{BASCI-IEQ66} and \eqref{INX66} can be proven.
\end{IEEEproof}

\section{Proof of Theorem \ref{THE-4}}\label{APPG}
\begin{IEEEproof}
To show the convergence of V-DISA \eqref{Inexact}, we introduce the following notations.
\begin{align*}
\tilde{\M{x}}_{\varepsilon^k=0}^{k+1}&=\mathrm{prox}^{\Gamma^{-1}}(\M{x}^k-\Gamma\nabla F(\M{x}^k)-\Gamma {\M{B}}\tr \M{y}^k),\\
\overline{\M{y}}^{k+1}&=\M{y}^{k}+\M{Q}^{-1}\M{B}\tilde{\M{x}}_{\varepsilon^k=0}^{k+1},\\
\overline{\M{x}}^{k+1}&=\tilde{\M{x}}_{\varepsilon^k=0}^{k+1}+\Gamma \M{B}\tr(\M{y}^k-\overline{\M{y}}^{k+1}).
\end{align*}
Let $\overline{\M{w}}^k=[(\overline{\M{x}}^k)\tr,(\overline{\M{y}}^k)\tr]\tr$, $\overline{\M{v}}^k=[(\tilde{\M{x}}_{\varepsilon^k=0}^k)\tr,(\overline{\M{y}}^k)\tr]\tr$. It follows from the nonexpansivity of $\mathrm{prox}^{\Gamma^{-1}}_{ G}(\cdot)$ that
\begin{align}\label{PP2}
\|\tilde{\M{x}}^{k+1}-\tilde{\M{x}}_{\varepsilon^k=0}^{k+1}\|\leq \tau\| \M{d}^k\|\leq\tau \sqrt{m}\varepsilon^k.
\end{align}
Note that
\begin{align*}
\|\M{y}^{k+1}-\overline{\M{y}}^{k+1}\|&=\|\M{Q}^{-1}\M{B}(\tilde{\M{x}}^{k+1}-\tilde{\M{x}}_{\varepsilon^k=0}^{k+1})\|\\
&\leq\|\M{Q}^{-1}\M{B}\| \|\tilde{\M{x}}^{k+1}-\tilde{\M{x}}_{\varepsilon^k=0}^{k+1}\|,\\
\|\M{x}^{k+1}-\overline{\M{x}}^{k+1}\|&=\|(\tilde{\M{x}}^{k+1}-\tilde{\M{x}}_{\varepsilon^k=0}^{k+1})+\Gamma\M{B}\tr(\overline{\M{y}}^{k+1}-\M{y}^{k+1})\|\\
&\leq(1+\tau\|\M{B}\tr\|\|\M{Q}^{-1}\M{B}\|)\|\tilde{\M{x}}^{k+1}-\tilde{\M{x}}_{\varepsilon^k=0}^{k+1}\|.
\end{align*}
By \eqref{PP2}, we have
\begin{align}\label{EQ:APPH-Proof-2}
&\|\M{w}^{k+1}-\overline{\M{w}}^{k+1}\|_{\M{H}}\nonumber\\
&\leq  \underbrace{\tau\sqrt{2m}\|\M{H}\|^{\frac{1}{2}}(\|\M{Q}^{-1}\M{B}\|(1+\tau\|\M{B}\|)+1)}_{:=\psi>0} \varepsilon^k.
\end{align}
It follows from \eqref{BASCI-IEQ66} that
\begin{align}\label{PRFG2}
&\Theta(\tilde{\M{x}}_{\varepsilon^k=0}^{k+1},\M{y})-\Theta(\M{x},\overline{\M{y}}^{k+1})+\frac{1}{2}\|\M{w}^k-\overline{\M{v}}^{k+1}\|^2_{\M{M}}\nonumber\\
&\leq \frac{1}{2}\|\M{w}^k-\M{w}\|_{\M{H}}^2-\frac{1}{2}\|\overline{\M{w}}^{k+1}-\M{w}\|_{\M{H}}^2.
\end{align}
Let $\M{w}=\M{w}^*\in \mathcal{M}^*$ in \eqref{PRFG2}. By \textbf{VI 2}, it deduces that
\begin{align*}
\|\M{w}^k-\M{w}^*\|_{\M{H}}^2\geq \|\overline{\M{w}}^{k+1}-\M{w}^*\|_{\M{H}}^2+\|\M{w}^k-\overline{\M{v}}^{k+1}\|^2_{\M{M}}.
\end{align*}
Then, for $\forall k\geq0$, $\|\M{w}^k-\M{w}^*\|_{\M{H}}^2\geq \|\overline{\M{w}}^{k+1}-\M{w}^*\|_{\M{H}}^2$. Hence, it holds that for $\forall \M{w}^*\in \mathcal{M}^*$
\begin{align}\label{Gproof1}
\|\M{w}^{k+1}-\M{w}^*\|_{\M{H}}&\leq \|\overline{\M{w}}^{k+1}-\M{w}^*\|_{\M{H}}+\|\M{w}^{k+1}-\overline{\M{w}}^{k+1}\|_{\M{H}}\nonumber\\
&\leq \|\M{w}^k-\M{w}^*\|_{\M{H}} + \|\M{w}^{k+1}-\overline{\M{w}}^{k+1}\|_{\M{H}}\nonumber\\
&\leq \|\M{w}^k-\M{w}^*\|_{\M{H}} +\psi \varepsilon^k.
\end{align}
Summing the above inequality over $k=0,1,\cdots,K-1$, one obtains that for any $K\geq 1$,
\begin{align*}
&\sum_{k=0}^{K-1}\big(\|\M{w}^{k+1}-\M{w}^*\|_{\M{H}}-\|\M{w}^k-\M{w}^*\|_{\M{H}}\big)\leq\sum_{k=0}^{K-1}\psi \varepsilon^k\\
&\Rightarrow \|\M{w}^{K}-\M{w}^*\|_{\M{H}} \leq \|\M{w}^{0}-\M{w}^*\|_{\M{H}} + \sum_{k=0}^{K-1}\psi \varepsilon^k.
\end{align*}
Since $\M{H}\succ 0$ and $\{\varepsilon^k\}$ is summable, for any $\M{w}^*\in \mathcal{M}^*$, $\{\|\M{w}^k-\M{w}^*\|\}$ and $\{\|\M{v}^k-\M{w}^*\|\}$ are bounded.

Let $\M{w}=\M{w}^*\in \mathcal{M}^*$ in \eqref{BASCI-IEQ66}. One gets that
\begin{align}\label{TH4RATA}
\|\M{w}^k-\M{v}^{k+1}\|^2_{\M{M}}\leq &\|\M{w}^k-\M{w}^*\|_{\M{H}}^2-\|\M{w}^{k+1}-\M{w}^*\|_{\M{H}}^2\nonumber\\
&+2\langle \tilde{\M{x}}^{k+1}-\M{x}^*,\M{d}^k\rangle.
\end{align}
Summing the inequality over $k=0,1,\cdots,\infty$, it gives that
\begin{align*}
&\sum_{k=0}^{\infty}\|\M{w}^k-\M{v}^{k+1}\|^2_{\M{M}}\leq\sum_{k=0}^{\infty}(\|\M{w}^k-\M{w}^*\|_{\M{H}}^2-\|\M{w}^{k+1}-\M{w}^*\|_{\M{H}}^2)\\
&\quad\quad\quad\quad\quad\quad\quad\quad\quad\quad+\sum_{k=0}^{\infty}2\langle \tilde{\M{x}}^{k+1}-\M{x}^*,\M{d}^k\rangle\\
&\leq \|\M{w}^{0}-\M{w}^*\|_{\M{H}}^2 + \sum_{k=0}^{\infty}2\|\tilde{\M{x}}^{k+1}-\M{x}^*\|\varepsilon^k<\infty.
\end{align*}
Thus, it holds that $\|\M{w}^k-\M{v}^{k+1}\|\rightarrow 0,k\rightarrow\infty$. Let $\mathbf{w}^{\infty}$ be an accumulation point of $\{\mathbf{w}^k\}$ and $\{\mathbf{w}^{k_j}\}$ be a subsequence converging to $\mathbf{w}^{\infty}$. Similar as Theorem \ref{THE-1}, we can show that $\mathbf{w}^{\infty}\in \mathcal{M}^*$. Therefore, by \eqref{Gproof1}, it gives that
$
\|\M{w}^{k+1}-\M{w}^{\infty}\|_{\M{H}}\leq \|\M{w}^k-\M{w}^{\infty}\|_{\M{H}} +\psi \varepsilon^k.
$
Since $\sum_{k=0}^{\infty}\varepsilon^k<\infty$, by \cite[Lemma 3.2]{InexactProximalALM2} the quasi-Fej\'{e}r monotone sequence $\{\|\M{w}^k-\M{w}^{\infty}\|_{\mathbf{H}}\}$ converges to a unique limit point. Then, with $\M{w}^{\infty}$ being an accumulation point of $\{\M{w}^k\}$, it holds that $\|\M{w}^k-\M{w}^{\infty}\|^2_{\M{H}}\rightarrow0$, i.e., $\lim_{k\rightarrow \infty}\M{w}^k=\M{w}^{\infty}$.
\end{IEEEproof}

\section{Proof of Theorem \ref{THE-5}}\label{APPGG}
\begin{IEEEproof}
It follows from \eqref{Gproof1} that, for $\forall \M{w}^*\in\mathcal{M}^*$,
\begin{align}\label{EQ:APPI-Proff-1}
\|\M{w}^{K}-\M{w}^*\|_{\M{H}}\leq \|\M{w}^{0}-\M{w}^*\|_{\M{H}} +\psi\sum_{k=0}^{K-1}\varepsilon^k, \forall K\geq 1.
\end{align}
Since $\|\M{w}^k-\M{w}^*\|_{\M{H}}^2\geq \|\overline{\M{w}}^{k+1}-\M{w}^*\|_{\M{H}}^2, \forall k\geq0$,
combining it with \eqref{EQ:APPH-Proof-2} and \eqref{EQ:APPI-Proff-1}, we have
\begin{align*}
&\|\tilde{\M{x}}^{k+1}-\M{x}^*\|\leq \|\M{v}^{k+1}-\M{w}^*\|\\
&\leq\|\M{v}^{k+1}-\overline{\M{w}}^{k+1}\|+\sqrt{\|\M{H}^{-1}\|}\|\overline{\M{w}}^{k+1}-\M{w}^*\|_{\M{H}}\\
&\leq \|\M{v}^{k+1}-\overline{\M{w}}^{k+1}\|+\sqrt{\|\M{H}^{-1}\|}(\|\M{w}^{0}-\M{w}^*\|_{\M{H}} +\psi\sum_{j=0}^{k-1}\varepsilon^j),
\end{align*}
Note that $\|\M{v}^{k+1}-\overline{\M{w}}^{k+1}\|\leq\tau\sqrt{2m}(1+\|\M{Q}^{-1}\M{B}\|)\varepsilon^k$. Let $\sigma=\max\{\sqrt{\|\M{H}^{-1}\|}\psi,\tau\sqrt{2m}(1+\|\M{Q}^{-1}\M{B}\|)\}$.
We have
$
\|\tilde{\M{x}}^{k+1}-\M{x}^*\|\leq\sigma\|\M{w}^{0}-\M{w}^*\|_{\M{H}} +\sigma\psi\sum_{k=0}^{\infty}\varepsilon^k
$.
Thus, for any $K\geq 1$
\begin{align*}
&\sum_{k=0}^{K}\langle \tilde{\M{x}}^{k+1}-\M{x}^*,\M{d}^k\rangle\leq \sum_{k=0}^{K}\|\tilde{\M{x}}^{k+1}-\M{x}^*\|\cdot \|\M{d}^k\|\\
&\leq \tau\sigma\sqrt{2m}\big(\|\M{w}^{0}-\M{w}^*\|_{\M{H}} +\psi\sum_{k=0}^{\infty}\varepsilon^k\big) \big(\sum_{k=0}^{\infty} \varepsilon^k\big).
\end{align*}
Summing the inequality \eqref{TH4RATA} over $k=0,\cdots,K$, we have
\begin{align}\label{INRATA}
\sum_{k=0}^{K}\|\M{w}^k-\M{v}^{k+1}\|^2_{\M{M}}\leq \psi_1<\infty,
\end{align}
where $
\psi_1=\|\M{w}^{0}-\M{w}^*\|_{\M{H}}^2+2 \tau\sigma\sqrt{2m}\big(\|\M{w}^{0}-\M{w}^*\|_{\M{H}} +\psi\sum_{k=0}^{\infty}\varepsilon^k\big) \big(\sum_{k=0}^{\infty} \varepsilon^k\big)
$.
By the definition of $\tilde{\M{x}}^{k+1}$ and the update of $\M{y}^{k+1}$, it gives
\begin{align*}
&0\in \partial G(\tilde{\M{x}}^{k+1})+\nabla F(\M{x}^k)+{\M{B}}\tr\M{y}^k+\frac{1}{\tau}(\tilde{\M{x}}^{k+1}-\M{x}^k)-\M{d}^k,\\
&0=-{\M{B}}\tilde{\M{x}}^{k+1}+\M{Q}(\M{y}^{k+1}-\M{y}^{k}).
\end{align*}
Similar as the proof of Theorem \ref{THE-2}, we can prove that
\begin{align*}
\mathrm{dist}^2(0,\mathcal{T}(\M{v}^{k+1}))\leq\kappa^2_3\|\M{v}^{k+1}-\M{w}^k\|^2_{\M{M}}+4(\varepsilon^k)^2,
\end{align*}
where $\kappa^2_3=\frac{1}{c_1^2}\max\{4L^2+\frac{4}{\tau^2},4\|\M{B}\|^2+\|\M{Q}\|^2\}$. Thus, we can get from \eqref{INRATA} that \eqref{INRATA2} hold.

Summing \eqref{INX66} over $k=0,1,\cdots,K-1$, we obtain
\begin{align*}
&2\sum_{k=0}^{K-1}(\mathcal{L}(\tilde{\M{x}}^{k+1},\M{y})-\mathcal{L}(\M{x},\M{y}^{k+1}))\\
&\leq \|\M{w}^0-\M{w}\|_{\M{H}}^2 + 2\sum_{k=0}^{K-1}\|\tilde{\M{x}}^{k+1}-\M{x}\| \|\M{d}^k\|\\
&\leq \|\M{w}^0-\M{w}\|_{\M{H}}^2+ 2\tau\sigma\bar{\varepsilon}\sqrt{2m}\big(\|\M{w}^{0}-\M{w}\|_{\M{H}} +\psi\bar{\varepsilon}\big)\\
&= (1+2\tau\sigma\bar{\varepsilon}\sqrt{2m})\|\M{w}^0-\M{w}\|_{\M{H}}^2+2\tau\sigma\psi\bar{\varepsilon}^2\sqrt{2m}.
\end{align*}
It follows from the convexity of $F$, $G$ and the definition of $(\M{X}^{K},\M{Y}^K)$ that $K(\mathcal{L}(\M{X}^{K},\M{y})-\mathcal{L}(\M{x},\M{Y}^{K}))\leq\sum_{k=0}^{K-1}(\mathcal{L}(\tilde{\M{x}}^{k+1},\M{y})-\mathcal{L}(\M{x},\M{y}^{k+1}))$. Let $\omega_1=2\tau\sigma\bar{\varepsilon}\sqrt{2m}$ and $\omega_2=\|\M{x}^0-\M{x}^*\|_{\Gamma^{-1}}^2+\rho^2\|\M{Q}\|$. Similar as the proof of Theorem \ref{THE-PDGAP}, we have
$$
\Phi(\M{X}^K)-\Phi(\M{x}^*)+\rho\|\M{B}\M{X}^K \|\leq \frac{(1+\omega_1)\omega_2+\omega_1\bar{\varepsilon}}{2K}.
$$
Letting $\rho=\max\{1+\|\M{y}^*\|,2\|\M{y}^*\|\}$, applying \cite[Lemma 2.3]{ADMM-D1}, we can get \eqref{TH6-28} and \eqref{TH6-29}.
\end{IEEEproof}

\vspace{-15 mm}
\begin{IEEEbiography}[{\includegraphics[width=1.3in,height=1.2in,clip,keepaspectratio]{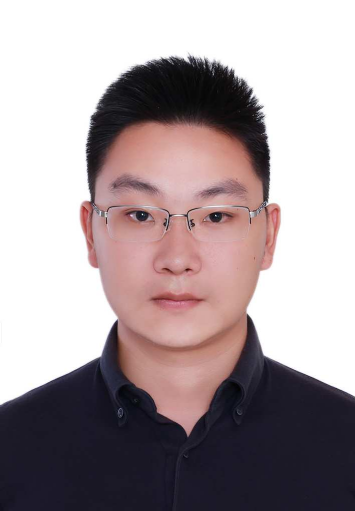}}]{Luyao Guo}
received the B.S. degree in information and computing science
from Shanxi University, Taiyuan, China, in 2020. He is currently pursuing the Ph.D.
degree in applied mathematics with the Jiangsu Provincial Key Laboratory of Networked Collective
Intelligence, School of Mathematics, Southeast University, Nanjing, China. His current research focuses on distributed optimization and learning.
\end{IEEEbiography}
\vspace{-15 mm}
\begin{IEEEbiography}[{\includegraphics[width=1.3in,height=1.2in,clip,keepaspectratio]{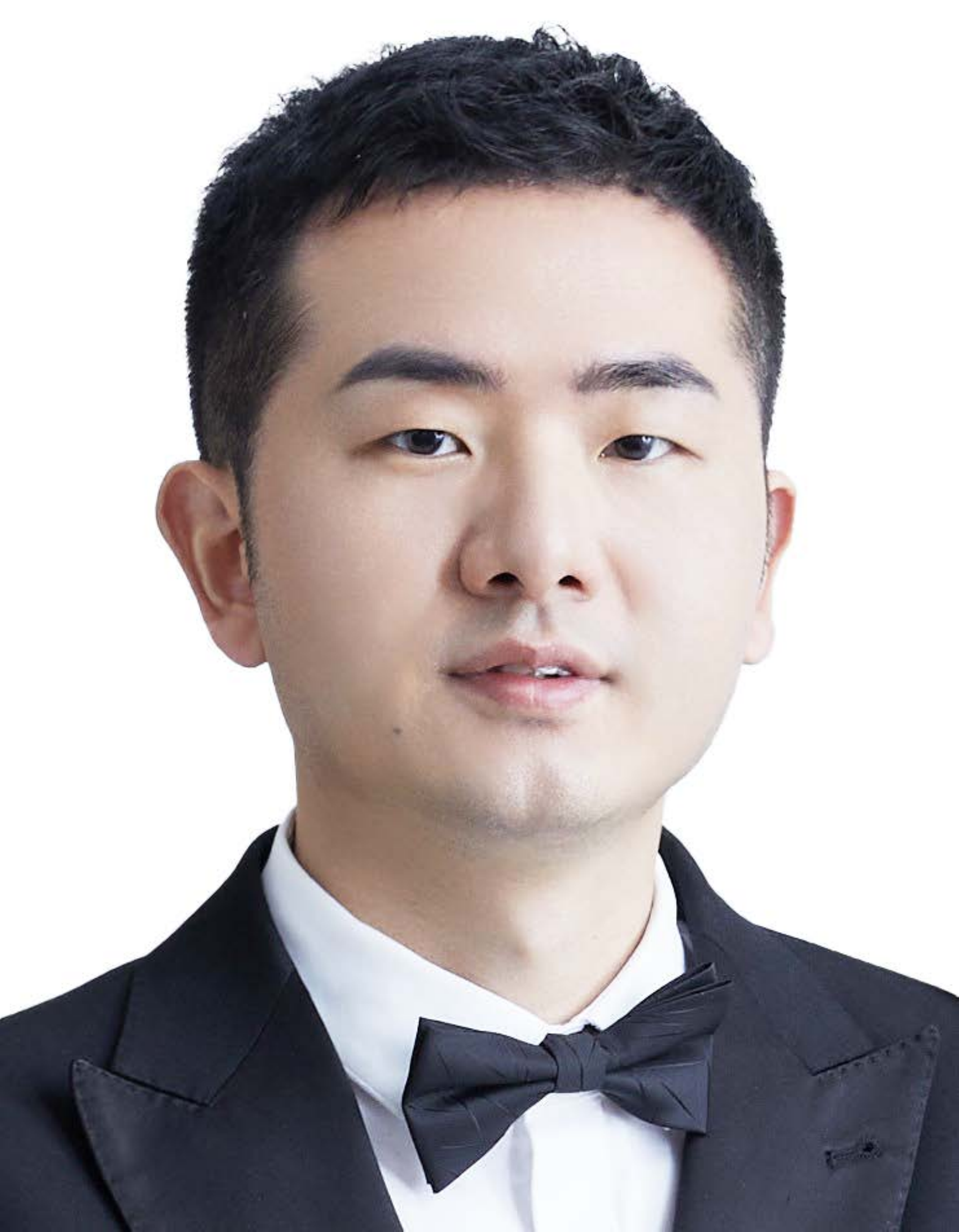}}]{Xinli Shi (Senior Member, IEEE)}
received the B.S. degree in software engineering, the M.S. degree in applied mathematics, and the Ph.D. degree in control science and engineering from Southeast University, Nanjing, China, in 2013, 2016, and 2019, respectively. He held a China Scholarship Council Studentship for one-year study with the University of Royal Melbourne Institute of Technology, Melbourne, VIC, Australia, in 2018. He is currently an Associate Professor with the School of Cyber Science and Engineering, Southeast University. His current research interests include distributed optimization, nonsmooth analysis, and network control systems. Dr. Shi was the recipient of the Outstanding Ph.D. Degree Thesis Award from Jiangsu Province, China.
\end{IEEEbiography}
\vspace{-15 mm}
\begin{IEEEbiography}[{\includegraphics[width=1.3in,height=1.2in,clip,keepaspectratio]{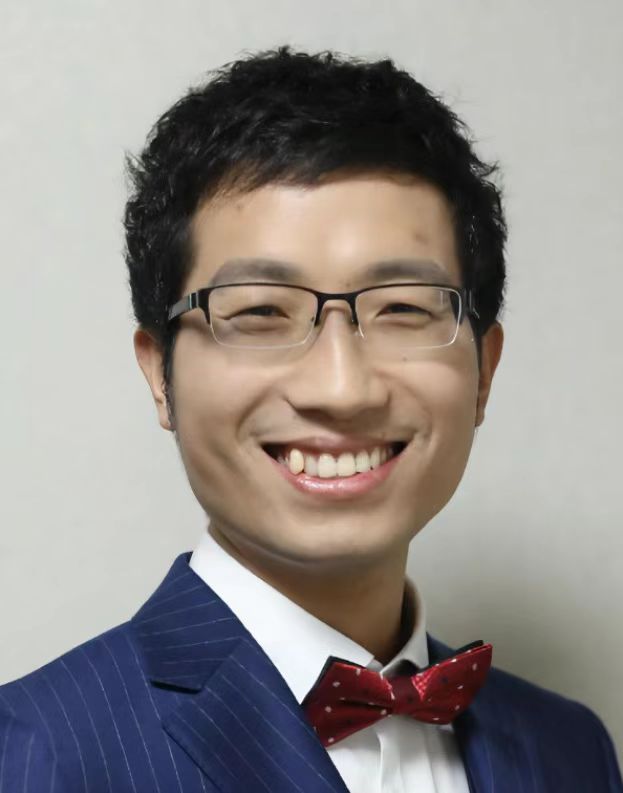}}]{Shaofu Yang (Member, IEEE)}
received the B.S. and M.S. degrees in applied mathematics from the Department of Mathematics, Southeast University, Nanjing, China, in 2010 and 2013, respectively, and the Ph.D. degree in engineering from the Department of Mechanical and Automation Engineering, The Chinese University of Hong Kong, Hong Kong, in 2016. He was a Post-Doctoral Fellow with the City University of Hong Kong, Hong Kong, in 2016. He is currently an Associate Professor with the School of Computer Science and Engineering, Southeast University. His current research interests include distributed optimization and learning, game theory, and their applications.
\end{IEEEbiography}
\vspace{-15 mm}
\begin{IEEEbiography}[{\includegraphics[width=1.3in,height=1.2in,clip,keepaspectratio]{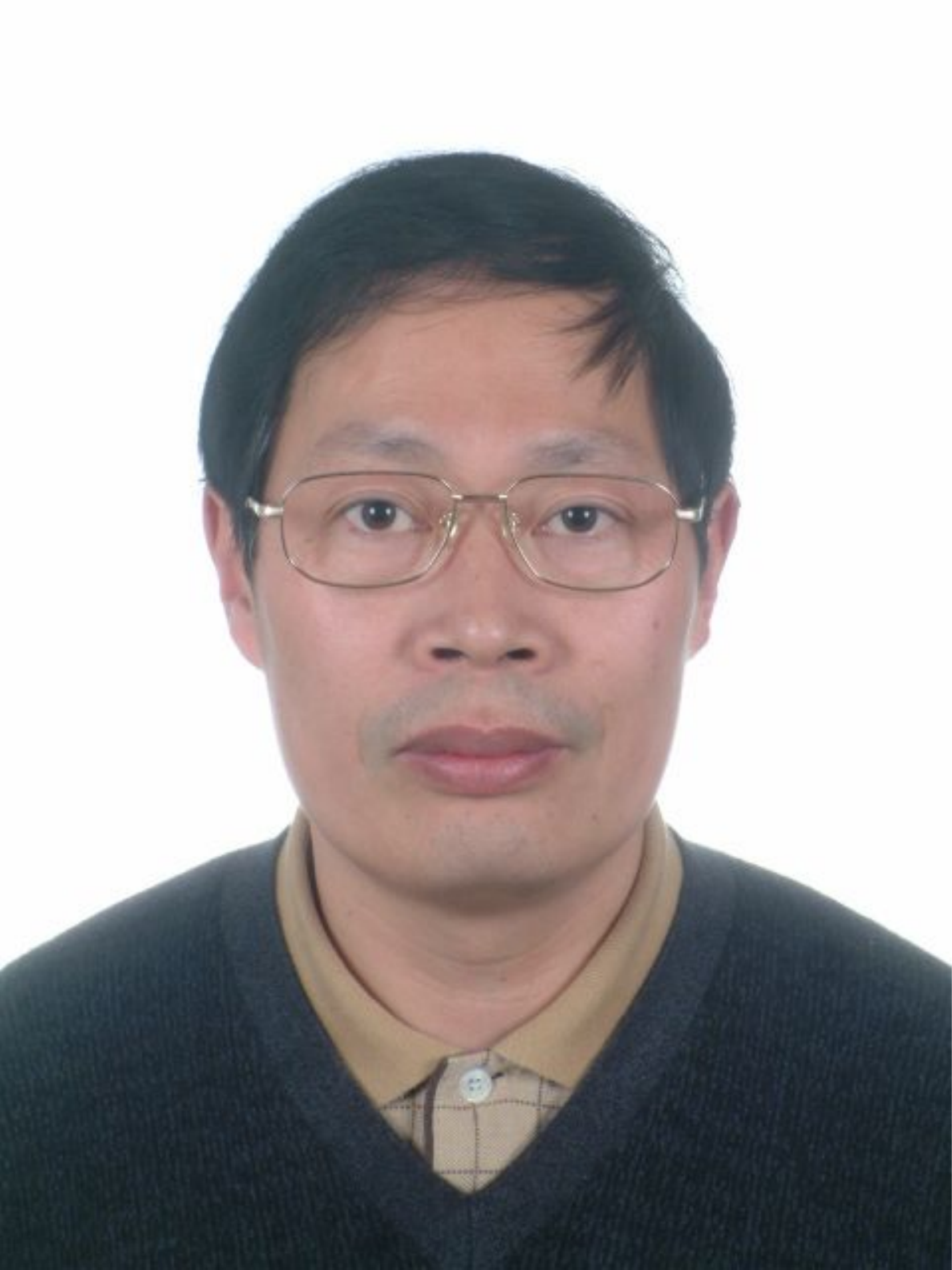}}]{Jinde Cao (Fellow, IEEE)}
received the B.S. degree from Anhui Normal University, Wuhu, China, the M.S. degree from Yunnan University, Kunming, China, and the Ph.D. degree from Sichuan University, Chengdu, China, all in mathematics/applied mathematics, in 1986, 1989, and 1998, respectively. He was a Postdoctoral Research Fellow at the Department of Automation and Computer-Aided Engineering, Chinese University of Hong Kong, Hong Kong, from 2001 to 2002. He is an Endowed Chair Professor, the Dean of the School of Mathematics and the Director of the Research Center for Complex Systems and Network Sciences at Southeast University (SEU). He is also the Director of the National Center for Applied Mathematics at SEU-Jiangsu of China and the Director of the Jiangsu Provincial Key Laboratory of Networked Collective Intelligence of China. Prof. Cao was a recipient of the National Innovation Award of China, Obada Prize and the Highly Cited Researcher Award in Engineering, Computer Science, and Mathematics by Clarivate Analytics. He is elected as a member of Russian Academy of Sciences, a member of the Academy of Europe, a member of Russian Academy of Engineering, a member of the European Academy of Sciences and Arts, a member of the Lithuanian Academy of Sciences, a fellow of African Academy of Sciences, and a fellow of Pakistan Academy of Sciences.
\end{IEEEbiography}

\end{document}